\documentclass[12pt, a4paper]{article}

\usepackage[english]{babel}

\usepackage{authblk}
\usepackage{placeins}

\usepackage{graphicx}
\usepackage{multirow}
\usepackage{bm}
\usepackage{pgfplots,tikz}
\usepackage{pgfplotstable}
\usepackage{hyperref}

\usepackage{mathtools}

\usepackage{amsmath}
\usepackage{amssymb}
\usepackage[latin1]{inputenc}
\usepackage{caption}
\usepackage{subcaption}
\usepackage{graphicx}
\usepackage{graphics}
\usepackage{xcolor}
\usepackage{float}
\usepackage{listings}
\usepackage{textcomp}

\usepackage{epstopdf}
\usepackage{enumitem}
\usepackage{comment}
\usepackage{cite}

\newcommand{\dJac}{\mathrm{\mathbf{D}}}

\usepackage{caption}
\usepackage{color}

\newcommand*\diff{\mathop{}\!\mathrm{d}}

\newcommand{\pdiff}[2]{\frac{\partial #1}{\partial #2}}

\DeclareMathOperator*{\argmin}{argmin}

\newcommand{\MA}{\leavevmode Monge-Amp\`ere}
\newcommand{\MAe}{\leavevmode Monge-Amp\`ere equation}

\usepackage{pgfplots}
\usetikzlibrary{arrows.meta}
\usetikzlibrary{pgfplots.patchplots}
\usetikzlibrary{calc}
\usetikzlibrary{positioning}
\usetikzlibrary{decorations.pathreplacing}
\usetikzlibrary{decorations.markings}
\usetikzlibrary{intersections}

\usepackage{calc}

\usepackage{afterpage} 
\usepackage[toc,page]{appendix}
\usepackage{empheq}
\usepackage[nottoc]{tocbibind}

\newcommand\elll{l}

\providecommand{\keywords}[1]
{
	\small	
	\textbf{\textit{Keywords---}} #1
}

\title{An Iterative Least-Squares Method for the Hyperbolic Monge-Amp\`ere Equation with Transport Boundary Condition}
\date{}
\author[1,*]{M.W.M.C. Bertens}
\author[1]{{}M.J.H. Anthonissen}
\author[1]{\\J.H.M. ten Thije Boonkkamp}
\author[1,2]{W.L. IJzerman}
\affil[1]{CASA, Department of Mathematics and Computer Science, Eindhoven University of Technology, PO Box 513, 5600 MB Eindhoven, The Netherlands}
\affil[2]{Signify Research, High Tech Campus 7, 5656 AE Eindhoven, The Netherlands}
\affil[*]{Corresponding author: m.w.m.c.bertens@tue.nl}
\begin{document}
		
\maketitle
	
\begin{abstract}
A least-squares method for solving the hyperbolic \MAe{} with transport boundary condition is introduced. The method relies on an iterative procedure for the gradient of the solution, the so-called mapping. By formulating error functionals for the interior domain, the boundary, both separately and as linear combination, three minimization problems are solved iteratively to compute the mapping. After convergence, a fourth minimization problem, to compute the solution of the \MAe{}, is solved.
The approach is based on a least-squares method for the elliptic \MAe{}\cite{Prins_2015}, and is improved upon by the addition of analytical solutions for the minimization on the interior domain and by the introduction of two new boundary methods. Lastly, the iterative method is tested on a variety of examples. It is shown that, when the iterative method converges, second-order global convergence as function of the spatial discretization is obtained.
\end{abstract}


\section{Introduction}
 In this paper we introduce a least-squares method for the hyperbolic \MAe{} with transport boundary condition. We are motivated by applications to optical design. In~\cite{CorienThesis} it was found that designing lenses and reflectors for some single-optical-surface systems is equivalent to solving the \MA{} (MA) equation with transport boundary condition. The optical surface satisfies either the elliptic (+) or hyperbolic (-) MA equation, i.e.,
 \begin{subequations}
 \begin{align}
   	\label{eqn:MA1}
    	\det(\dJac^2 u(\mathbf{x}) ) = \pm \frac{E(\mathbf{x})}{I(\nabla u (\mathbf{x}) )},  \qquad \mathbf{x} \in \mathcal{X},    
 \end{align}
 where $\dJac^2 u$ is the Hessian matrix of the optical surface $z = u(\mathbf{x})$, $E \geq 0$ the emittance of the source with domain $\mathcal{X} \subset \mathbb{R}^2$ and $I > 0$ the illuminance on the target with domain $\mathcal{Y} = \nabla u(\mathcal{X}) \subset \mathbb{R}^2$.
 The accompanying transport boundary condition is given by
 \begin{align}
    	\nabla u (\partial \mathcal{X}) = \partial \mathcal{Y}.
 \end{align}
  \end{subequations}
In the elliptic case this boundary condition follows from convexity or concavity arguments of the optical surface~\cite{CorienThesis}. For the hyperbolic MA equation, the transport boundary condition cannot be derived in the same manner as $z = u(\mathbf{x})$ is a saddle surface, and one instead relies on optical arguments, e.g., the edge-ray principle~\cite{Ries1994}. Further disparities between the elliptic and hyperbolic variant are common. This for one is due to the connection of the elliptic variant to the rich field of optimal transport (OT)~\cite[p.~282]{villani_optimal_2009}, which is absent for the hyperbolic MA equation. OT was originally established by Monge, who was concerned with rearranging mass from one distribution to another\cite[p.~xiv]{Santambrogio_2015}. Brenier proved that, assuming regularity conditions~\cite{Brenier_1987}, the corresponding optimal (point-wise) transport map $\mathbf{m} = \nabla u$ satisfies the elliptic \MAe{}~\eqref{eqn:MA1}, where the plus sign is assumed and $E$ and $I$ should be interpreted as densities. See~\cite[p.~323-332]{villani_optimal_2009}, for example, for regularity, uniqueness and existence of solutions to the elliptic MA equation. Assumptions on energy conservation allow design of optical systems to be cast in the framework of OT, but only for the elliptic variant and not for the hyperbolic MA equation. As a consequence, results regarding regularity, uniqueness and existence for the hyperbolic MA equation are scarce.
 The most relevant results for the hyperbolic MA equation
 	 \begin{align}
  	 	\label{eqn:intro:HMA}
    	\det(\dJac^2 u(\mathbf{x}) ) = - f^2(\mathbf{x}), \qquad \mathbf{x} \in \mathcal{X} \subset \mathbb{R}^2,
 \end{align}
 with $f:\mathbb{R}^2 \mapsto (0, \infty)$, follow from the method of characteristics~\cite{Tunitskii_1997, GlobalSmoothSolutions} and only hold for Cauchy boundary conditions. Equation~\eqref{eqn:intro:HMA} has received little interest in numerical and computational journals. To the best of our knowledge, the equation has been solved twice, once on a triangular computational domain~\cite[p.~614]{Westcott_1976} and secondly on a rectangular computational domain by Bertens et al.~\cite{Bertens2022}. The former method, based on a finite difference scheme, assumes Cauchy boundary conditions on an initial curve, and does not treat boundary conditions on the rest of the domain. The latter method, derived using the method of characteristics, shows that the requirements on the boundary conditions are strenuous. 
 The method of characteristics shows, by parametrizing the characteristics with the $x$-coordinate and assuming Cauchy boundary conditions on an initial curve $x =$ const, that the remaining boundary conditions depend on the location of the characteristics. Consequently, the required boundary conditions are fundamentally different from the transport boundary condition.

 In this work, we therefore do not consider the method of characteristics, but instead resort to a least-squares method which has been proven to work for various elliptic problems, among which the \MAe{}~\cite{Prins_2015}, the generalized \MAe{} and the generated Jacobian equation~\cite{Romijn_2021}. The least-squares algorithm is an iterative method which does not directly solve for the unknown $u$, but instead first constructs the mapping $\mathbf{m} = \nabla u$ and afterwards approximates $u$. The general outline of the least-squares method is as follows: first, we approximate the Jacobi matrix of $\mathbf{m}$ in the interior of the domain by minimizing an error functional. Secondly, $\mathbf{m}$ restricted to the boundary of the domain is approximated. By minimizing another error functional involving the newly found Jacobi matrix and the boundary approximation, we obtain a new approximation for the mapping. We proceed by repeating these three steps iteratively until $\mathbf{m}$ no longer changes and subsequently calculate $u$ by minimizing a fourth functional. One of the benefits of this method is that each of the three stages can be adapted for the problem at hand. For example, the minimization for $\mathbf{m}$ in~\cite{Prins_2015} relies on a finite difference scheme while in~\cite{yadav2019monge} it uses a finite volume scheme.      
 This three-stage approach allows us to introduce two new boundary methods, viz. a \textit{segmented projection method} and a \textit{segmented arc length method}, which both lead to better results and higher computational efficiency than the original \textit{projection method}\cite{Prins_2015}. Even more importantly, the iterative method using the \textit{segmented projection method} converges in some cases when the original \textit{projection method} does not. And, as we will show, the \textit{segmented arc length method} converges for all examples.
 Furthermore, we improve upon the first minimization procedures, viz. the procedure for approximating the Jacobi matrix in the interior of the domain.
Numerical experiments have shown that grid lines in target space can intersect, preventing proper numerical convergence of our algorithm. Therefore, we introduce a method to prevent these so-called grid shocks.

The content of this paper is as follows. We discuss the theory of the least-squares method for the \MAe{} in Section~\ref{sec:LSSIntroduction}. In Section~\ref{sec:LSS_approach} the least-squares method is introduced. Afterwards, we adapt part of the method, viz. the optimization in the interior domain, in Section~\ref{sec:P_optim}. Next, in Section~\ref{sec:BoundaryMethods}, we introduce various boundary methods to replace the existing projection method and in Section~\ref{sec:crossing_grid_lines} we introduce a grid shock correction method. In Section~\ref{sec:numericalResults} we compare the boundary methods, show their weaknesses and strengths and elaborate on the convergence of the algorithm for various test cases. Lastly, we end with a discussion of the results followed by conclusions in Section~\ref{sec:conclusion}.
 
\section{The least-squares formulation}
\label{sec:LSSIntroduction}
We are interested in the two-dimensional hyperbolic \MAe{} with transport boundary condition, given by
\begin{subequations}
\label{eqn:MAbasis}
\begin{align}
\label{eqn:MAbasis_a}
\det{\left(\dJac^2 u(\mathbf{x})\right)} + f^2(\mathbf{x}, \nabla u(\mathbf{x})) & = 0, \quad \mathbf{x} \in \mathcal{X}, \\
\label{eqn:MAbasis_b}
\nabla u(\partial \mathcal{X}) & = \partial \mathcal{Y},
\end{align}
\end{subequations}
where $u = u(\mathbf{x})$ is the unknown, $\dJac^2 u$ the Hessian matrix of $u$, $f^2 > 0$ and $\mathcal{X}, \mathcal{Y} \subset \mathbb{R}^2$ connected domains.
We require the boundaries $\partial \mathcal{X}$ and $\partial \mathcal{Y}$ to be orientable. The transport boundary condition~\eqref{eqn:MAbasis_b} can be interpreted as 
\begin{subequations}
\begin{align}
  	\begin{cases}
\forall \mathbf{x} \in \partial \mathcal{X}: \nabla u(\mathbf{x}) \in \partial \mathcal{Y}, \\
\forall \mathbf{y} \in \partial \mathcal{Y} \,\,\, \exists \mathbf{x} \in \partial \mathcal{X}: \nabla u(\mathbf{x}) = \mathbf{y},
\end{cases}
\end{align}
\end{subequations}
where the latter condition is recognized as surjectivity of $\nabla u$. Bijectivity is generally not implied, not even when restricted to the boundary, as will become apparent by the example discussed in Section~\ref{sec:Annulus}.
Hyperbolicity of~\eqref{eqn:MAbasis_a} follows from the discriminant of the characteristic condition, which can be obtained by rewriting~\eqref{eqn:MAbasis_a} as
\begin{align}
F(\mathbf{x}, u, p, q, r, s, t) = r t - s^2 + f^2 = 0,
\end{align}
where $p = u_{x_1}$, $q = u_{x_2}$, $r = u_{x_1x_1}$, $s = u_{x_1x_2}$ and $t = u_{x_2x_2}$. The characteristic condition is given by~\cite[p.~10]{Bertens2022}
\begin{align}
\label{eqn:characteristicCondition}
F_r \mu^2 - F_s \mu + F_t = 0,
\end{align}
for the unknown function $\mu$, representing the slope of the characteristics. For the MA equation to be hyperbolic, two real characteristics need to exist for every point in the domain, hence the slopes of the two characteristics, and thus the roots of~\eqref{eqn:characteristicCondition}, need to be real and distinct. Henceforth, the discriminant of~\eqref{eqn:characteristicCondition} should be strictly positive.
It follows that the discriminant $\Delta$ of~\eqref{eqn:characteristicCondition} is given by
\begin{align}
  \Delta = F_s^2 - 4 F_r F_t = 4 s^2 - 4 t r = 4 f^2,
\end{align}
which is, by assumption, strictly positive. Hence, equation~\eqref{eqn:MAbasis_a} is hyperbolic. 

\subsection{Least-squares approach}
\label{sec:LSS_approach}
In~\cite{Prins_2015} a least-squares method was introduced to solve the elliptic \MAe{} given by $
\det{\left(\dJac \mathbf{m} \right)} = f^2(\mathbf{x},\nabla u(\mathbf{x}))$ for $\mathbf{x} \in \mathcal{X}$ and $\dJac \mathbf{m}$ the Jacobi matrix of $\mathbf{m}$.
The main idea of the least-squares method is to reformulate the \MAe{} in terms of the mapping $\mathbf{m}: \mathcal{X} \rightarrow \mathcal{Y}$, representing $\nabla u$, and solve for $\mathbf{m}$. Subsequently, $u$ is reconstructed from $\mathbf{m}$ in a least-squares sense. To solve the hyperbolic problem  we replace the right-hand side of the elliptic \MAe{} by $- f^2(\mathbf{x},\nabla u(\mathbf{x}))$ and substitute $\dJac \mathbf{m} = \dJac^2 u$, thus obtaining
\begin{subequations}
\label{eqn:MAbasis2}
\begin{align}
	\label{eqn:MAbasis2_a}
	\det{\left(\dJac \mathbf{m}(\mathbf{x})\right)} + f^2(\mathbf{x}, \mathbf{m}(\mathbf{x})) & = 0, \quad \mathbf{x} \in \mathcal{X}, \\
	\label{eqn:MAbasis2_b}
	\mathbf{m}(\partial \mathcal{X}) & = \partial \mathcal{Y}.
\end{align}
\end{subequations}
We formulate a minimization problem for $\mathbf{m}$ which we solve numerically. For this, we introduce the auxiliary functions $\mathbf{P} : \mathcal{X} \rightarrow \mathbb{R}^{2\times 2}$ and $\mathbf{b}: \partial \mathcal{X} \rightarrow \partial \mathcal{Y}$ which are used to approximate $\dJac \mathbf{m}$ on the whole domain and $\mathbf{m}$ on the boundary, respectively. This is achieved by the least-squares method, i.e., subsequently minimizing three separate functionals given by
\begin{subequations}
\label{eqn:def:Js}
\begin{align}
\label{eqn:def:J_I}
    J_\textrm{I}(\mathbf{m}, \mathbf{P}) & = \frac{1}{2} \iint_{\mathcal{X}} \|\dJac\mathbf{m}-\mathbf{P}\|^2 \diff \mathbf{x}, \\
\label{eqn:def:J_B}
    J_\textrm{B}(\mathbf{m}, \mathbf{b}) & = \frac{1}{2} \oint_{\partial \mathcal{X}} |\mathbf{m} - \mathbf{b}|^2 \diff s, \\
\label{eqn:def:J}        
    J(\mathbf{m}, \mathbf{P}, \mathbf{b}) & = \alpha J_\textrm{I}(\mathbf{m}, \mathbf{P}) + (1- \alpha) J_\textrm{B}(\mathbf{m}, \mathbf{b}),
\end{align}
\end{subequations}
where $|\cdot|$ is the standard 2-norm, $\| \cdot \|$ is the Frobenius norm defined by $\|\mathbf{A}\|^2 = \mathrm{Tr}(\mathbf{A} \mathbf{A}^\text{T})$ for a matrix $\mathbf{A}$ and $0 < \alpha < 1$ is a control parameter to either place weights on the boundary and the interior.
Starting with an initial guess $\mathbf{m}^0$, the iterative optimization procedure for $n = 0, 1, 2, \dots$ reads
\begin{subequations}
   	\label{eqn:iteration}
   	\begin{align}
   		\mathbf{P}^{n+1} & = \argmin_{\mathbf{P}\in \mathcal{P}(\mathbf{m}^n)} J_\textrm{I}(\mathbf{m}^n, \mathbf{P}), \\
		\mathbf{b}^{n+1} & = \argmin_{\mathbf{b}\in \mathcal{B}} J_\textrm{B}(\mathbf{m}^n, \mathbf{b}), \\
   		\label{eqn:iteration_c}
   		\mathbf{m}^{n+1} & = \argmin_{\mathbf{m} \in \mathcal{V}} J(\mathbf{m}, \mathbf{P}^{n+1}, \mathbf{b}^{n+1}).
   	\end{align}
\end{subequations}
The spaces $\mathcal{P}(\mathbf{m}^n)$, $\mathcal{B}$ and $\mathcal{V}$ follow from three key observations. First, because $\mathbf{m} = \nabla u$, the Jacobi matrix $\dJac\mathbf{m} = \dJac^2 u$ is symmetric and $\det(\dJac\mathbf{m}) = -f^2(\mathbf{x}, \mathbf{m}(\mathbf{x}))$. Secondly, by the transport boundary condition, for all $\mathbf{x} \in \partial \mathcal{X}: \mathbf{m}(\mathbf{x}) \in \partial \mathcal{Y}$. As we require $\mathbf{m}$ to be twice continuously differentiable later on, we impose this requirement. The three sets are then given by
\begin{subequations}
\begin{align}
    	\mathcal{P}(\mathbf{m}) & = \left\{ \mathbf{P} \in  [C^1(\mathcal{X})]^{2 \times 2} \mid \det(\mathbf{P}(\mathbf{x})) = -f^2(\mathbf{x}, \mathbf{m}(\mathbf{x})), \mathbf{P} = \mathbf{P}^\text{T} \right\}, \\
    \mathcal{B} & = \left\{ \mathbf{b} \in [C(\mathcal{\partial X})]^2 \mid \mathbf{b}(\mathbf{x}) \in \partial \mathcal{Y}  \right\}, \\
    \mathcal{V} & = [C^2(\mathcal{X})]^2.
\end{align}
\end{subequations}

We first outline the minimization of $J$, as it remains unchanged w.r.t.~\cite{Prins_2015}, and in the next sections we elaborate on the minimization of $J_\mathrm{I}$ and $J_\mathrm{B}$. Taking the variational derivative of \eqref{eqn:def:J} and applying the fundamental lemma of calculus of variations~\cite[p.~185]{HilbertCourantVol1} yields that for the optimal $\mathbf{m}$, each of its components should satisfy a Poisson equation with Robin boundary condition given by
\begin{subequations}
   	\label{eqn:m1PDE}
\begin{align}
    \Delta m_1 & = \nabla \boldsymbol{\cdot} \mathbf{p}_1,  && \mathbf{x} \in \mathcal{X}, \\
    (1-\alpha) m_1 + \alpha \nabla m_1 \boldsymbol{\cdot} \hat{\mathbf{n}}  & = (1-\alpha) b_1 + \alpha \mathbf{p}_1 \boldsymbol{\cdot} \hat{\mathbf{n}},  && \mathbf{x} \in \partial \mathcal{X},
\end{align}
\end{subequations}
for the first component $m_1$ and 
\begin{subequations}
	\label{eqn:m2PDE}
\begin{align}
    \Delta m_2 & = \nabla \boldsymbol{\cdot} \mathbf{p}_2,  &&  \mathbf{x} \in \mathcal{X}, \\
    (1-\alpha) m_2 + \alpha \nabla m_2 \boldsymbol{\cdot} \hat{\mathbf{n}} & = (1-\alpha) b_2 + \alpha \mathbf{p}_2 \boldsymbol{\cdot} \hat{\mathbf{n}},  && \mathbf{x} \in \partial \mathcal{X},     
\end{align}
\end{subequations}
for the second component $m_2$. The functions $\mathbf{p}_i$ ($i = 1,2$) denote the $i^\text{th}$ column of the matrix $\mathbf{P}$ and $\hat{\mathbf{n}}$ is the unit outward normal vector to $\partial \mathcal{X}$.

Upon convergence of~\eqref{eqn:iteration} we reconstruct $u$ from $\mathbf{m}$ by minimizing another least-squares functional, viz.
\begin{align}
u = \argmin_{\psi \in C^2(\mathcal{X})} \frac{1}{2} \iint_\mathcal{X} | \nabla \psi - \mathbf{m}|^2 \diff \mathbf{x}.
\end{align}
Using calculus of variations once more, we obtain the Poisson equation with Neumann boundary conditions for $u$, which reads
\begin{subequations}
    \label{eqn:u_bvp}
\begin{alignat}{2}
   	\label{eqn:u_bvp_a}
    \Delta u & = \nabla \boldsymbol{\cdot} \mathbf{m}, && \qquad \mathbf{x} \in \mathcal{X}, \\
      	\label{eqn:u_bvp_b}
    \nabla u \boldsymbol{\cdot} \hat{\mathbf{n}} & = \mathbf{m} \boldsymbol{\cdot} \hat{\mathbf{n}}, && \qquad \mathbf{x} \in \partial \mathcal{X}.
\end{alignat}
\end{subequations}
For~\eqref{eqn:u_bvp} to admit a solution, the compatibility condition~\cite[p.~184]{Abdallah1987}
\begin{align}
	\label{eqn:compatibilityPoissonEq}
	\iint_{\mathcal{X}} \nabla \boldsymbol{\cdot} \mathbf{m} \diff \mathbf{x} - \oint_{\partial \mathcal{X}} \mathbf{m} \boldsymbol{\cdot} \hat{\mathbf{n}} \diff s = 0,
\end{align}
is automatically satisfied due to the divergence theorem.

We solve the three Poisson equations using finite differences (FD), more specifically, standard second-order central differences for both the first and second order derivatives. For grid points on the boundary we introduce ghost points, which we eliminate using the normal derivatives in the Robin boundary condition. 
The system we obtain from discretizing~\eqref{eqn:m1PDE} and~\eqref{eqn:m2PDE} needs to be solved in each iteration. In order to increase computational efficiency, we compute the LU-decomposition in the initialization of the algorithm.    
Note that the solution for $u$ is not unique due to the (transport) boundary condition~\cite[p.~A1438]{Froese2012}, which is also reflected by~\eqref{eqn:u_bvp_b}, so we enforce uniqueness by fixing one function value of $u$, i.e., let $\mathbf{x} \in \mathcal{X}$ be arbitrary, we then impose the condition $u(\mathbf{x}) = 0$. In practice we assume $\mathcal{X} = [x_\mathrm{m}, x_\mathrm{M}] \times [y_\mathrm{m}, y_\mathrm{M}]$ and we impose $u(x_\mathrm{m}, y_\mathrm{m}) = 0$.  Alternatively, one could prescribe the average value of $u$ on the domain \cite[p.~177]{LotteProefschrift}.

\subsection{$\mathbf{P}$-optimization}
\label{sec:P_optim}
The matrix $\dJac \mathbf{m}$ cannot be determined exactly during the iterative process. Because the integrand of $J_{\textrm{I}}$, i.e., $\|\dJac\mathbf{m} - \mathbf{P}\|^2$, does not depend on derivatives of $\mathbf{P}$ we employ a piece-wise minimization. To this end we approximate $\dJac \mathbf{m}$ using standard finite difference. 
Let $\mathbf{x}_{ij} = ((x_1)_{i}, (x_2)_{j}) \in \mathcal{X}$ be the grid points of a Cartesian grid with $i = 1, \dots, N_{x_1}$ and $j = 1, \dots N_{x_2}$ denoting the first and second coordinate, respectively. We write $\mathbf{m}_{ij} \approx \mathbf{m}(\mathbf{x}_{ij})$ and similar for the other variables. We approximate $(\dJac \mathbf{m})_{ij}$ by $\mathbf{D}_{ij}$ using central and one-sided second-order finite differences in the interior and at the boundary, respectively.
This implies that $\mathbf{D}$ is in general not symmetric, while $\dJac \mathbf{m}$ and $\mathbf{P}$ are.
By virtue of the point-wise minimization we proceed to drop the subscripts, e.g., we write $\mathbf{m}$ instead of $\mathbf{m}_{ij}$, for brevity.

Let $F(p_{11}, p_{22}, p_{12}) = \|\mathbf{D} - \mathbf{P}\|^2$; expanding it yields
\begin{align}
    \label{eqn:def:F}
    F(p_{11}, p_{22}, p_{12}) = \frac{1}{2} \Big( (p_{11} - d_{11})^2 + (p_{12} - d_{12})^2 + (p_{12} - d_{21})^2 + (p_{22} - d_{22})^2 \Big).
\end{align}
We replace $\mathbf{D}$ by its symmetric part $\mathbf{D}_\mathrm{s} = \tfrac{1}{2}(\mathbf{D} + \mathbf{D}^\text{T})$, or written in its components, we introduce $d_\textrm{s} = \tfrac{1}{2} ( d_{12} + d_{21} )$ and
\begin{align}
\mathbf{D}_\textrm{s} = \begin{pmatrix}
d_{11} & d_\textrm{s} \\ d_\textrm{s} & d_{22}
\end{pmatrix}.
\end{align}
Furthermore, we replace $F$ by $F_\textrm{s} = \tfrac{1}{2} \|\mathbf{P} - \mathbf{D}_\textrm{s}\|^2$, i.e., 
\begin{align}
F_\textrm{s}(p_{11}, p_{22}, p_{12}) = \frac{1}{2} \Big( (p_{11} - d_{11})^2 + 2 (p_{12} - d_\textrm{s})^2 + (p_{22} - d_{22})^2 \Big).
\end{align}
To justify the replacement, note that
\begin{align}
\begin{split}
 F(p_{11}, p_{22}, p_{12}) = F_\textrm{s}(p_{11}, p_{22}, p_{12}) +  \frac{1}{4}(d_{12} - d_{21})^2,
\end{split}
\end{align}
hence, $(p_{11}, p_{22}, p_{12})$ minimizes $F$ if and only if it minimizes $F_\textrm{s}$. To obtain the minimizers, we minimize $F_\textrm{s}$ under the condition $\mathbf{P} \in \mathcal{P}(\mathbf{m})$ using Lagrange multipliers. The Lagrangian is thus given by
\begin{align}
\label{eqn:def:Lagrangian}
\Lambda(p_{11}, p_{22}, p_{12}, \lambda) = F_\textrm{s}(p_{11}, p_{22}, p_{12}) + \lambda\left(p_{11} p_{22} - p_{12}^2 + f^2 \right).
\end{align}
By setting the partial derivatives of $\Lambda$ with respect to $p_{11}$, $p_{22}$, $p_{12}$ and $\lambda$ to zero, we find that the critical points of $\Lambda$ have to satisfy
\begin{subequations}
\label{eqn:PSystemComps}
\begin{align}
	\label{eqn:PSystemComps_a}    	
p_{11} + \lambda p_{22} & = d_{11}, \\
   	\label{eqn:PSystemComps_b}
\lambda p_{11} + p_{22} & = d_{22}, \\
   	\label{eqn:PSystemComps_c}
(1 - \lambda) p_{12} & = d_\textrm{s}, \\
   	\label{eqn:PSystemComps_d}
p_{11} p_{22} - p_{12}^2 & = -f^2.
\end{align}
\end{subequations}
This system can be solved analytically and the results are given by Prins et al.~\cite[p.~B942-B947]{Prins_2015} for the elliptic \MAe{}, with $-f^2$ replaced by $f^2$ in ~\eqref{eqn:PSystemComps_d}. Unfortunately, the list of solutions is not complete as for the case $d_{11} = -d_{22}$, two roots of~\eqref{eqn:PSystemComps} are missing.
We propose a different solution strategy here. First, two remarks are in place. While minimizing $F_\textrm{s}$, the matrix $\mathbf{D}_\mathrm{s}$ and the function value of $f$ are given and both $\mathbf{P}$ and $\lambda$ have to be computed. Hence, we provide a classification in terms of $\mathbf{D}_\mathrm{s}$ and the corresponding solutions of~\eqref{eqn:PSystemComps}.
Furthermore, because the matrix $\mathbf{D}_\mathrm{s}$ is an approximation, $\det(\mathbf{D}_\mathrm{s}) \neq -f^2$ and in general $\det(\mathbf{D}_\mathrm{s}) \geq 0$ could possibly occur.
We first write the linear equations of~\eqref{eqn:PSystemComps} as 
 \begin{align}
       	\boldsymbol{\Lambda} \mathbf{p} = \mathbf{d}, \quad 
       	\boldsymbol{\Lambda} = \begin{pmatrix}
       		1 & \lambda & 0 \\ \lambda & 1 & 0 \\ 0 & 0 & 1 - \lambda
       	\end{pmatrix}, \quad 
    \mathbf{p} = \begin{pmatrix}
       	p_{11} \\ p_{22} \\ p_{12}
    \end{pmatrix}, \quad
    \mathbf{d} = \begin{pmatrix}
       	d_{11} \\ d_{22} \\ d_{12}
    \end{pmatrix}.
 \end{align}
	The vector $\mathbf{p}$ is uniquely determined when $\boldsymbol{\Lambda}$ is regular, i.e., when $0 \neq \det(\boldsymbol{\Lambda}) = (1-\lambda)^2 (1+\lambda)$. We should therefore distinguish between the cases $\lambda = 1$, $\lambda = -1$ and $\lambda \neq \pm 1$.
	
	Although we should consider the cases $\lambda = 1$, $\lambda = -1$ and $\lambda \neq \pm 1$ separately, $\mathbf{d}$ and $f$ are given and $\lambda$ and $\mathbf{p}$ are to be calculated.
Therefore we consider three cases based on $\mathbf{D}_\mathrm{s}$, viz., \textbf{\small\textit{Case 1:}} $d_{11} = d_{22}$ and $d_\mathrm{s} = 0$, \textbf{\small\textit{Case 2:}} $d_{11} = - d_{22}$ and \textbf{\small\textit{Case 3:}} all other $\mathbf{D}_\mathrm{s}$. We consider $\mathbf{D}_\mathrm{s} = \mathbf{0}$ as a special case of $d_{11} = d_{22}$ and $d_\mathrm{s} = 0$.

We start with some general results, to be used in the subsequent derivations.         
First, let $\mathrm{Tr}(\mathbf{A})$ denote the trace of a matrix $\mathbf{A}$. Using~\eqref{eqn:PSystemComps} we find
\begin{subequations}
  	\label{eqn:det_Ds_Tr_Ds}
\begin{align}
   	\label{eqn:det_Ds}
  		\delta_\mathrm{s} := \det(\mathbf{D}_\mathrm{s}) & = 
	\lambda \mathrm{Tr}(\mathbf{P})^2 - (\lambda-1)^2 f^2,\\
   	\label{eqn:Tr_Ds}
\mathrm{Tr}(\mathbf{D}_\mathrm{s}) & = 
	(\lambda + 1) \mathrm{Tr}(\mathbf{P}).
\end{align}
\end{subequations}
Solving the second equation for $\mathrm{Tr}(\mathbf{P})$ and subsequently substituting it in the first equation yields
\begin{align}
	\label{eqn:lambda_4th_order_2}
	f^2 (\lambda^2 - 1)^2 + \delta_\mathrm{s} (\lambda+1)^2 - \mathrm{Tr}(\mathbf{D}_\mathrm{s})^2 \lambda = 0.
\end{align}
Next, we consider the roots of~\eqref{eqn:lambda_4th_order_2} and the corresponding solutions $\mathbf{P}$. \\

\noindent{}\textbf{\small\textit{Case 1:}} $d_{11} = d_{22}$ and $d_\mathrm{s} = 0$, which we write as $\mathbf{D}_\mathrm{s} = d \mathbf{I}$ with $d\in\mathbb{R}$. 
We will show that this condition is equivalent with $\lambda = 1$. So, let $\mathbf{D}_\mathrm{s} = d \mathbf{I}$. We show that $\lambda = 1$ by forcing a contradiction, so, assume $\lambda \neq 1$. Then subtracting~\eqref{eqn:PSystemComps_a} from~\eqref{eqn:PSystemComps_b} gives $p_{11} = p_{22}$ and by~\eqref{eqn:PSystemComps_c} we have $p_{12} = 0$. Substitution of $p_{11} = p_{22}$ and $p_{12} = 0$ in~\eqref{eqn:PSystemComps_d} yields $p_{11}^2 = - f^2 < 0$, being a contradiction. Therefore $\lambda = 1$. 
Conversely, substitution of $\lambda = 1$ in $\boldsymbol{\Lambda}$ gives
\begin{align}
	\boldsymbol{\Lambda} = \begin{pmatrix}
		1 & 1 & 0 \\ 
		1 & 1 & 0 \\
		0 & 0 & 0 								
	\end{pmatrix}.
\end{align}
In this case the null space of $\boldsymbol{\Lambda}$ is given by $\mathcal{N}(\boldsymbol{\Lambda}) = \langle\mathbf{v}_1, \mathbf{v}_2\rangle$ with $\mathbf{v}_1 = (1, -1, 0)^\text{T}$ and $\mathbf{v}_2 = (0, 0, 1)^\text{T}$. So $\boldsymbol{\Lambda} \mathbf{p} = \mathbf{d}$ only has a solution if $\mathbf{d}$ lies in the column space of $\boldsymbol{\Lambda}$, i.e., if $d_{11} = d_{22}$ and $d_\mathrm{s} = 0$ or $\mathbf{D}_\mathrm{s} = d \mathbf{I}$ with $d \in \mathbb{R}$. Henceforth we have that $\lambda = 1$ is equivalent with $\mathbf{D}_\mathrm{s} = d \mathbf{I}$ and thus $\lambda = 1$ only occurs in \textbf{\small\textit{Case 1}}. The general solution to $\boldsymbol{\Lambda} \mathbf{p} = \mathbf{d}$ is now given by
\begin{align}
	\label{eqn:case_1_p_general}
	\mathbf{p} = (p, \,\, d-p, \,\,0)^\text{T} + \mu_1 \mathbf{v}_1 + \mu_2 \mathbf{v}_2, \qquad p, \mu_1, \mu_2 \in \mathbb{R}.
\end{align}
We aim to minimize $F_s$. Substitution of~\eqref{eqn:case_1_p_general} in $F_s$ gives
\begin{align}
F_\textrm{s}(p_{11}, p_{22}, p_{12}) = \tfrac{1}{2} \Big( (p + \mu_1- d)^2 + 2 \mu_2^2 + (p + \mu_1)^2 \Big),
\end{align}
thus showing $\mu_2 = 0$. Furthermore, the minimum is independent of the choice for $\mu_1$ as can be seen by writing $\mathbf{p} = (\tilde{p}, \,\, d-\tilde{p}, \,\,0)^\text{T}$ with $\tilde{p} = p + \mu_1$. For simplicity we choose $\mu_1 = 0$. Subsequent substitution of $\mathbf{p}$ into~\eqref{eqn:PSystemComps_d} gives $p(d - p) = - f^2$. This second order polynomial in $p$ has two real roots, viz.
\begin{align}
	p = \frac{1}{2}\left(d \pm \sqrt{d^2 + 4 f^2}\right).
\end{align}
So in total we find the two solutions
\begin{align}
	p_{11} = \frac{1}{2}\left(d \pm \sqrt{d^2 + 4 f^2}\right), \qquad p_{22} = d - p_{11}, \qquad p_{12} = 0.
\end{align}
In case $d = 0$, i.e., in case $\mathbf{D}_\mathrm{s} = \mathbf{0}$, the above derivation still holds so we consider $\mathbf{D}_\mathrm{s} = \mathbf{0}$ an instance of \textbf{\small\textit{Case 1}}.
\newline

\noindent{}\textbf{\small\textit{Case 2:}} $d_{11} = - d_{22}$, which we write as $\mathbf{d} = (d, \,\, -d, \,\, d_\mathrm{s})^\text{T}$ with $d, d_\mathrm{s}\in\mathbb{R}$. We have that $\mathrm{Tr}(\mathbf{D}_\mathrm{s}) = 0$ and $\delta_\mathrm{s} = - (d^2 + d_\mathrm{s}^2)$. For this case the fourth order polynomial~\eqref{eqn:lambda_4th_order_2} can be written as
\begin{align}
	(\lambda + 1)^2 \left((\lambda - 1)^2 + \frac{\delta_\mathrm{s}}{f^2}\right) = 0.
\end{align}
It follows that we have the three unique roots, $\lambda = -1$ (with multiplicity 2) and $\lambda = 1 \pm \sqrt{|\delta_\mathrm{s}|}/f$. 
\begin{itemize}
	\item[\textbullet]  In case $\lambda = -1$ we have
\begin{align}
	\boldsymbol{\Lambda} = \begin{pmatrix}
		1 & -1 & 0 \\ 
		-1 & 1 & 0 \\
		0 & 0 & 2
	\end{pmatrix},
\end{align}
and the corresponding null space $\mathcal{N}(\boldsymbol{\Lambda}) = \langle \mathbf{v}_3 \rangle$ with $\mathbf{v}_3 = (1, 1, 0)^\text{T}$.
For $\mathbf{p}$ to be a solution to $\boldsymbol{\Lambda} \mathbf{p} = \mathbf{d}$ we require $\mathbf{d}$ to be in the column space of $\boldsymbol{\Lambda}$. It follows that $\mathbf{d} = (d, \,\, -d, \,\, d_\mathrm{s})^\text{T}$, $d, d_\mathrm{s}\in\mathbb{R}$. Henceforth $\lambda = -1$ only occurs for \textbf{\small\textit{Case 2}}.

The general solution to  $\boldsymbol{\Lambda} \mathbf{p} = \mathbf{d}$ is therefore given by
\begin{align}
	\mathbf{p} = (p, \,\, p-d, \,\, \tfrac{1}{2} d_\mathrm{s})^\text{T} + \mu_3 \mathbf{v}_3, \qquad p, \mu_3 \in \mathbb{R}.
\end{align}
Writing $\mathbf{p} = (\tilde{p}, \,\, \tilde{p}-d, \,\, \tfrac{1}{2} d_\mathrm{s})^\text{T}$ with $\tilde{p} = p + \mu_3$ shows that the actual solution $\mathbf{p}$ does not change by choosing $\mu_3$, so we simply choose $\mu_3 = 0$. By~\eqref{eqn:PSystemComps_d} it follows that
\begin{align}
	p(p - d) - \frac{1}{4} d_\textrm{s}^2 + f^2 = 0.
\end{align}
Consequently, solving for $p$ we find that for $|\delta_\mathrm{s}|  - 4 f^2 \geq 0$ we have
\begin{align}
	p_{11} = \tfrac{1}{2}\left(d \pm \sqrt{|\delta_\mathrm{s}|  - 4 f^2 }\right), \qquad
	p_{22} = p_{11} - d, \qquad 
	p_{12} = \tfrac{1}{2}d_\textrm{s}.
\end{align}
When $|\delta_\mathrm{s}|  - 4 f^2 < 0$ the solution $\mathbf{p}$ is complex. Because we are only interested in real-valued solutions, we do not consider $\lambda = -1$ when $|\delta_\mathrm{s}|  - 4 f^2 < 0$.  

\item[\textbullet]  In the case $\lambda = 1 \pm \sqrt{|\delta_\mathrm{s}|}/f$, the matrix $\boldsymbol{\Lambda}^{-1}$ is uniquely defined \linebreak(see~\eqref{eqn:Lambda_inverse} for an explicit expression) and by $\mathbf{p} = \boldsymbol{\Lambda}^{-1} \mathbf{d}$ we obtain
\begin{align}
	\label{eqn:pMin_lambda_neq_pm1_d11_d22_negate}
	p_{11} = \mp \frac{d f}{\sqrt{|\delta_\mathrm{s}|}}, \qquad
	p_{22} = - p_{11}, \qquad
	p_{12} = \mp \frac{d_\textrm{s} f}{\sqrt{|\delta_\mathrm{s}|}}.
\end{align}
The solutions~\eqref{eqn:pMin_lambda_neq_pm1_d11_d22_negate} are new with respect to those found by Prins et al.~\cite{Prins_2015} and are not specific to the hyperbolic \MAe{}.
\end{itemize}~
\newline

\noindent{}\textbf{\small\textit{Case 3:}} All other $\mathbf{D}_\mathrm{s}$, i.e., both $\mathbf{D}_\mathrm{s} \neq d \mathbf{I}$ and $\mathbf{d}_\mathrm{s} \neq (d, \,\, -d, \,\, d_\mathrm{s})^\text{T}$ for all $d, d_\mathrm{s}\in\mathbb{R}$. By \textbf{\small\textit{Case 1}} we have $\lambda \neq 1$ and by \textbf{\small\textit{Case 2}} we have $\lambda \neq -1$. Therefore $\det\boldsymbol{\Lambda} = (1-\lambda)^2 (1+\lambda)\neq 0$. Consequently $\boldsymbol{\Lambda}$ is invertible and its inverse is given by
\begin{align}
	\label{eqn:Lambda_inverse}
	\boldsymbol{\Lambda}^{-1} =
	\frac{1}{1 - \lambda^2}
\begin{pmatrix}
		1 & - \lambda & 0 \\ 
		- \lambda & 1 & 0 \\
		0 & 0 & 1 + \lambda							
	\end{pmatrix}.
\end{align}
The values for $\lambda$ are obtained by solving~\eqref{eqn:lambda_4th_order_2}. The roots of this fourth order polynomial can be determined analytically using Ferrari's method~\cite[p.~22]{Tignol2001} and are given in~\cite[p.~B945]{Prins_2015}. For $\mathbf{p}$ we subsequently find  $\mathbf{p} = \boldsymbol{\Lambda}^{-1} \mathbf{d}$, or more explicitly
\begin{align}
\label{eqn:case_1_p_2}
	p_{11} = \frac{\lambda d_{22} - d_{11}}{\lambda^2 - 1}, \qquad
	p_{22} = \frac{\lambda d_{11} - d_{22}}{\lambda^2 - 1}, \qquad
	p_{12} = \frac{d_\textrm{s}}{1 - \lambda}.
\end{align}

\subsection{$J_\mathrm{B}$-Optimization}
\label{sec:BoundaryMethods}
In~\cite[p.~131-133]{CorienThesis} a \textit{projection method} (PM) has been proposed for the minimization of~\eqref{eqn:def:J_B}.
As our numerical results will show, this method proves insufficient for some examples. Therefore we developed two improved methods, viz., the \textit{segmented projection method} (SPM) and the \textit{segmented arc length method} (SALM). Before we introduce the boundary methods, we first introduce some notation. 

Let $\mathbf{x}_{ij} = ((x_1)_{i}, (x_2)_{j}) \in \mathcal{X}$ be the grid points of a Cartesian grid with $i = 1, \dots, N_{x_1}$ and $j = 1, \dots N_{x_2}$ denoting the first and second coordinate, respectively. Let $\mathbf{x}_\elll$ be the grid points restricted to $\partial \mathcal{X}$ for $\elll =  1, \dots N$. We index $\mathbf{x}_\elll =  1, \dots N$ in the clockwise direction such that $\mathbf{x}_1 = \mathbf{x}_{1,1}$, i.e. the first point on $\partial \mathcal{X}$ equals the point $\mathbf{x}_{ij}$ with $i = j = 1$. 
We approximate $\mathbf{m}_{ij} \approx \mathbf{m}(\mathbf{x}_{ij})$, $\mathbf{m}_{l} \approx \mathbf{m}(\mathbf{x}_{l})$ and similarly for the other variables.

The main idea behind SPM and SALM is to partition the boundaries of the source and target domains in segments. We then uniquely enforce one source segment to be mapped to one target segment. We   
then distribute $\mathbf{b}_\elll$, corresponding to $\mathbf{m}_\elll$ by either a projection (SPM) or by a ratio of arc lengths (SALM). 

Let the boundary segments of $\mathcal{X}$ be the curves $\Gamma_k^\mathcal{X} \subset \partial \mathcal{X}$ such that for $N_\Gamma$ boundary segments we have $\cup_{k=1}^{N_\Gamma} \Gamma_k^\mathcal{X} = \partial \mathcal{X}$. We denote $\Gamma_{N_\Gamma + 1}^\mathcal{X} = \Gamma_{1}^\mathcal{X}$ and assume the intersections $\Gamma_{k_1}^\mathcal{X} \cap \Gamma_{k_2}^\mathcal{X}$ contain precisely one element if $k_2 = k_1 + 1$ and no elements otherwise.
Furthermore, we require each $\Gamma_k^\mathcal{X}$ to be parametrizable. We assume similar properties for $\Gamma^\mathcal{Y}_k$.
We aim to map each boundary segment of $\partial \mathcal{X}$ to a boundary segment of $\partial \mathcal{Y}$, hence we enforce $\mathbf{m}(\Gamma_k^\mathcal{X}) = \Gamma_k^\mathcal{Y}$ for $k = 1, \dots, N_\Gamma$, from which it follows that
\begin{align}
   	\mathbf{m}(\partial \mathcal{X}) = \mathbf{m}(\cup_{k=1}^{N_\Gamma} \Gamma^\mathcal{X}_k) = \cup_{k=1}^{N_\Gamma} \mathbf{m}(\Gamma^\mathcal{X}_k) = \cup_{k=1}^{N_\Gamma} \Gamma^\mathcal{Y}_k = \partial \mathcal{Y},
\end{align} which is the required transport boundary condition.

\begin{figure}[htbp]
   	\centering
   	\begin{tikzpicture} 
   		\begin{axis}[
   			x=1.8cm, y=2cm,
   			xmin=0.2, xmax=5.5,
   			ymin=-0, ymax=2.2,
   			axis line style={draw=none},
   			tick style={draw=none},
   			ticks=none,
   			]
   			
   			\addplot [domain=0:1, samples = 6, red, only marks] ({xf(x)},{yf(x)});
   			\addplot [domain=0:1, samples = 6, red] ({xf(x)},{yf(x)});
   			
   			\addplot [domain=0:1, samples=50, blue] ({xf(x)}, {yf(x)});
   			
   			\node[blue, right, xshift=3] at (axis cs:{xf(0.1)}, {yf(0.1)}) {$\Gamma^\mathcal{Y}_k$};
   			\node[red, left] at (axis cs:{xf(0)}, {yf(0)}) {$\mathbf{y}_1$};
   			\node[red, below] at (axis cs:{xf(0.4)}, {yf(0.4)}) {$\mathbf{y}_i$};
   			\node[red, below, yshift=-3pt, xshift=6pt] at (axis cs:{xf(0.6)}, {yf(0.6)}) {$\mathbf{y}_{i+1}$};
   			\node[red, below] at (axis cs:{xf(1)}, {yf(1)}) {$\mathbf{y}_{N_\textrm{b}}$};
   			
   			\draw[black, dashed]  (axis cs:{xf(0)}, {yf(0)}) .. controls (axis cs: 0.35,1.8) .. (axis cs: 1.3,2) node [right, yshift=-16pt, xshift=-13pt] {$\Gamma^\mathcal{Y}_{k-1}$};
   			\draw[black, dashed]  (axis cs:{xf(1)}, {yf(1)}) .. controls (axis cs: 5.65, 0.8) and (axis cs: 5, 1.2) ..  (axis cs: 5.2, 1.4) node [left, above, xshift=-16pt, yshift=-30pt] {$\Gamma^\mathcal{Y}_{k+1}$};

   			\draw [decorate, decoration={brace,amplitude=10pt,raise=1pt}, red] (axis cs:{xf(0.4)}, {yf(0.4)}) -- (axis cs:{xf(0.6)}, {yf(0.6)}) node [red, midway, yshift=15pt, outer sep=3pt]{$\tau_i$};
   		\end{axis}
   		
   	\end{tikzpicture}
   	\caption{Schematic overview of the discretization of $\partial \mathcal{Y}$.}
   	\label{fig:SALM}
\end{figure}
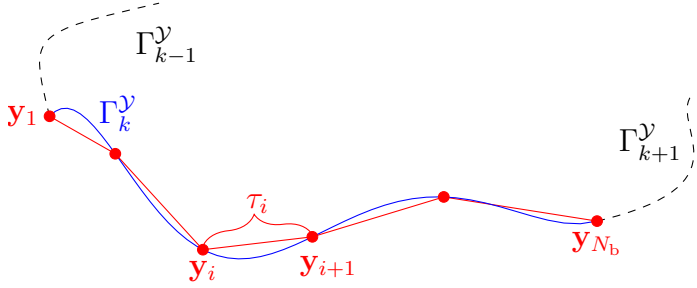

Figure~\ref{fig:SALM} shows a part of $\partial \mathcal{Y}$ and (parts of) three boundary segments. In the following we fix $k$ and for brevity drop the subscript in $\Gamma_k^\mathcal{X}$ and $\Gamma_k^\mathcal{Y}$. 

Let $\mathbf{y}_i$ for $i = 1, \dots, N_\mathrm{b}$ be a discretization of the boundary segment $\Gamma^\mathcal{Y}$ such that for a given counter clockwise parametrization $\mathbf{y}(s): [0,1] \rightarrow \Gamma^\mathcal{Y}$, we have $\mathbf{y}_1 = \mathbf{y}(0)$ and $\mathbf{y}_{N_\mathrm{b}} = \mathbf{y}(1)$. We choose to parametrize $\partial \mathcal{X}$ and $\partial \mathcal{Y}$ in opposite directions because $\mathbf{m}$ reverts the direction if it is a solution to the hyperbolic \MA{}.

\pgfmathdeclarefunction{xf}{1}{\pgfmathparse{2 * (#1 + 1/2) * (#1 + 1/2)}}
\pgfmathdeclarefunction{yf}{1}{\pgfmathparse{2 * pow(cos(2.5 * deg(#1)),3) - 0.75 + 1/3 * pow(#1, 16) + 2 * #1}}

~\newline    
\noindent{}\textbf{\textit{Projection method.}}
We briefly explain PM as introduced by \cite[p.~131-133]{CorienThesis}. We perform the following for each approximation $\mathbf{m}_\elll$ individually. Let $N_\Gamma = 1$, i.e., we consider the whole boundary as one boundary segment. Furthermore, let $\mathbf{y}_{N_\textrm{b}+1} = \mathbf{y}_1$, we connect adjacent points $\mathbf{y}_i$ and $\mathbf{y}_{i+1}$ by straight line segments. The projection of $\mathbf{m}_\elll$ onto the line connecting $\mathbf{y}_i$ and $\mathbf{y}_{i+1}$ is given by
\begin{subequations}
	\begin{align}
		\mathbf{m}^\mathrm{P}_i (t_i) & = \mathbf{y}_i + t_i (\mathbf{y}_{i+1} - \mathbf{y}_i), \\
		t_i & = \frac{(\mathbf{m}_\elll-\mathbf{y}_i) \boldsymbol{\cdot} (\mathbf{y}_{i+1} - \mathbf{y}_i)}{|\mathbf{y}_{i+1} - \mathbf{y}_i|^2}.
	\end{align}
\end{subequations}
As only $0\leq t_i \leq 1$ corresponds to a point on the line segment between $\mathbf{y}_i$ and $\mathbf{y}_{i+1}$, we limit $t_i$ according to $\hat{t}_i = \min(1, \max(0, t_i))$. Among all possible line segments, we choose $\mathbf{b}_\elll$ corresponding to $\mathbf{m}_\elll$, such that the distance is smallest, i.e.,
\begin{subequations}
	\begin{align}
		i_\mathrm{min} & = \argmin_i \{|\mathbf{m}^\mathrm{P}_i (\hat{t}_i) - \mathbf{m}_\elll|\}, \\
		\mathbf{b}_\elll &  = \mathbf{m}^\mathrm{P}_{i_\mathrm{min}} (\hat{t}_\mathrm{min}).
	\end{align}
\end{subequations}
~\newline    
\noindent{}\textbf{\textit{Segmented projection method.}}
Let $\Gamma_k^\mathcal{X} \subseteq \partial \mathcal{X}$ and $\Gamma_k^\mathcal{Y} \subseteq \partial \mathcal{Y}$ be boundary segments of the source and target, respectively.  For this method, we apply PM to the individual boundary segments instead of the whole boundary at once. Furthermore, we set 
\begin{align}
	\mathbf{b}(\Gamma_{k}^\mathcal{X} \cap \Gamma_{k+1}^\mathcal{X}) = \Gamma_{k}^\mathcal{Y} \cap \Gamma_{k + 1}^\mathcal{Y}, \qquad 1 \leq k \leq N_\Gamma
\end{align}
meaning, we map the end points of the source segments to the end points of the corresponding target segments. In practice, these end points are the corners of the source and target domains.

~\newline    
\noindent{}\textbf{\textit{Segmented arc length method.}} 
The core idea of this method is as follows: if $\mathbf{m}(\Gamma^\mathcal{X}) = \Gamma^\mathcal{Y}$, then the arc length of the curve $\mathbf{m}(\Gamma^\mathcal{X})$ should be equal to the arc length of the curve $\Gamma^\mathcal{Y}$. Numerically we approximate this condition by approximating the arc length of both $\Gamma^\mathcal{Y}$ and the distance between the points $\{\mathbf{m}(\mathbf{x}_\elll) \mid \mathbf{x}_\elll \in \Gamma^\mathcal{X}\}$.

We start with the arc length of the curve $\Gamma^\mathcal{Y}$. We approximate the arc length between $\mathbf{y}_{i}$ and $\mathbf{y}_{i+1}$ along $\Gamma^\mathcal{Y}$ by the length of the line segment connecting $\mathbf{y}_{i}$ and $\mathbf{y}_{i+1}$. We denote the approximation by
\begin{align}
\tau_i = |\mathbf{y}_{i+1} - \mathbf{y}_{i}|, \quad i = 1, \dots N_\mathrm{b}-1.
\end{align}
The approximate cumulative arc length between $\mathbf{y}_{1}$ and $\mathbf{y}_{i}$ in the direction of increasing $s$ is then given by
\begin{align}
t_i = \sum_{j=1}^{i-1} \tau_j, \quad i = 1, \dots, N_\mathrm{b}.
\end{align}
The total arc length from $\mathbf{y}_1$ to $\mathbf{y}_{N_\mathrm{b}}$ is then approximated by $L = t_{N_\textrm{b}}$. We use the cumulative arc lengths to introduce a piece-wise linear interpolation $\mathbf{b}_\mathrm{int}$ approximating $\mathbf{y}(s)$, viz.
\begin{align}
\mathbf{b}_\mathrm{int}(t) = \mathbf{y}_{i} + \frac{t - t_{i}}{t_{i+1} - t_{i}} (\mathbf{y}_{i+1} - \mathbf{y}_{i}), \quad t_{i} \leq t \leq t_{i+1},
\end{align}
where the scalar factor is a scaled coordinate between $\mathbf{y}_{i}$ and $\mathbf{y}_{i+1}$. Note that by construction $\mathbf{b}_\mathrm{int}$ satisfies
\begin{align}
    \mathbf{b}_\mathrm{int}(t_i) = \mathbf{y}_{i}, \quad i = 1, \dots, N_\mathrm{b},
\end{align}
and is an approximation of $\partial \mathcal{Y}$.

Next we consider the points $\{\mathbf{m}(\mathbf{x}_\elll) \mid \mathbf{x}_\elll \in \Gamma^\mathcal{X}\}$.
Let $N_\mathrm{m}$ be the number of grid points on $\Gamma^\mathcal{X}$ such that $\mathbf{x}_\elll \in \Gamma^\mathcal{X}$ for $\elll = 1, \dots, N_\mathrm{m}$. Furthermore, let
\begin{align}
\sigma_\elll = | \mathbf{m}_{\elll+1} - \mathbf{m}_{\elll} |, \quad \elll, \dots, N_\textrm{m}-1,
\end{align}
be an approximation of the arc length from $\mathbf{m}(\mathbf{x}_{\elll})$ to $\mathbf{m}(\mathbf{x}_{\elll+1})$ along $\partial \mathcal{Y}$.
This again introduces a cumulative arc length and a total arc length, respectively, given by
\begin{align}
s_\elll  = \sum_{j=1}^{\elll-1} \sigma_j, \quad \elll = 1, \dots, N_\mathrm{m}, \qquad
\tilde{L}  = s_{N_\mathrm{m}}.
\end{align}
Because $\mathbf{m}_\elll$ is an approximation and $\Gamma^\mathcal{Y}$ is approximated by straight line segments, $\tilde{L} \neq L$ in general. Hence, $s_\elll \neq L$ may occur such that the end points of $\Gamma^\mathcal{X}$ may not be mapped to the end points of $\Gamma^\mathcal{Y}$. We fix this by letting
\begin{align}
\tilde{s}_\elll = \frac{L}{\tilde{L}} s_\elll.
\end{align}
It follows that $\mathbf{b}_\mathrm{int}(\tilde{s}_\elll)$ forms a proper approximation for $\mathbf{m}_{\elll}$ restricted to $\Gamma^\mathcal{Y}$, viz.
\begin{align}
\label{eqn:b_points_SALM}
\mathbf{b}_{\elll} = \mathbf{b}_\mathrm{int}(\tilde{s}_\elll), \quad  \elll = 1 , \dots, N_\mathrm{m}.
\end{align}

\subsection{Grid shock correction}
\label{sec:crossing_grid_lines}
\begin{figure}[t!]
	\centering
	\includegraphics[width = 0.45\linewidth]{"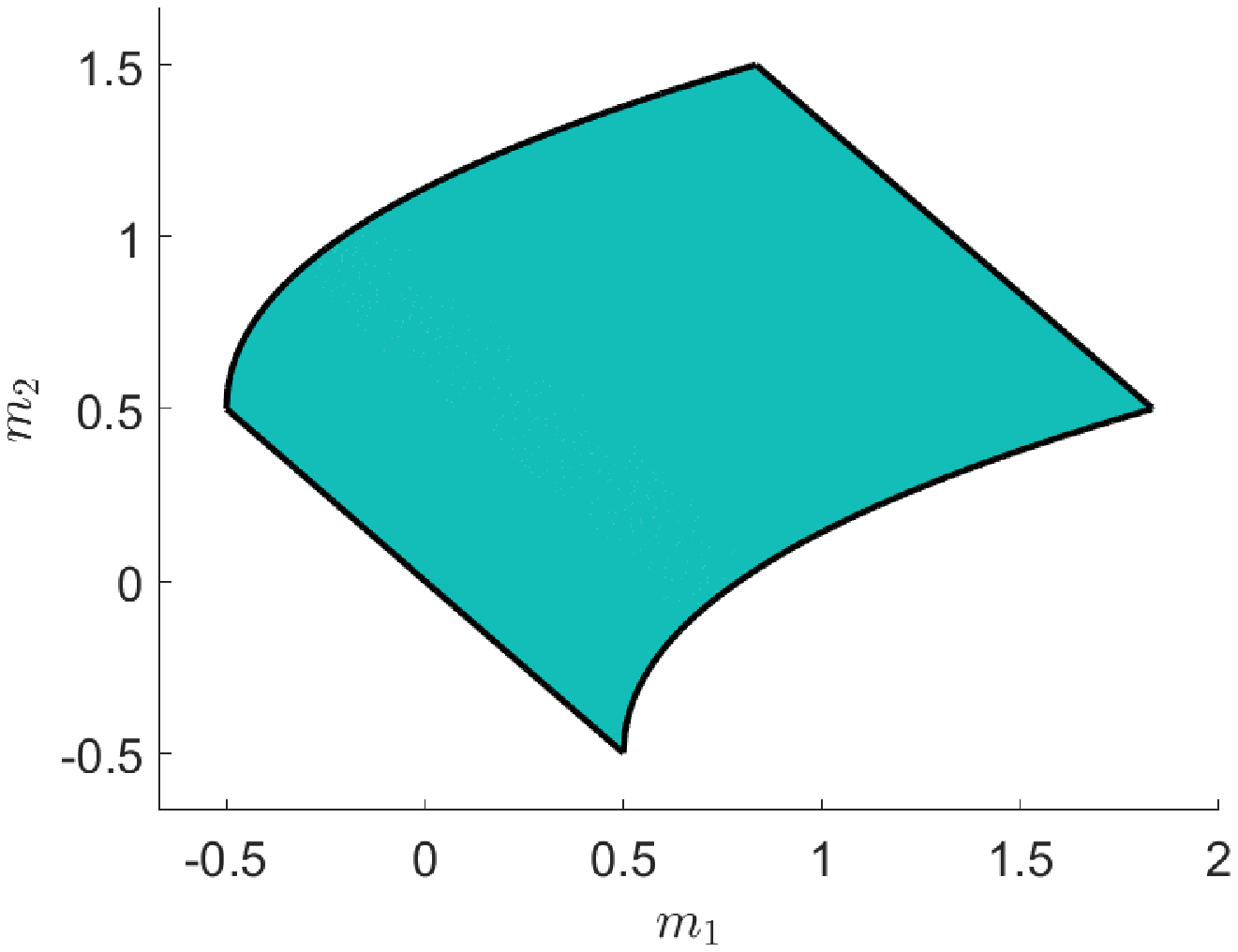"}
	\hspace{20pt}
	\includegraphics[width = 0.45\linewidth]{"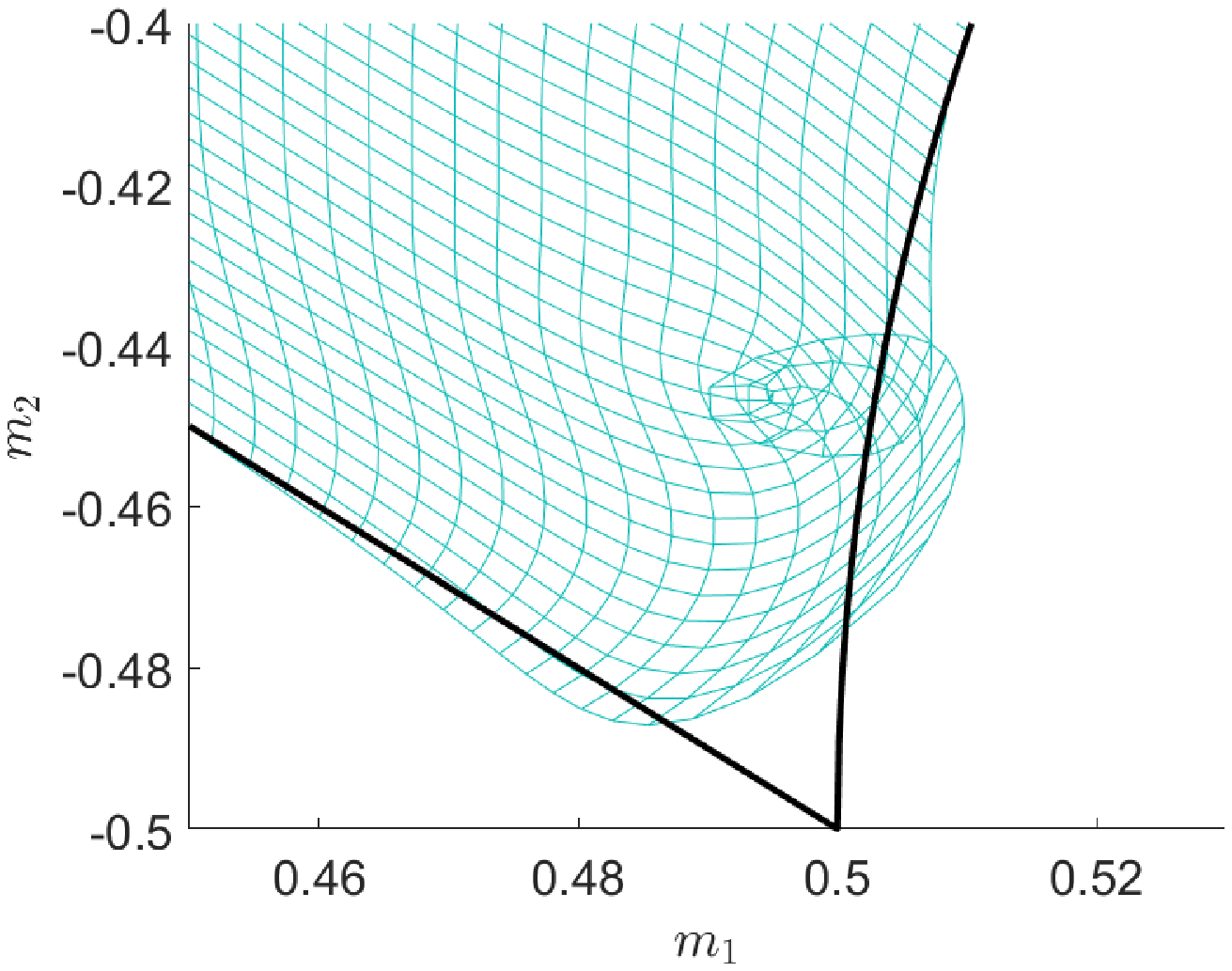"}
	\caption{Example of grid shock. The global numerical approximation is shown on the left, and a zoomed-in version on the right. }
	\label{fig:crossingGridline}
\end{figure}

Using the methods outlined above, it is possible that the approximation $\mathbf{m}^n$ of $\mathbf{m}$ contains crossing grid lines, also known as grid shocks~\cite{CORDOVA1988}. This phenomenon is shown in Figure~\ref{fig:crossingGridline} for an example we discuss in Section~\ref{sec:DeformedSquare}, with grid parameters $N_{x_1} = N_{x_2} = 321$ after $n = 15,000$ iterations. Though the solution on the left may look visually correct, the grid shock, as seen on the right, prevents proper numerical convergence of our algorithm.
   	\begin{figure}[t!]
   	\centering
   	\begin{tikzpicture}
   		\draw [dashed] (1,0) |- (5,0) node [right] {$\partial \mathcal{X}$};
   		
   		\node at (1.5, 0  )[circle,fill,inner sep=1.5pt]{};   
   		\node at (1.5, 0  )[below]{$\mathbf{x}_{i-1, 1}$};
   		
   		\node at (3,   0  )[circle,fill,inner sep=1.5pt]{};   
   		\node at (3,   0  )[below]{$\mathbf{x}_{i, 1}$};   
   		
   		\node at (4.5, 0  )[circle,fill,inner sep=1.5pt]{};   
   		\node at (4.5, 0  )[below]{$\mathbf{x}_{i+1, 1}$};
   		
   		\node at (1.5, 1.5)[circle,fill,inner sep=1.5pt]{}; 
   		\node at (1.5, 1.5)[left]{$\mathbf{x}_{i-1, 2}$};
   		
   		\node at (3,   1.5)[circle,fill,inner sep=1.5pt]{};   
   		\node at (3,   1.5)[below right]{$\mathbf{x}_{i, 2}$};     
   		
   		\node at (4.5, 1.5)[circle,fill,inner sep=1.5pt]{};  
   		\node at (4.5,   1.5)[below right]{$\mathbf{x}_{i+1, 2}$};     			
   		
   		\node at (3,   3  )[circle,fill,inner sep=1.5pt]{};  
   		\node at (3,   3  )[right]{$\mathbf{x}_{i, 3}$};      	        	        	         	        
   		
   		\draw [-] (1.5, 1.5) |- (4.5, 1.5);
   		\draw [-] (3, 0) |- (3,3);
   	\end{tikzpicture}
   	\qquad\quad
   	\begin{tikzpicture}
   		\draw [dashed] (4.95,-0.5) -- (5,0) -- (5.2, 1.75) -- (5.5, 3.5) -- (5.6, 4) node [above right] {$\partial \mathcal{Y}$};
   		
   		\node at (5, 0   )[circle,fill,inner sep=1.5pt]{};   
   		\node at (5, 0   )[right]{$\mathbf{y}_k$};   
   		
   		\node at (5.2, 1.75)[circle,fill,inner sep=1.5pt]{};  
   		\node at (5.2, 1.75   )[right]{$\mathbf{y}_{k+1}$};   
   		
   		\node at (5.5, 3.5 )[circle,fill,inner sep=1.5pt]{};   
   		\node at (5.5, 3.5 )[right]{$\mathbf{y}_{k+2}$};   
   		
   		\node at (2.7,   0  )[circle,fill,inner sep=1.5pt]{};   
   		\node at (2.7,   0  )[below]{$\mathbf{m}_{i, 3}$};   
   		
   		\node at (1.5, 1.4)[circle,fill,inner sep=1.5pt]{}; 
   		\node at (1.5, 1.4)[left]{$\mathbf{m}_{i+1, 2}$};
   		
   		\node at (3,   1.5)[circle,fill,inner sep=1.5pt]{};   
   		\node at (3,   1.5)[below left]{$\mathbf{m}_{i, 2}$};     
   		
   		\node at (4.3, 1.3)[circle,fill,inner sep=1.5pt]{};  
   		\node at (4.3, 1.3)[below]{$\mathbf{m}_{i-1, 2}$};     			
   		
   		\node at (2.7,   3  )[circle,fill,inner sep=1.5pt]{};  
   		\node at (2.7,   3  )[above]{$\mathbf{m}_{i, 1}$};      
   		
   		\node at (5.1, 0.875)[circle,fill,inner sep=1.5pt]{};  
   		\node at (5.1, 0.875)[right]{$\mathbf{b}_{i, 1}$};         	        	         	        
   		
   		\draw [-] (1.5, 1.4) -- (3, 1.5) -- (4.3, 1.3);
   		\draw [-] (2.7, 0) -- (3,1.5) -- (2.7,3);
   		
   		\node at (1.4,   2.6 )[circle,fill,inner sep=1.5pt]{};  
   		\node at (1.4,   2.6 )[above]{$\mathbf{m}_{i+1, 1}$};    
   		\node at (4.4,   3 )[circle,fill,inner sep=1.5pt]{};  
   		\node at (4.4,   3 )[above]{$\mathbf{m}_{i-1, 1}$};    
   	\end{tikzpicture}
   	\caption{Schematic overview of stencil used for detecting grid shocks. On the right $\mathbf{m}_{i-1, 2}$ is closer to $\mathbf{b}_{i,1}$ than $\mathbf{m}_{i, 1}$ is, so a grid shock occurs.}
   	\label{fig:GridShockCorrection}
\end{figure}
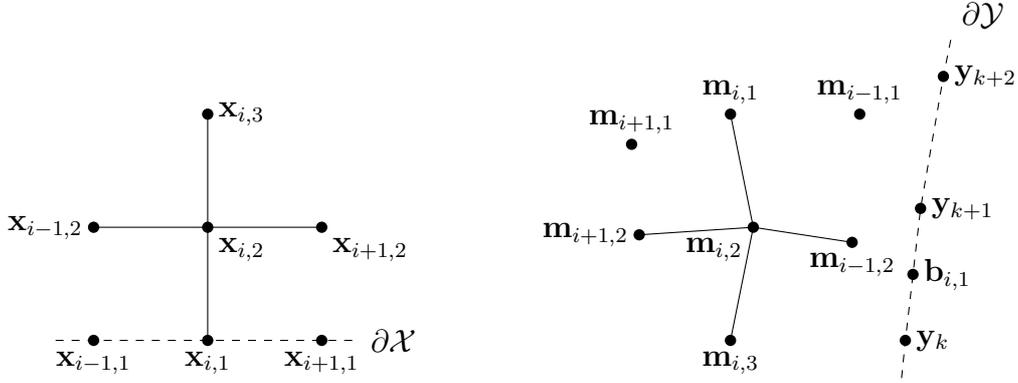
To resolve grid shocks, consider a point $\mathbf{x}_{ij} \in \partial \mathcal{X}$ as shown for $j = 1$ in Figure~\ref{fig:GridShockCorrection} on the left, and the corresponding image $\mathbf{m}_{ij}$ shown on the right. If both $\mathbf{m}_{ij}$ and $\mathbf{b}_{ij}$ are exact, then $|\mathbf{m}_{ij} - \mathbf{b}_{ij}| = 0$ and $|\mathbf{m}_{kl} - \mathbf{b}_{ij}| > 0$ for all $(k,l) \neq (i, j)$. Because both $\mathbf{m}$ and $\mathbf{b}$ are approximated, $|\mathbf{m}_{ij} - \mathbf{b}_{ij}| \neq 0$ in general. To detect grid shocks, we compute the distance $|\mathbf{m}_{ij} - \mathbf{b}_{ij}|$ for all $(k, l)$ such that $|(k,l)^\text{T} - (i,j)^\text{T}| \leq 2$. If the minimum distance if found for $(k, l) \neq (i, j)$, then we assume a grid shock orrcurs and we recompute $\mathbf{m}$.
We do so by making $\alpha$ in~\eqref{eqn:def:J} dependent on the coordinate, i.e., $\alpha = \alpha(\mathbf{x})$ and $\alpha_{ij} = \alpha(\mathbf{x}_{ij})$ and subsequently reduce $\alpha_{ij}$ on the boundary which puts more emphasis on the minimization of $|\mathbf{m}_{ij} - \mathbf{b}_{ij}|$.

Introducing the $\mathbf{x}$-dependency, the coefficients $\alpha$ and $1-\alpha$ in~\eqref{eqn:def:J} formally have to be moved inside the integrals of $J_\textrm{I}$ and $J_\textrm{B}$. After doing so, we compute the first variation of $J(\mathbf{m}, \mathbf{P}, \mathbf{b}, \alpha) = \frac{1}{2} \iint_{\mathcal{X}} \alpha \|\dJac\mathbf{m}-\mathbf{P}\|^2 \diff \mathbf{x} +  \frac{1}{2} \oint_{\partial \mathcal{X}} (1-\alpha) |\mathbf{m} - \mathbf{b}|^2 \diff s$ w.r.t. $\mathbf{m}$ and apply the fundamental lemma of calculus of variations. Consequently we obtain
\begin{subequations}
   	\label{eqn:m1PDE_2}
   	\begin{align}
   		\nabla \alpha \boldsymbol{\cdot} \nabla m_1 + \alpha \Delta m_1 & = \mathbf{p}_{1} \boldsymbol{\cdot} \nabla \alpha + \alpha \nabla \boldsymbol{\cdot} \mathbf{p}_1 ,  && \mathbf{x} \in \mathcal{X}, \\
   		(1-\alpha) m_1 + \alpha \nabla m_1 & = (1-\alpha) b_1 + \alpha \mathbf{p}_1 \cdot \hat{\mathbf{n}},  && \mathbf{x} \in \partial \mathcal{X},
   	\end{align}
\end{subequations}
for the first component of $\mathbf{m}$.
If $\alpha$ is constant in the interior, we have $\nabla \alpha = \mathbf{0}$ in the interior and equations~\eqref{eqn:m1PDE_2} reduce to~\eqref{eqn:m1PDE}. By analogy, we have~\eqref{eqn:m2PDE} for the second component of $\mathbf{m}$.
Let $\alpha_1 \in (0,1)$. We set $\alpha(\mathbf{x}_{ij}) = \alpha_1$ for $\mathbf{x}_{ij}$ in the interior of $\mathcal{X}$. For a boundary point $\mathbf{x}_{ij} \in \partial \mathcal{X}$ we instead set
\begin{align}
   	\label{eqn:discont_alpha}
   	\alpha_{ij} = 
   	\begin{cases}
		\alpha_2, & \quad \text{if} \quad \min\limits_{(k,l) \in \mathcal{I}(i,j)} |\mathbf{m}_{kl} - \mathbf{b}_{ij}|  < |\mathbf{m}_{ij} - \mathbf{b}_{ij}|, \\
		\alpha_1 & \quad \text{otherwise},
   	\end{cases}
\end{align}
where $\mathcal{I}(i,j) = \{(k,l) \mid \mathbf{x}_{kl} \in \mathrm{int}(\mathcal{X}), \, |(k,l)^\text{T} - (i,j)^\text{T}| \leq 2\}$ is the space over which we minimize, $\mathrm{int}(\mathcal{X})$ the interior of $\mathcal{X}$ and $\alpha_2 \in (0,\alpha_1)$ a constant. We choose a distance of 2 and the values $\alpha_2 = 0.005$ and $\alpha_1 = 0.2$ since these have proven to work well in practice.

If in the $n^\text{th}$ iteration we obtain $\alpha(\mathbf{x}) \not\equiv \alpha_1$, we solve~\eqref{eqn:m1PDE} and~\eqref{eqn:m2PDE} for a second time with the updated $\alpha$ to perform a correction.

Recall, we use a finite difference method for inverting the Poisson equations~\eqref{eqn:m1PDE} and~\eqref{eqn:m2PDE}, yielding a system of equations. This system of equations depends on $\alpha$. Without grid shock correction, a LU-factorization can be calculated once and used for each subsequent iteration making the inversion of the system efficient. In case $\alpha_{ij} \neq \alpha_1$ for any $(i,j)$, the same LU-factorization can no longer be used due to the component $\alpha \nabla m_k \boldsymbol{\cdot} \hat{\mathbf{n}}$ in the Robin boundary conditions and a new LU-factorization has to be calculated for the iteration.

\section{Numerical results}
\label{sec:numericalResults}
In this section we present numerical results for five examples. For each example we know the exact solution and compare the numerical methods. We choose $\mathcal{X} = [x_1^m, x_1^M] \times [x_2^m, x_2^M]$, a rectangle which may vary per case. For each example we choose $\mathcal{Y}$ such that it has a unique feature to it. We measure the residual
\begin{align}
    \epsilon_r & = \Big|D_{x_1}[m_1]_{ij} D_{x_2}[m_2]_{ij}- D_{x_1}[m_2]_{ij} D_{x_2}[m_1]_{ij} + f^2(\mathbf{x}_{ij}, \mathbf{m}_{ij})\Big|_\infty,
\end{align}%
with $D_{x_1}$ and $D_{x_2}$ standard second-order (central in interior and one-sided on boundary) finite difference operators for the first-order derivatives with respect to $x_1$ and $x_2$, respectively. 
Furthermore, we measure the global discretization errors $\epsilon_u$, $\epsilon_{m_1}$ and $\epsilon_{m_2}$ defined by
\begin{align}
\begin{split}
    \epsilon_u & = \big|\big(u_{ij} - u_{11}\big) - \big(u(\mathbf{x}_{ij}) - u(\mathbf{x}_{11})\big)\big|_\infty, \\
    \epsilon_{m_1} & = \big|(m_1)_{ij} - m_1(\mathbf{x}_{ij})\big|_\infty, \\
    \epsilon_{m_2} & = \big|(m_2)_{ij} - m_2(\mathbf{x}_{ij})\big|_\infty, \\
\end{split}
\end{align}
where the terms $u_{11}$ and $u(\mathbf{x}_{11})$ are introduced due to the nonuniqueness of $u$ given $\mathbf{m}$; see the discussion following equations~\eqref{eqn:u_bvp}. Any fixed grid point could be used, here $\mathbf{x}_{11}$ is chosen. The choice for the $\infty$-norm is arbitrary in the sense that any standard norm would give similar results. However, the $\infty$-norm is more sensitive to differences in the local errors than, for example, the standard 2-norm.
Starting the least-squares algorithm requires an initial guess $\mathbf{m}^0$, so we introduce $\widetilde{\mathcal{Y}} = [y_1^m, y_1^M] \times [y_2^m, y_2^M]$, the smallest bounding box of $\mathcal{Y}$.
We then choose $\partial \mathbf{m}^0(\mathcal{X}) = \partial \widetilde{\mathcal{Y}}$, such that $\mathbf{m}^0_{ij}$ is equidistantly distributed, i.e., $\mathbf{m}^0_{ij}$ is the result a bilinear uniform interpolation of the bounding box of $\mathcal{Y}$ with $\det(\dJac \mathbf{m}^0) < 0$. The initial guess then reads
\begin{subequations}
    \label{eqn:initial_guess_m}
\begin{align}
    (m^0_1)_{ij} & = \frac{(x_1)_{ij} - x_1^m}{x_1^M - x_1^m} y_1^M + \frac{x_1^M - (x_1)_{ij}}{x_1^M - x_1^m} y_1^m, \\
    (m^0_2)_{ij} & = \frac{x_2^M - (x_2)_{ij}}{x_2^M - x_2^m} y_2^M + \frac{(x_2)_{ij} - x_2^m}{x_2^M - x_2^m} y_2^m.
\end{align}
\end{subequations}
The initial guess is a (discretized) solution of the hyperbolic \MAe{} with $f^2 = \mathrm{area}(\widetilde{\mathcal{Y}}) / \mathrm{area}(\mathcal{X})$, $\mathbf{m}= \nabla u$ and $u = \tfrac{1}{2}(x^2 - y^2) f$. Three more such initial guesses exist, viz., 
$u = \tfrac{1}{2}(y^2 - x^2) f$ and $u = \pm x y f$.

We segmentate the source boundary in segments, clockwise, according to
\begin{align}
\begin{split}
\Gamma_1^\mathcal{X} & = \{x_1^m\} \times [x_2^m, x_2^M], \qquad
\Gamma_2^\mathcal{X} = [x_1^m, x_1^M] \times \{x_2^M\}, \\
\Gamma_3^\mathcal{X} & = \{x_1^M\} \times [x_2^m, x_2^M], \qquad
\Gamma_4^\mathcal{X} = [x_1^m, x_1^M] \times \{x_2^m\},
\end{split}
\end{align}
and we write 
\begin{align}
	\Gamma^\mathcal{Y}_k = \left\{ \Gamma^\mathcal{Y}_k(s)  \, \Big| \, s \in [0,1] \right\},
\end{align}
for $k = 1, \dots, 4$.
Furthermore, we apply grid shock correction only for iteration step $n\geq 100$, as the distance between the boundary of the initial guess and the boundary of the target may be large for small $n$. 

Lastly, we stop the iteration~\eqref{eqn:iteration} based on the update of $\mathbf{m}^n$, i.e., based on
\begin{align}
   	\Delta m^{n} = | \mathbf{m}^{n} - \mathbf{m}^{n-1}|,
\end{align}
instead of, the already introduced measures, $J_\textrm{I}$ and $J_\textrm{B}$. This is because the values for $J_\textrm{I}$ and $J_\textrm{B}$ may stagnate over the iterations while $\Delta m^{n}$ is still changing. Conversely, if $\Delta m^n$ has stagnated, then so have the functionals $J_\textrm{I}$ and $J_\textrm{B}$. We stop the iterative process when $\Delta m^n$ reaches floating-point precision.
    
 \subsection{Annulus segment}
\begin{figure}[t!]
	\centering
	\includegraphics[width = 0.33\linewidth]{"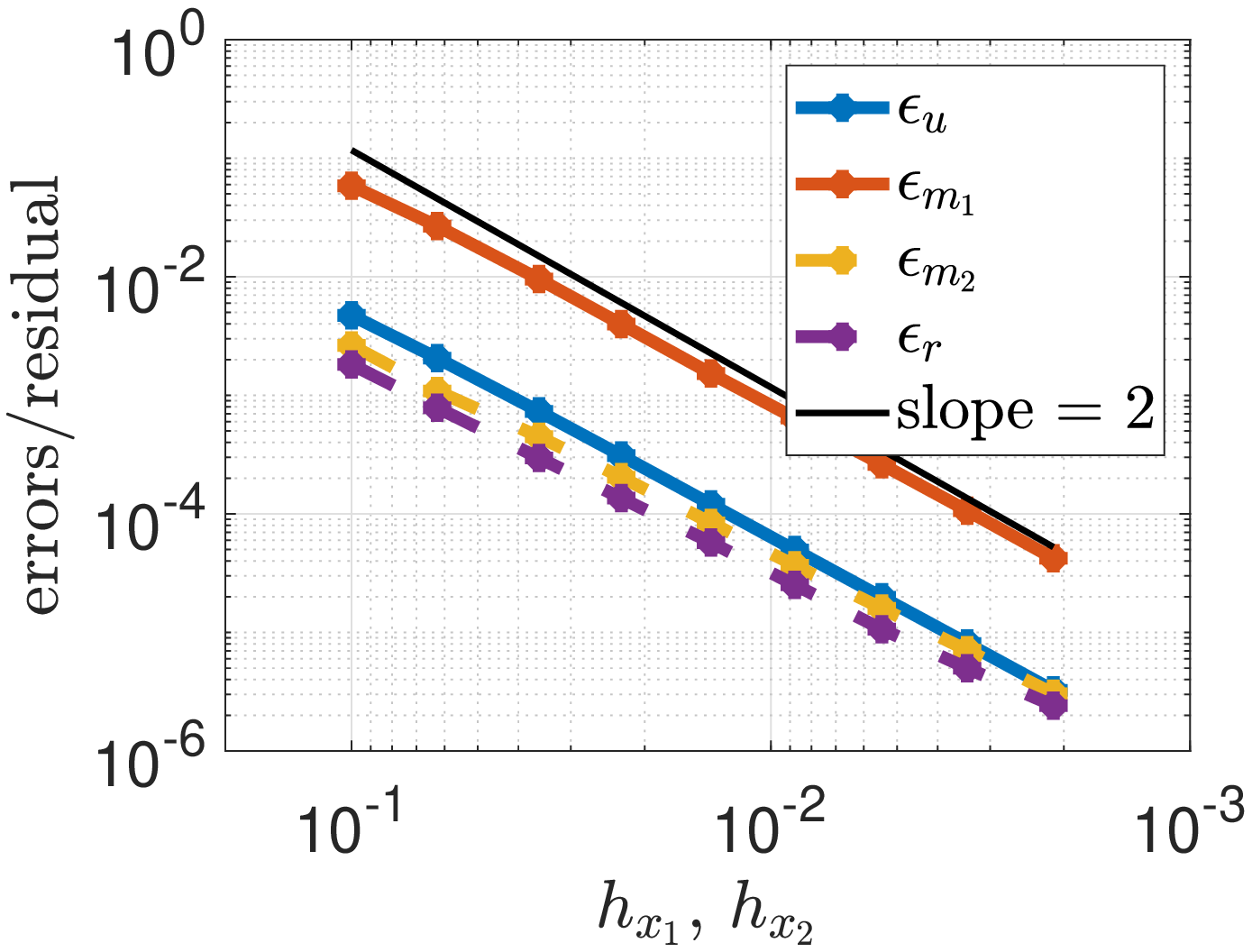"}
	\hspace{-5pt}
	\includegraphics[width = 0.33\linewidth]{"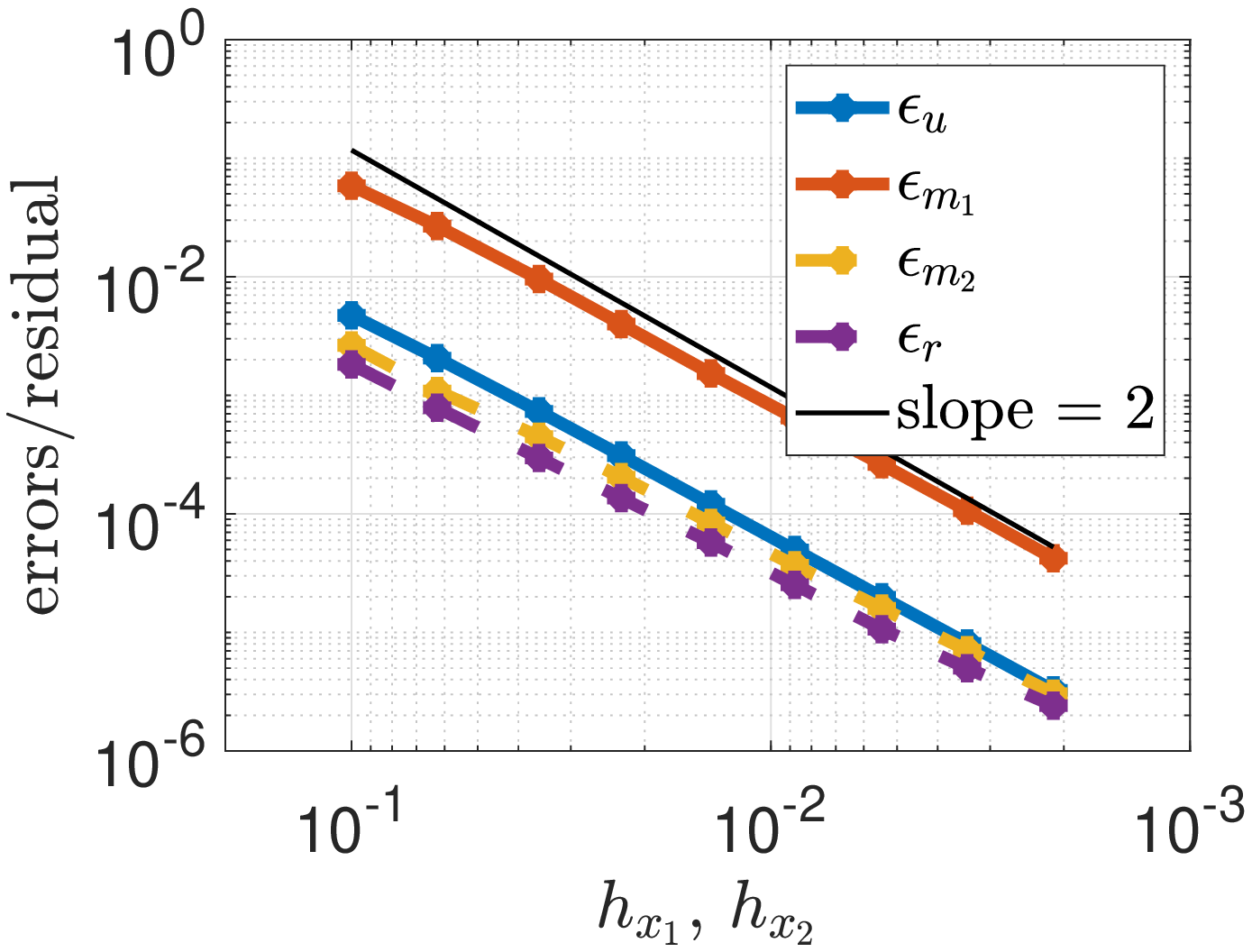"}
	\hspace{-5pt}
	\includegraphics[width = 0.33\linewidth]{"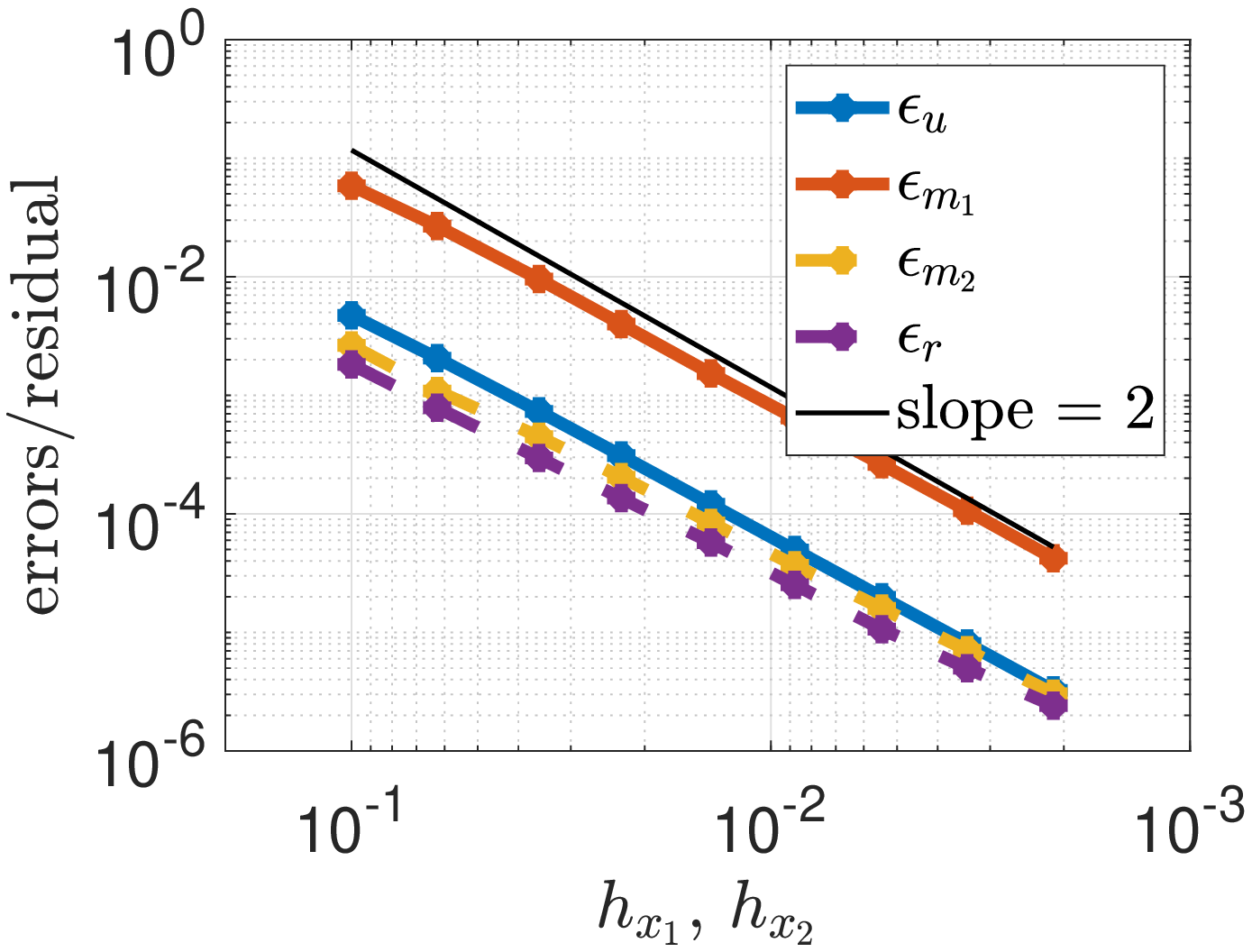"}
	\caption{ Global error and the residual for PM (left), SPM (middle) and SALM (right). }
	\label{fig:Conv_results_case_5}
\end{figure}

\begin{figure}[t!]
	\centering
	\includegraphics[width = 0.33\linewidth]{"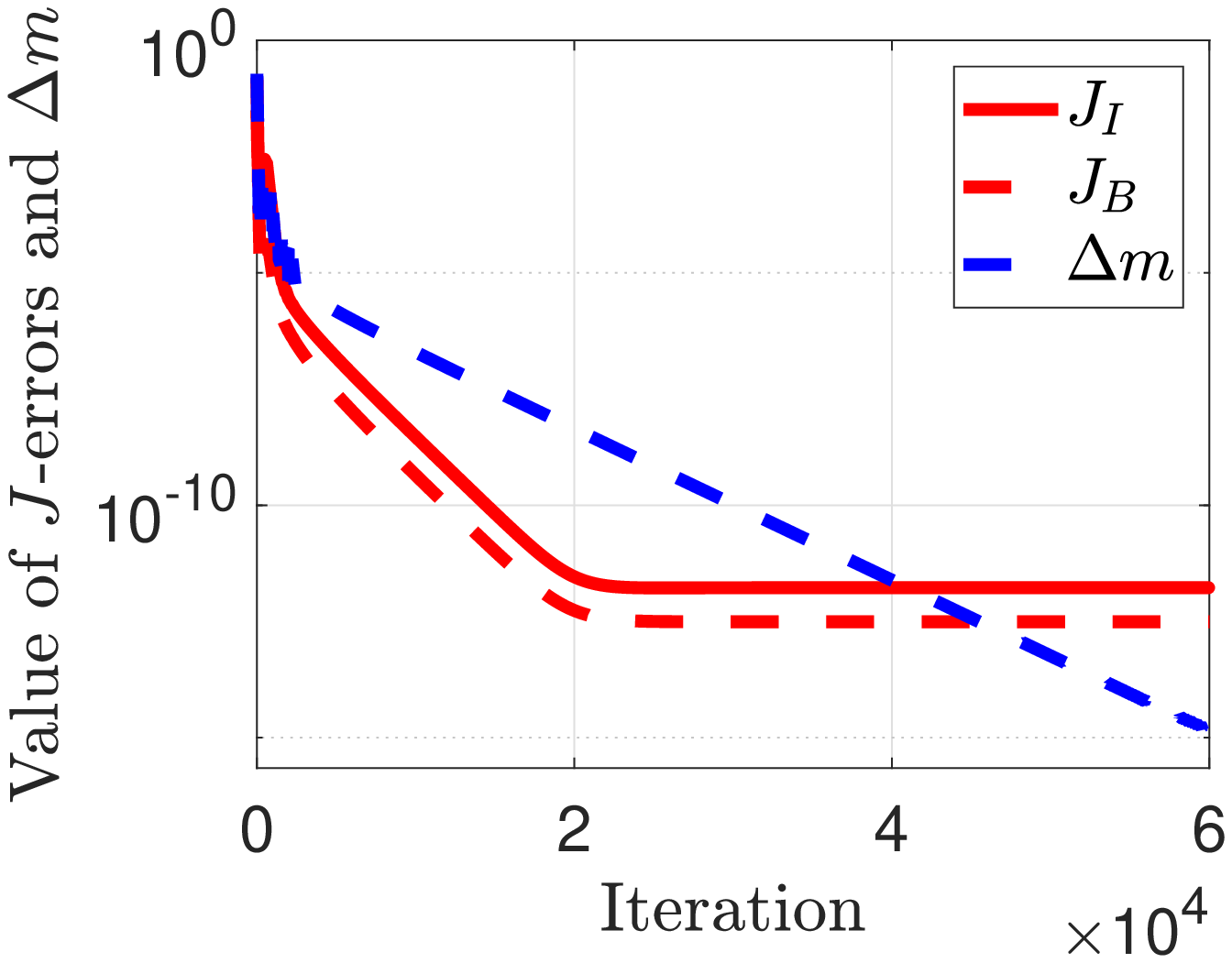"}
	\hspace{-5pt}
	\includegraphics[width = 0.33\linewidth]{"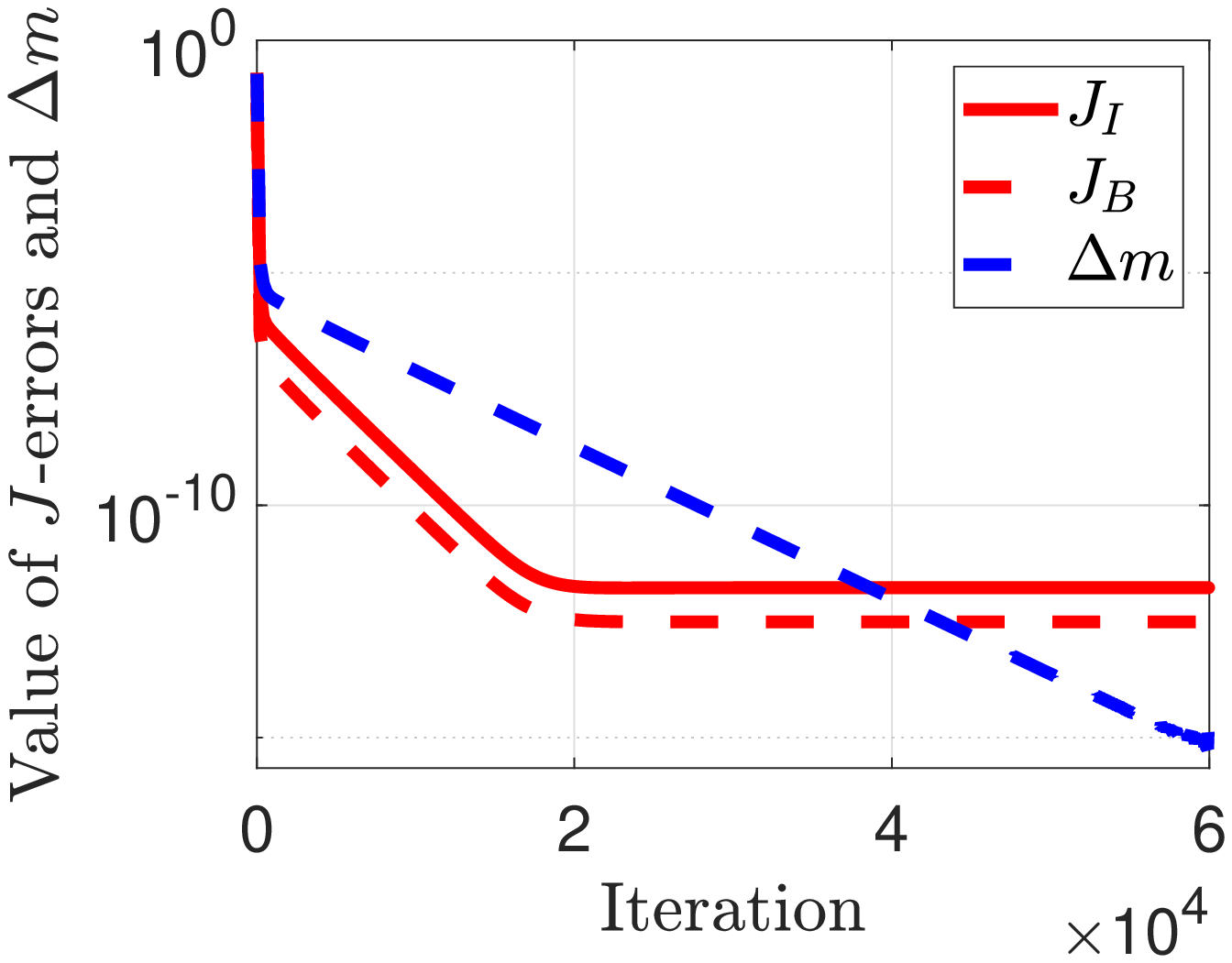"}
	\hspace{-5pt}
	\includegraphics[width = 0.33\linewidth]{"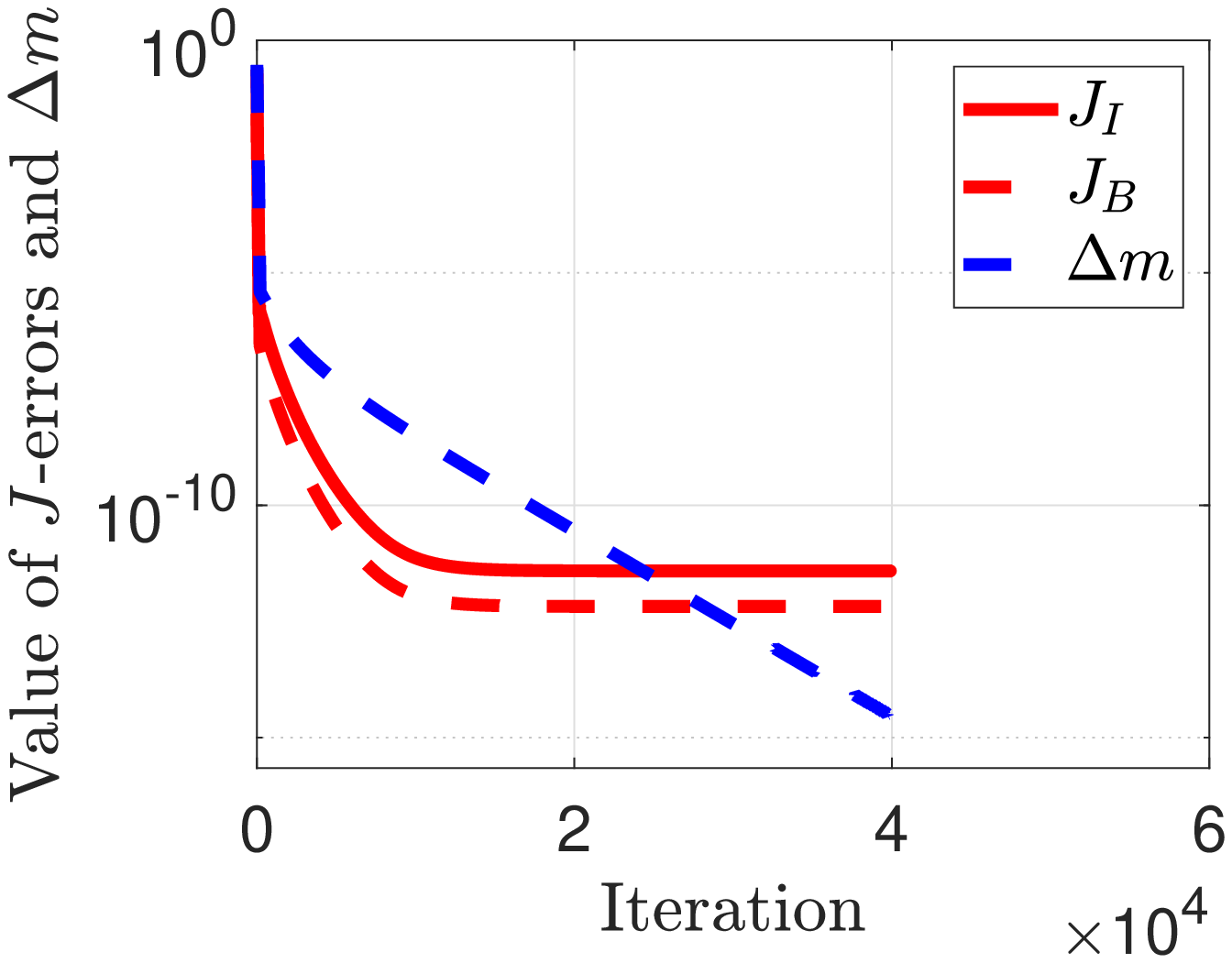"}
	\caption{ $J$-errors and $\Delta m$ over the iterations for PM (left), SPM (middle) and SALM (right). }
	\label{fig:case_5_jUpdates}
\end{figure}

\begin{figure}[t!]
\centering
\includegraphics[width = 0.45\linewidth]{"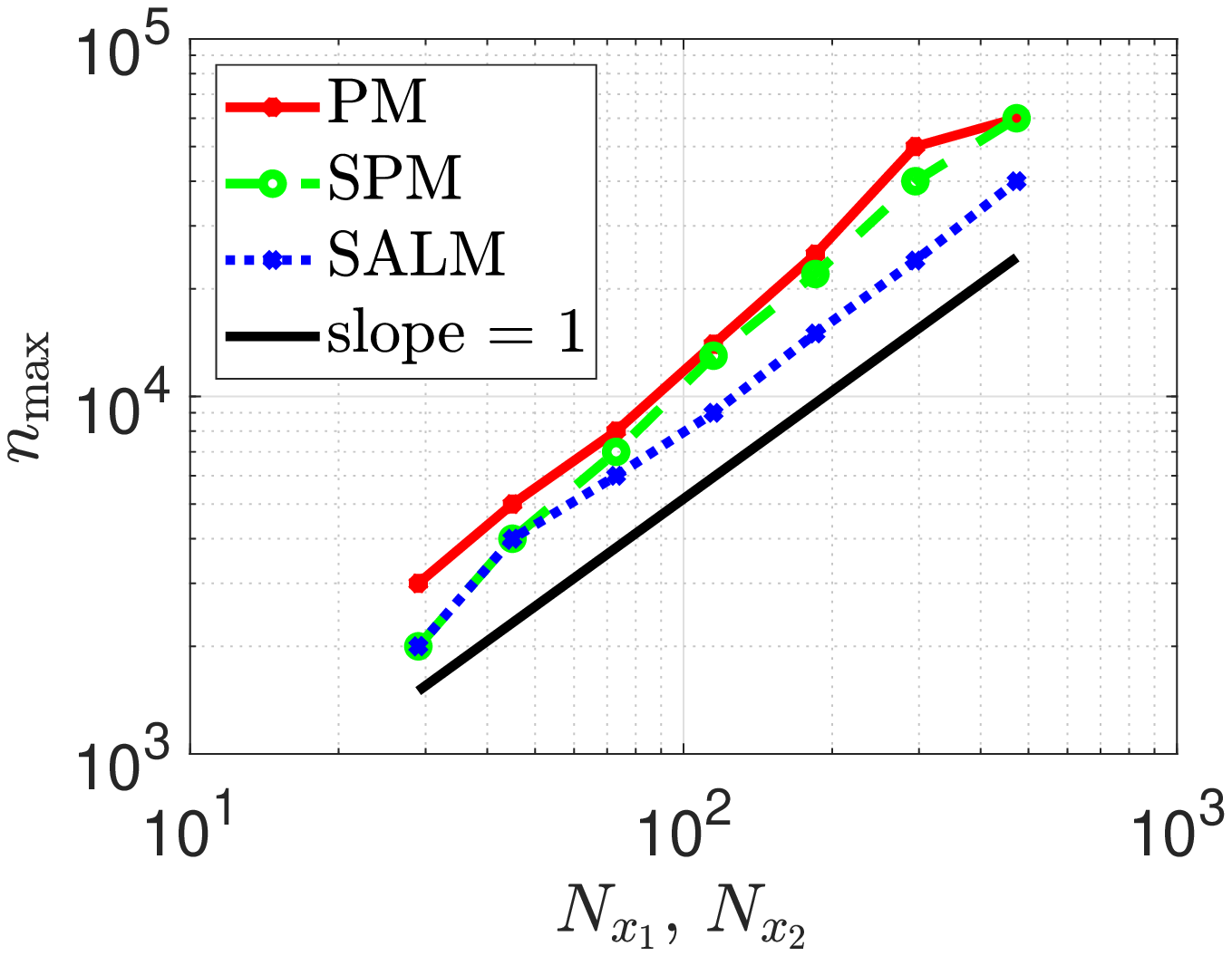"}
\qquad
\includegraphics[width = 0.45\linewidth]{"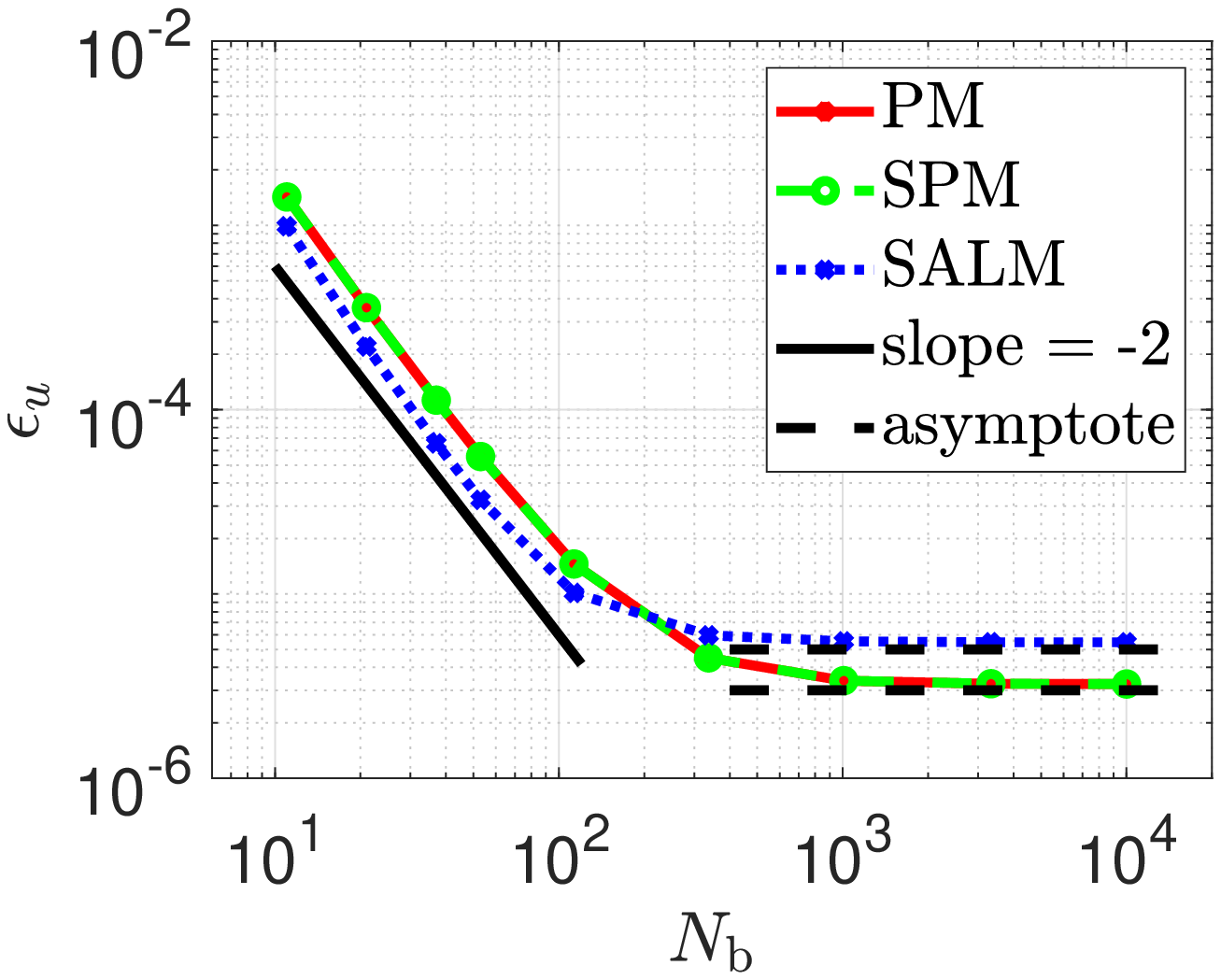"}
\caption{ The total number of iterations needed for convergence (left) and the influence of the boundary discretization on $\epsilon_u$ (right). }
\label{fig:case_5_nMax__timingnBConv}
\end{figure}

For the first example we consider $\mathcal{X} = [0, 1] \times[-1/2, 1/2]$, $\partial \mathcal{Y} = \cup_{k=1}^4 \Gamma_k^\mathcal{Y}$ with
\begin{subequations}
	\label{eqn:case_5_boundarySegments}
	\begin{align}
		\Gamma_1^\mathcal{Y}(s) & = (\cos(\tfrac{1}{2} - s), \, \sin(\tfrac{1}{2} - s)), \\
		\Gamma_2^\mathcal{Y}(s) & = (e^s \cos(\tfrac{1}{2}), \, -e^s \sin(\tfrac{1}{2})), \\
		\Gamma_3^\mathcal{Y}(s) & = ( e \cos(\tfrac{1}{2} - s), \, e\sin(s - \tfrac{1}{2})), \\
		\Gamma_4^\mathcal{Y}(s) & = (e^{1-s} \cos(\tfrac{1}{2}), \, e^{1-s} \sin(\tfrac{1}{2})),
	\end{align}
\end{subequations}	
\begin{figure}[b!]
		\centering
		\includegraphics[width=0.48\linewidth]{"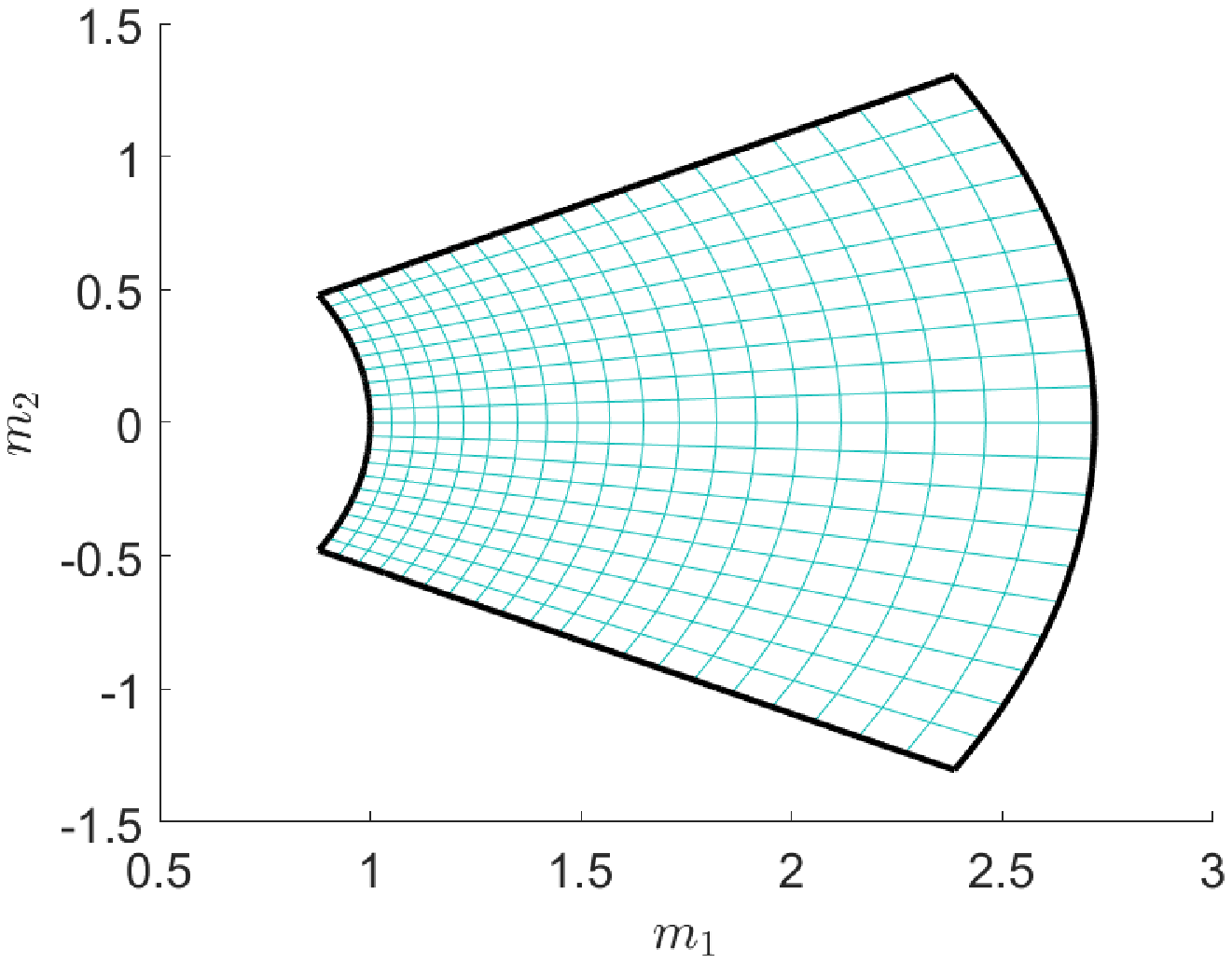"}
		\quad
		\includegraphics[width=0.48\linewidth]{"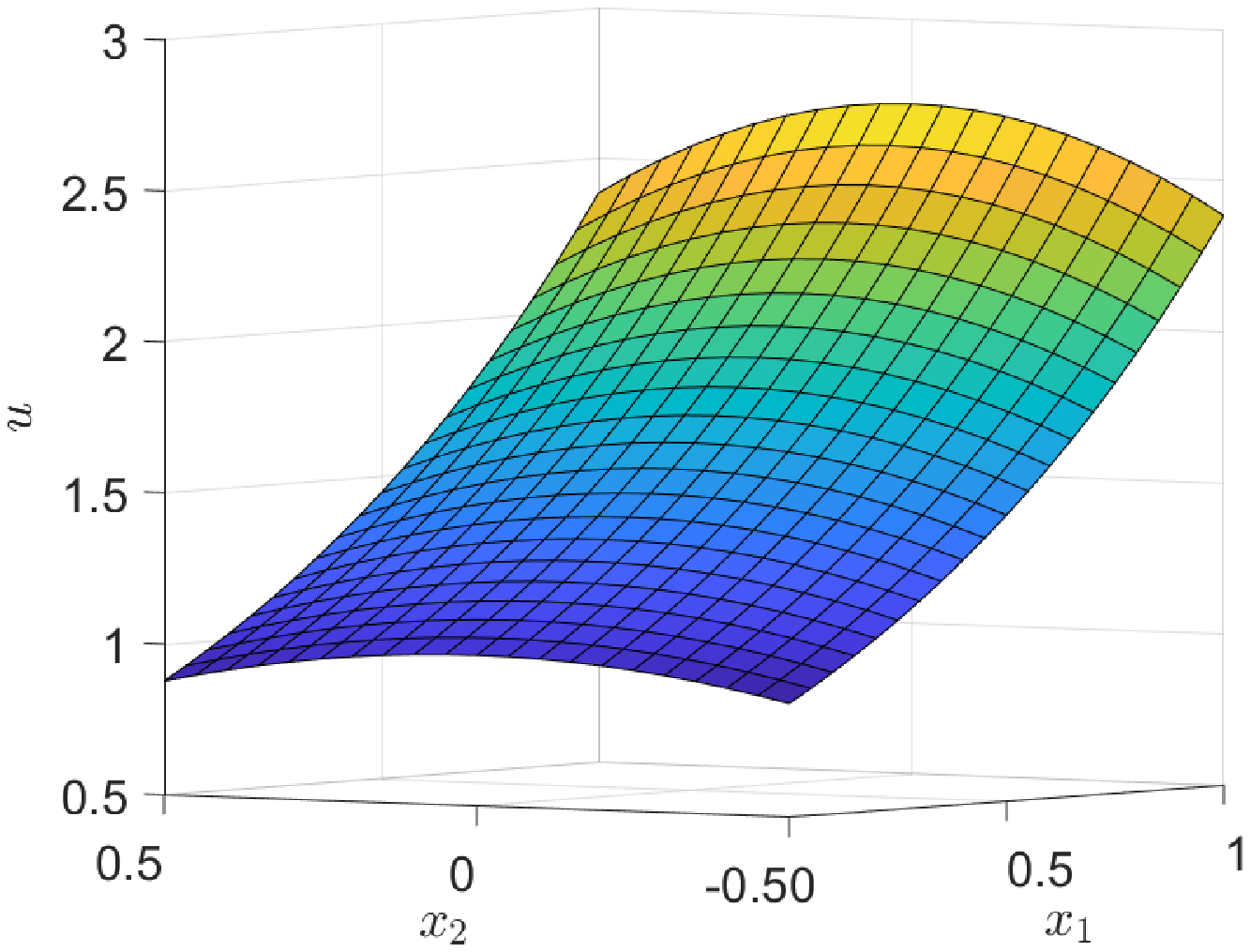"}
		\caption{ The exact mapping $\mathbf{m}$ (left) and solution $u$ (right) on a $21\times21$ grid.}
		\label{fig:Exact_case_5}
\end{figure}%
\noindent{}as shown together with the exact mapping on a $21\times 21$ grid in Figure~\ref{fig:Exact_case_5} on the left. We choose $\Gamma^\mathcal{Y}_k = \nabla u(\Gamma^\mathcal{X}_k)$ for all examples. Furthermore, this choice of $\Gamma_k^\mathcal{Y}$ implies that for SALM and SPM the corners of $\partial \mathcal{X}$ are mapped to the corners of $\partial \mathcal{Y}$. Let $f^2(x_1,x_2) = e^{2x_1}$. The solution is then given by 
\begin{align}
	u(x_1,x_2) = e^{x_1} \cos(x_2),
\end{align}
which is symmetric in $x_2 = 0$ as can be seen in Figure~\ref{fig:Exact_case_5}. Unless specified otherwise, we take $N_\mathrm{b} = 10^4$ and for each target segment $\Gamma^\mathcal{Y}_k$, with $k = 1, \dots, 4$, we construct $\mathbf{y}_i = \Gamma^\mathcal{Y}_k(s_i)$ with $s_i = (i-1) / (N_\mathrm{b}-1)$ and $i=1,\dots, N_\mathrm{b}$. The results for PM, SPM and SALM are shown in Figure~\ref{fig:Conv_results_case_5} for varying grid configurations with $N_{x_1} = N_{x_2}$. The three figures clearly show second-order convergence of the relevant errors and residual, which is in accordance with the discretization error of the finite difference approximations used. In terms of $\epsilon_u$, PM and SPM ($3\cdot 10^{-6}$) slightly outperform SALM ($6\cdot 10^{-6}$) for $N_{x_1} = N_{x_2} = 473$, though the difference is small. 
Figure~\ref{fig:case_5_jUpdates} shows the $J$-errors over the iterations on a grid of $N_{x_1} = N_{x_2} = 473$. For PM and SPM we obtained $J_\mathrm{I} \approx 2\cdot 10^{-12}$, $J_\textrm{B} \approx 3 \cdot 10^{-13}$ in approximately 60,000 iterations. SALM gave $J_\mathrm{I} \approx 4\cdot 10^{-12}$, $J_\textrm{B} \approx 7 \cdot 10^{-13}$ in 40,000, iterations. SALM consistently requires less iterations as seen on the left in Figure~\ref{fig:case_5_nMax__timingnBConv}, where the number of required iterations ($n_\mathrm{max}$) for various $N_{x_1} = N_{x_2}$ is shown.

Convergence of $u$ with respect to $N_\mathrm{b}$ is shown on the right of Figure~\ref{fig:case_5_nMax__timingnBConv} for $N_{x_1} = N_{x_2} = 473$. Two observations are in place. First, for increasing $N_\mathrm{b}$, the error $\epsilon_u$ reaches an asymptotic value (dashed black line). This phenomenon is to be expected and occurs when the discretization errors in $\mathbf{m}_{ij}$, $\mathbf{P}_{ij}$ and $u_{ij}$, and the finite differences $\mathbf{D}_{ij}$ become dominant, i.e., when the discretization errors due to the choice of $N_{x_1}$ and $N_{x_2}$ dominate the errors due to discretizing the boundary. Secondly, in the regime prior to the asymptote, the discretization error in $u$ due to the boundary discretization is second-order accurate for all three boundary methods.

\begin{figure}[t!]
	\centering
	\includegraphics[width = 0.33\linewidth]{"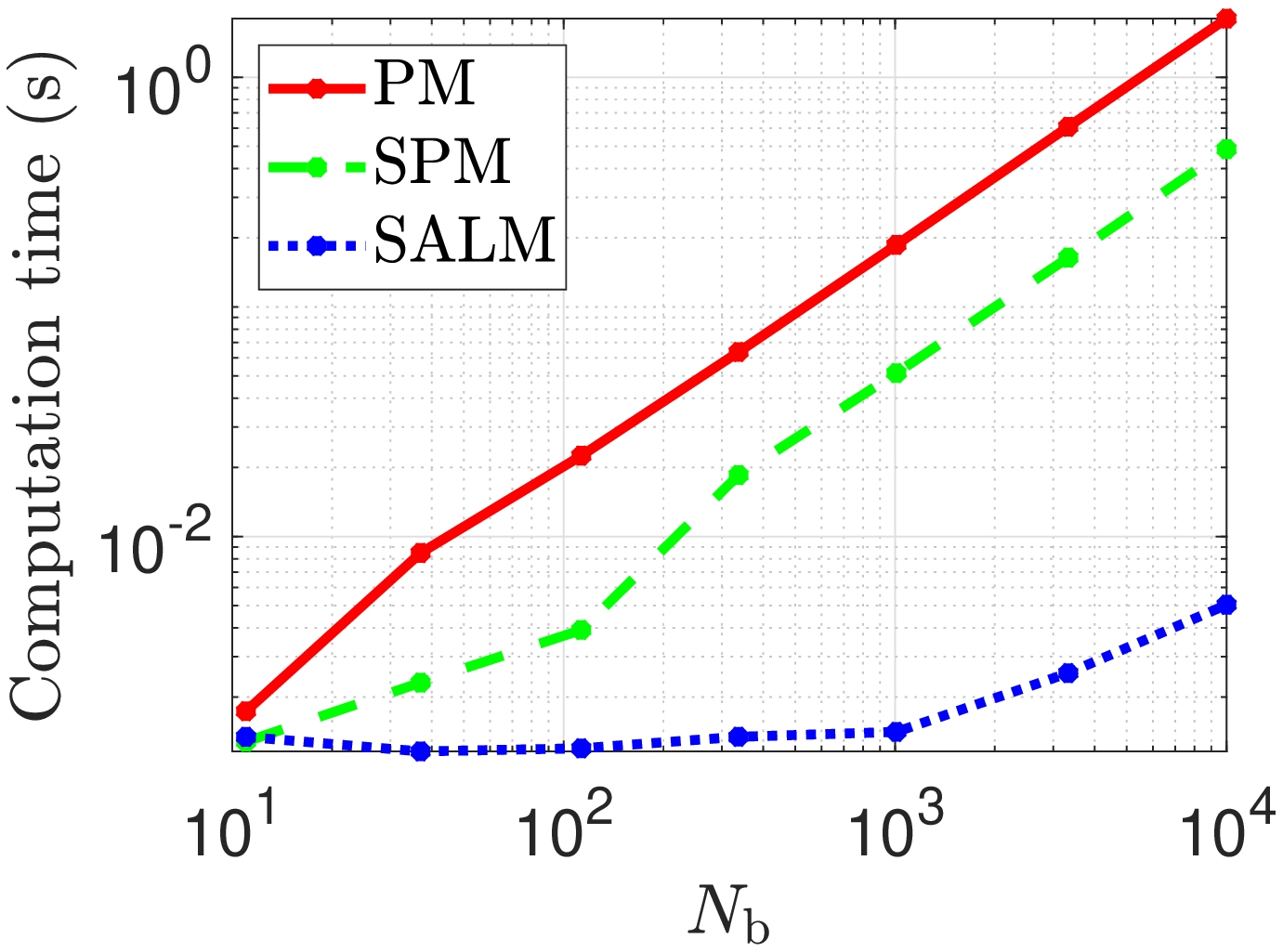"}
	\hspace{-5pt}
	\includegraphics[width = 0.33\linewidth]{"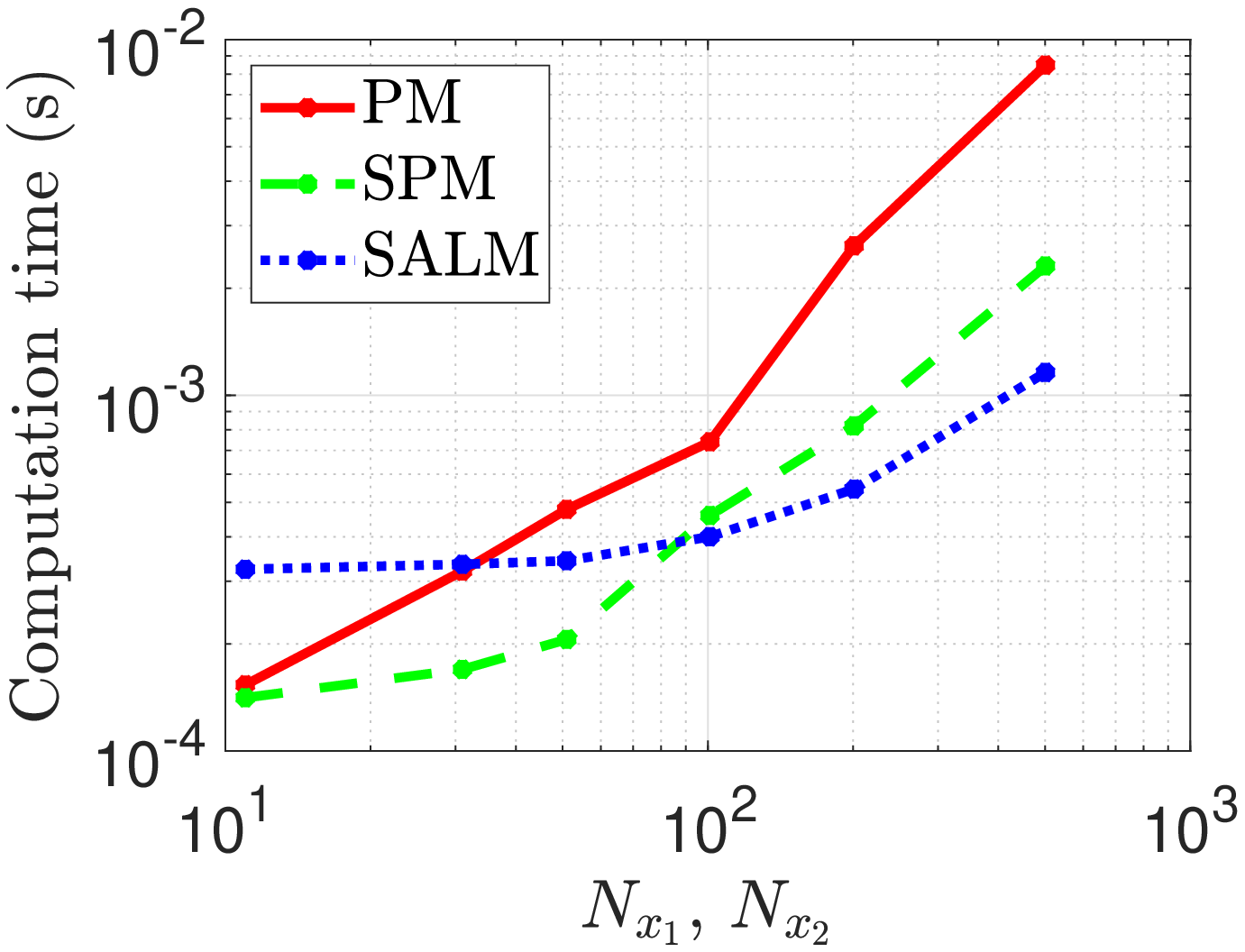"}
	\hspace{-5pt}
	\includegraphics[width = 0.33\linewidth]{"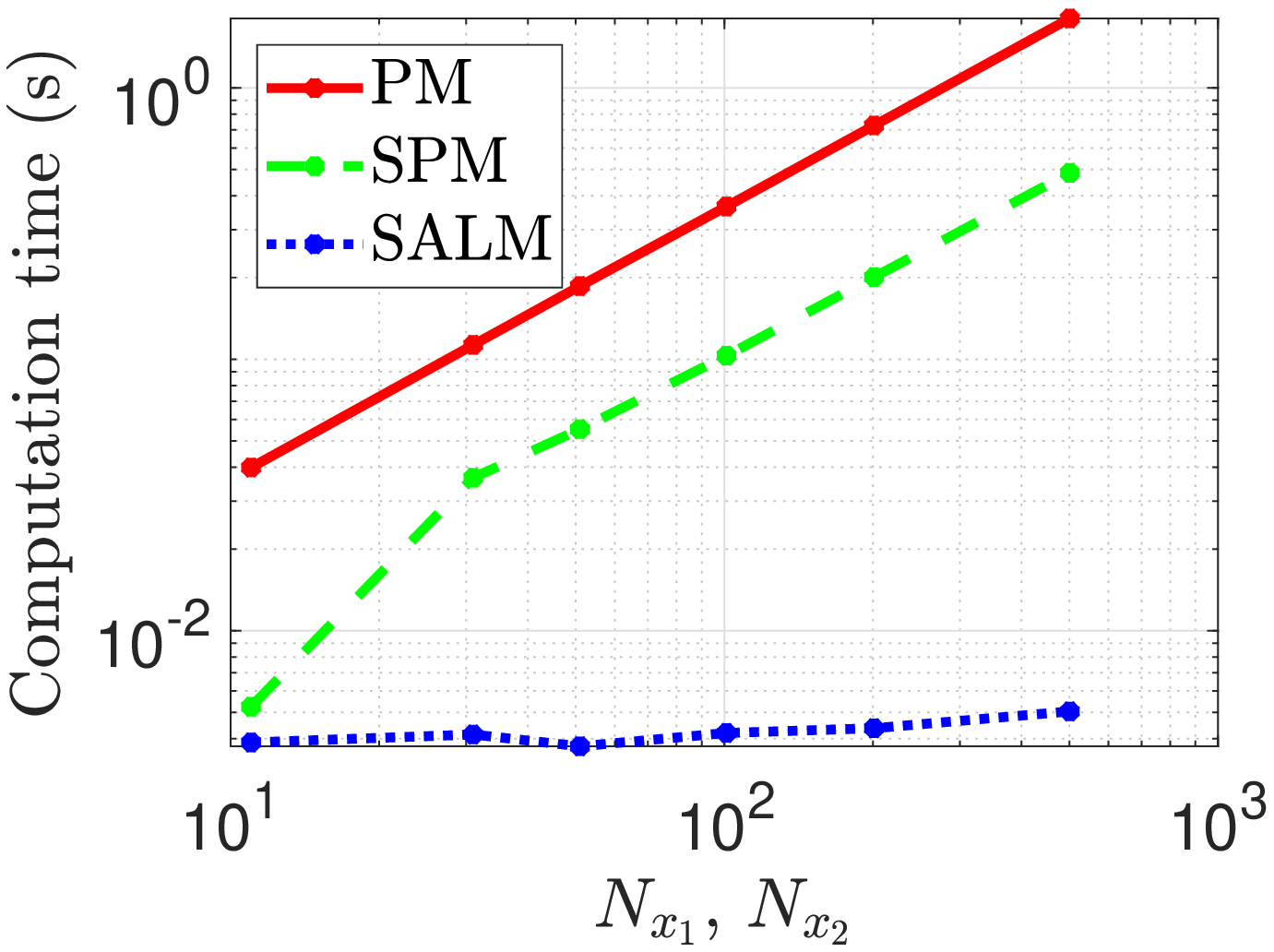"}
	\caption{ Timing results in seconds for the boundary procedures. On the left, $N_{x_1} = N_{x_2} = 501$ is fixed and in the middle and on the right $N_\textrm{b} = 37$ and $N_\textrm{b} = 10007$, respectively.}
	\label{fig:timing}
\end{figure}

Finally, the computational cost per iteration for SALM is lowest, second comes SPM and third PM. The projection methods calculate $C \max (N_{x_1}, N_{x_2})  N_\mathrm{b}$ projections and performs $C \max (N_{x_1}, N_{x_2} )$ searches over $N_\mathrm{b}$ points each, with $C \in \mathbb{N}_+$. Similarly, SALM performs a linear interpolation of $C \max (N_{x_1}, N_{x_2})$ points over $N_\mathrm{b}-1$ segments. 
Therefore, one would expect the average time per iteration to scale linearly in $N_x, N_y$ and $N_\mathrm{b}$ for PM and SPM when $N_{x_1} = N_{x_2}$, and $N_\mathrm{b}$ is fixed. For SALM, a linear relation is also expected, with a possible asymptote when either the computational load due to $\max (N_{x_1}, N_{x_2})$ or $N_\mathrm{b}$ dominates. This is also shown in Figure~\ref{fig:timing}, where, from left to right, $N_{x_1} = N_{x_2} = 501$ is fixed while $N_\textrm{b}$ varies, $N_\textrm{b} = 37$ is fixed and $N_{x_1} = N_{x_2}$ varies, and lastly, $N_\textrm{b} = 10007$ is fixed and $N_{x_1} = N_{x_2}$ varies. Additionally, it is observed that SALM, on average, significantly outperforms the projection methods. Lastly, SPM is approximately four times faster than PM because SPM projects one source segment on a target segment (four times) instead of the whole source boundary on the whole target boundary.

\FloatBarrier

\subsection{Deformed square}
\label{sec:DeformedSquare}
\begin{figure}[t]
   	\centering
   	\includegraphics[width=0.48\linewidth]{"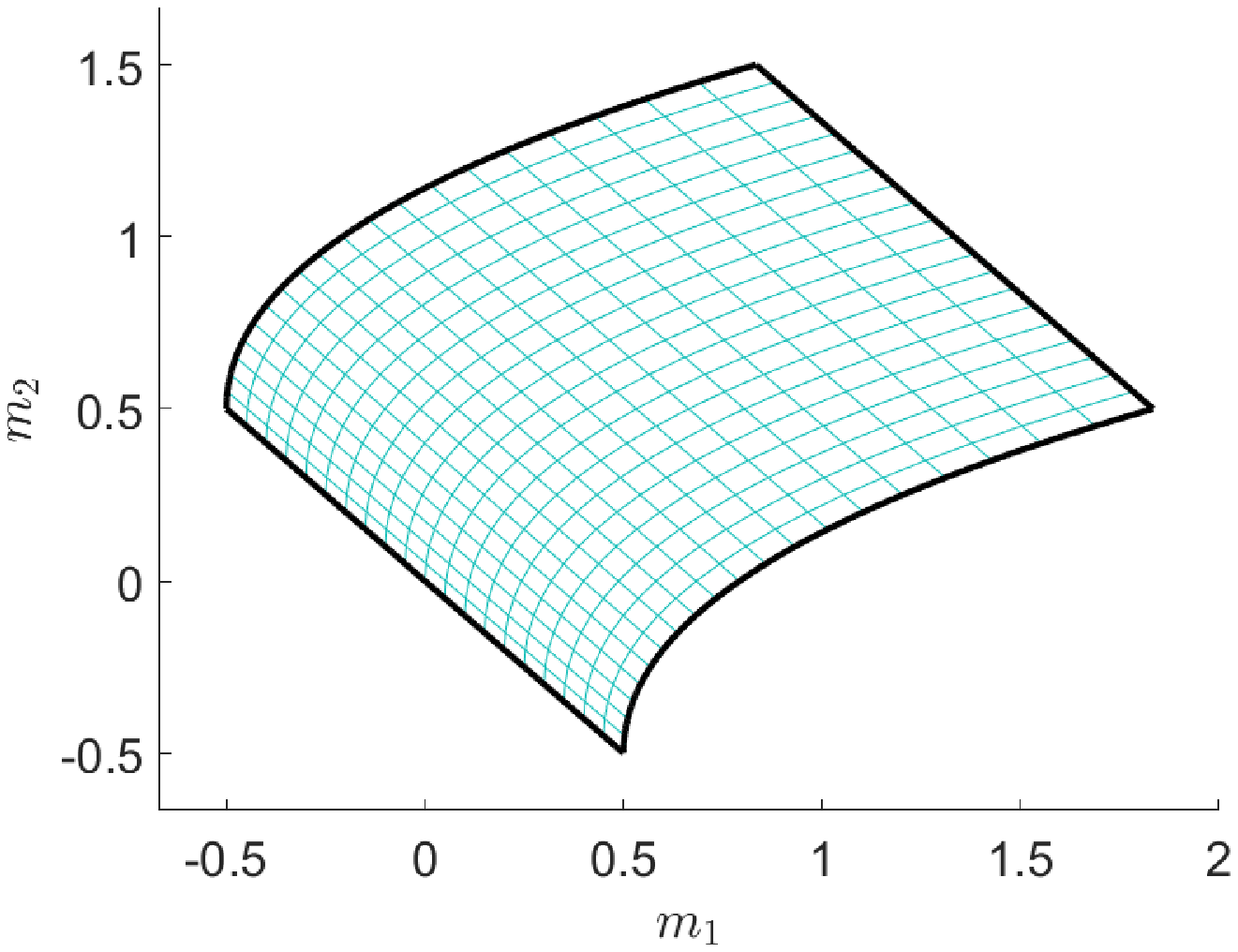"}
   	\quad
   	\includegraphics[width=0.48\linewidth]{"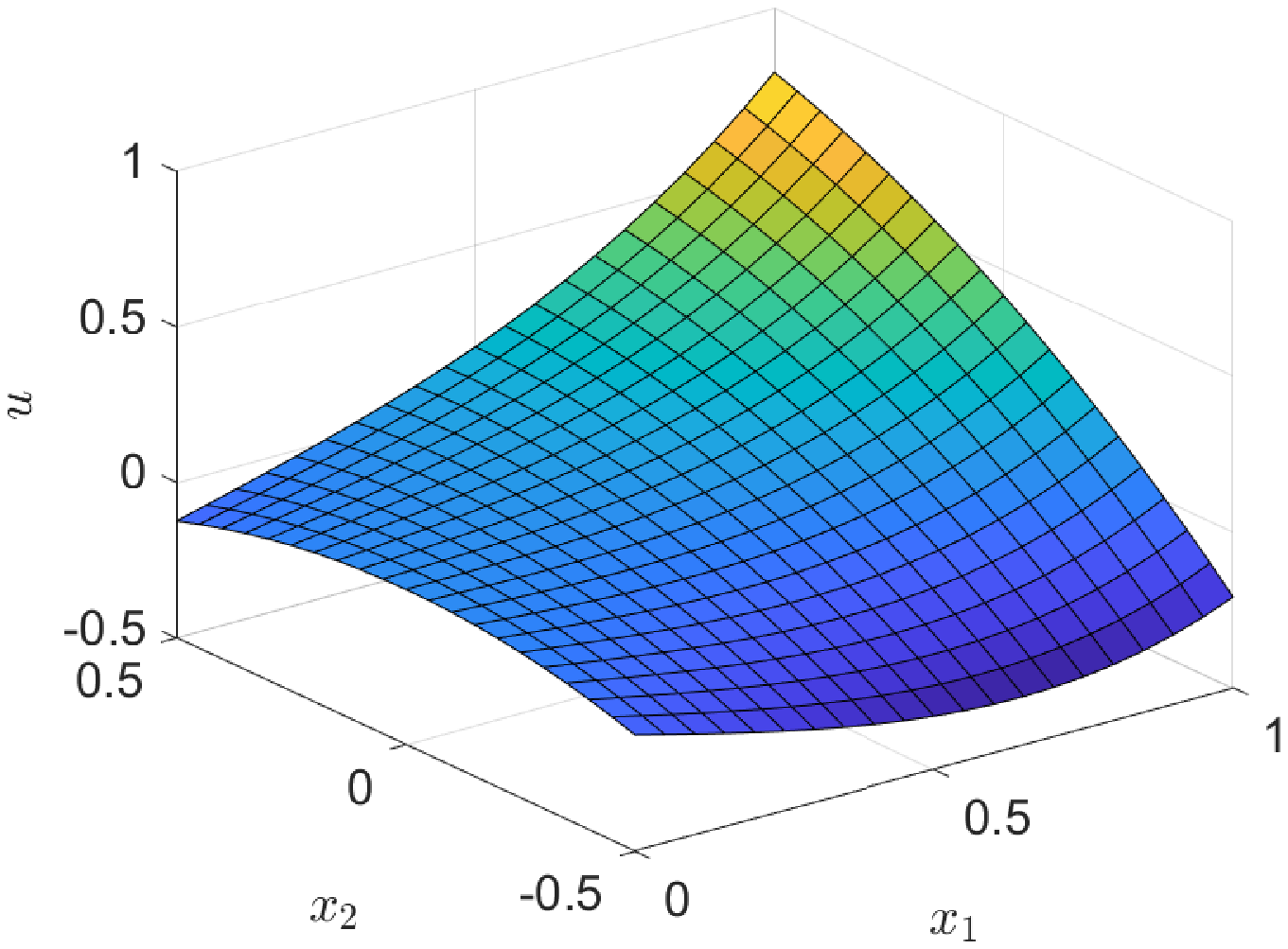"}
   	\caption{ The exact mapping $\mathbf{m}$ (left) and solution $u$ (right) on a $21\times21$ grid.}
   	\label{fig:Exact_case_3}
\end{figure}

\begin{figure}[t]
    \centering
    \includegraphics[width = 0.48\linewidth]{"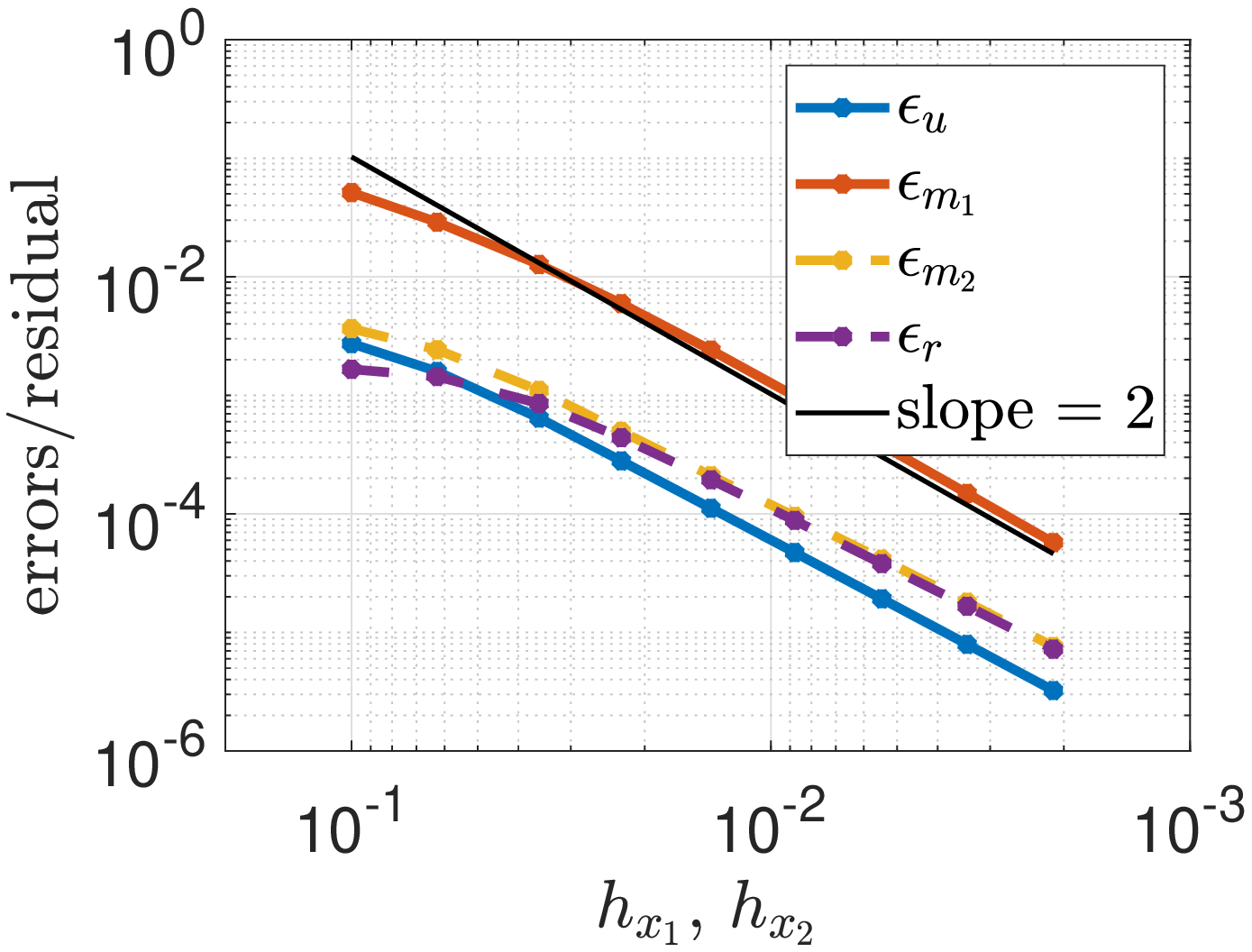"}
	\quad
    \includegraphics[width = 0.48\linewidth]{"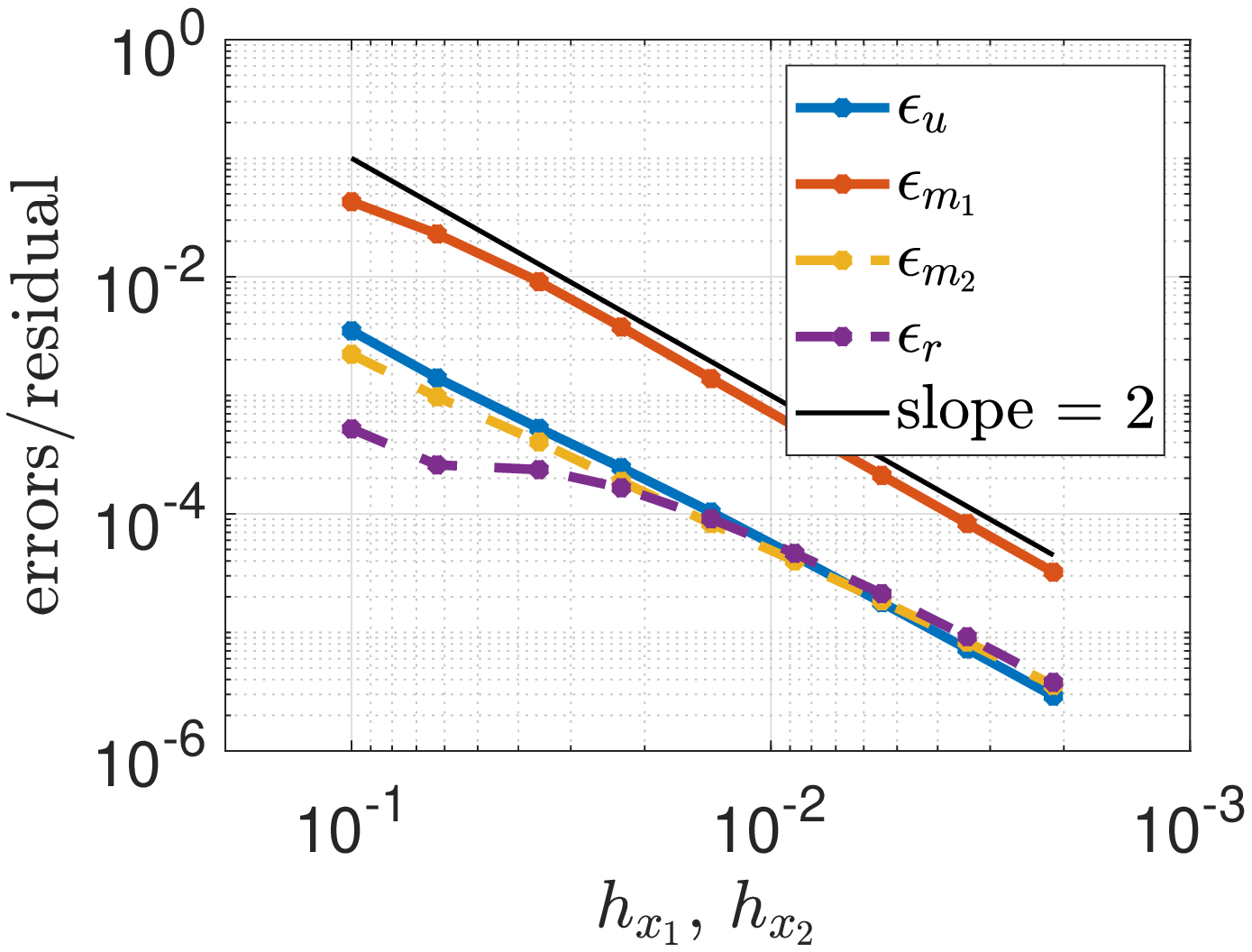"}
    \caption{ Convergence of the global error and the residual for SPM (left) and SALM (right). }
    \label{fig:Conv_results_case_3}
\end{figure}

\begin{figure}[b!]
    \centering
    \includegraphics[width = 0.48\linewidth]{"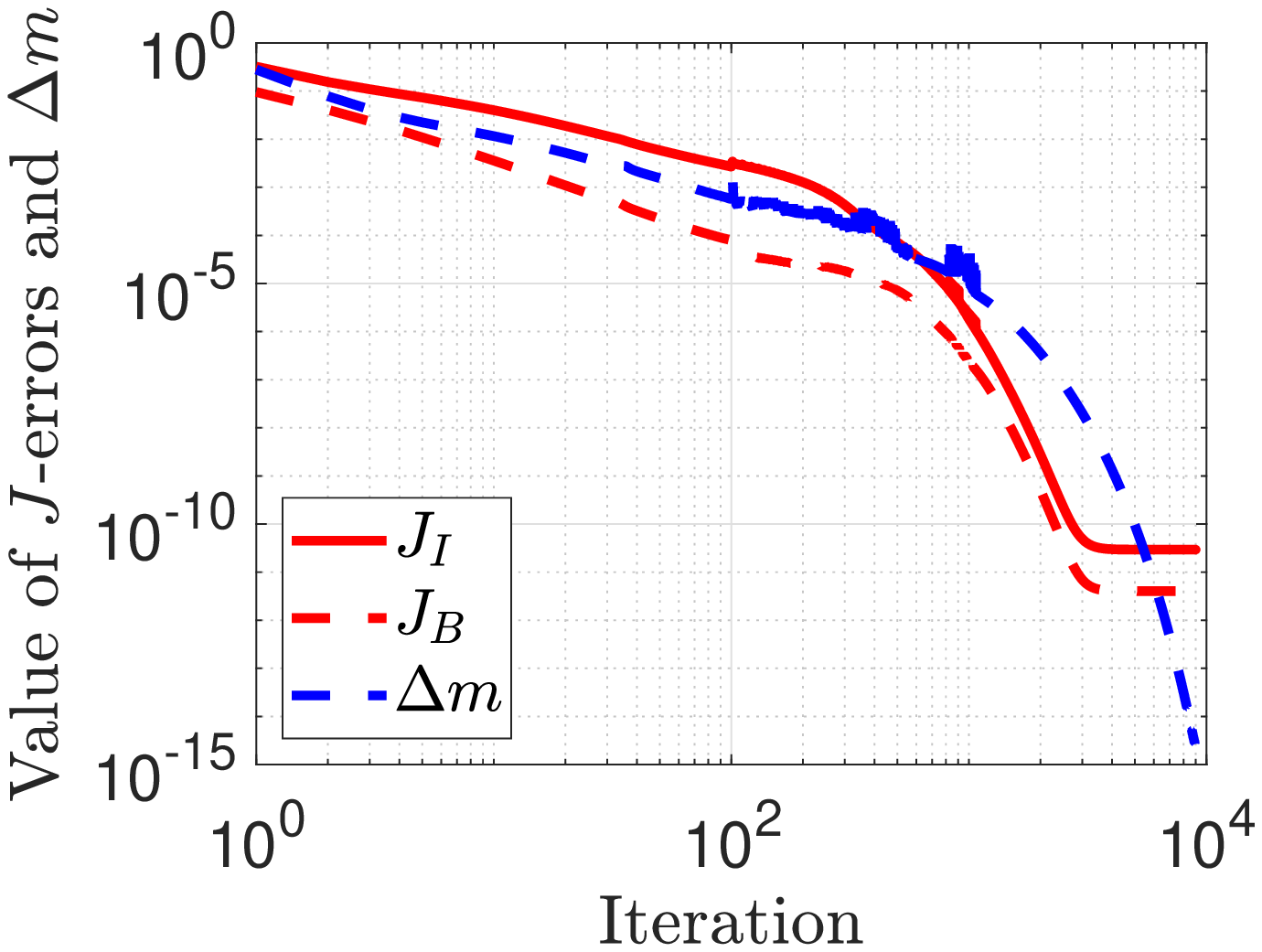"}
	\quad
    \includegraphics[width = 0.48\linewidth]{"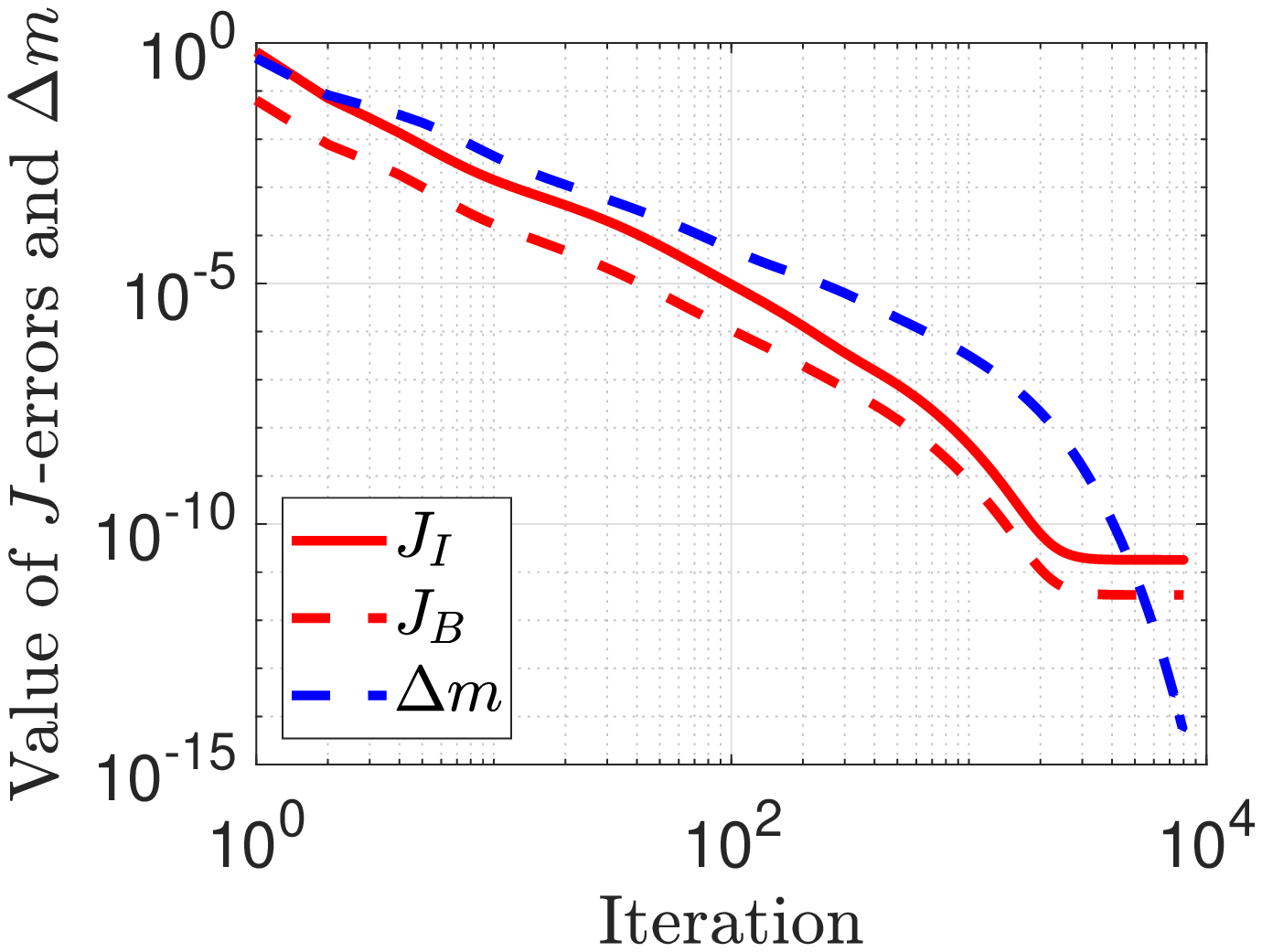"}
    \caption{ Comparison of $J_\textrm{I}$, $J_\textrm{B}$ and $\Delta m$ for SPM (left) and SALM (right) for a grid of $295 \times 295$. }
    \label{fig:JError_case_3}
\end{figure}

As a second example we consider a mapping and surface with no symmetries, viz. $\mathcal{X} = [0, 1] \times[-1/2, 1/2]$, $\partial \mathcal{Y} = \cup_{k=1}^4 \Gamma_k^\mathcal{Y}$ with
\begin{subequations}
\begin{align}
\Gamma_1^\mathcal{Y}(s) & =  (s - \tfrac{1}{2}, \, -s + \tfrac{1}{2}), \\
\Gamma_2^\mathcal{Y}(s) & =  (\tfrac{1}{3}s^3 + s^2 + \tfrac{1}{2}, \, s - \tfrac{1}{2}) , \\
\Gamma_3^\mathcal{Y}(s) & =  (-s + \tfrac{11}{6}, \, s + \tfrac{1}{2}), \\
\Gamma_4^\mathcal{Y}(s) & =  (- \tfrac{1}{3}s^3 + 2 s^2 - 3 s + \tfrac{5}{6}, \, \tfrac{3}{2} - s),
\end{align}
\end{subequations}
as shown on the left of Figure~\ref{fig:Exact_case_3}. Let $f^2(x_1,x_2) = (x_1+1)^2$, then the exact solution is given by
\begin{align}
u(x_1,x_2) & = \frac{1}{12} x_1^4  + \frac{1}{3} x_1^3 + x_1 x_2 - \frac{1}{2} x_2^2,
\end{align}
which is shown on the right of Figure~\ref{fig:Exact_case_3}.
In Figure~\ref{fig:Conv_results_case_3} errors and residuals are shown for SPM (left) and SALM (right) for varying grid configurations with $N_{x_1} = N_{x_2}$.
Both figures show second-order convergence which is in accordance with the discretization error of the finite difference approximations used.
Figure~\ref{fig:JError_case_3} shows the behaviour of $J_\textrm{I}$, $J_\textrm{B}$ and $\Delta m$ for SPM and SALM. For SPM, on the left $\Delta m$ exhibit oscillations starting at 100 iterations. This is due to the grid shock correction, enabled in the 100$^\text{th}$ iteration. As it turns out, this is one example where grid shocks occur using SPM. Without the grid shock correction, $\Delta m$ would still go to computer precision but the grid shock (as visualized on the right of Figure~\ref{fig:crossingGridline}) would remain and subsequently $J_\textrm{I}$, $J_\textrm{B}$ and the errors $\epsilon_u$, $\epsilon_{m_1}$ and $\epsilon_{m_2}$ and the residual $\epsilon_r$ would be three orders of magnitude higher.

\begin{figure}[b!]
    \centering
    \includegraphics[width = 0.48\linewidth]{"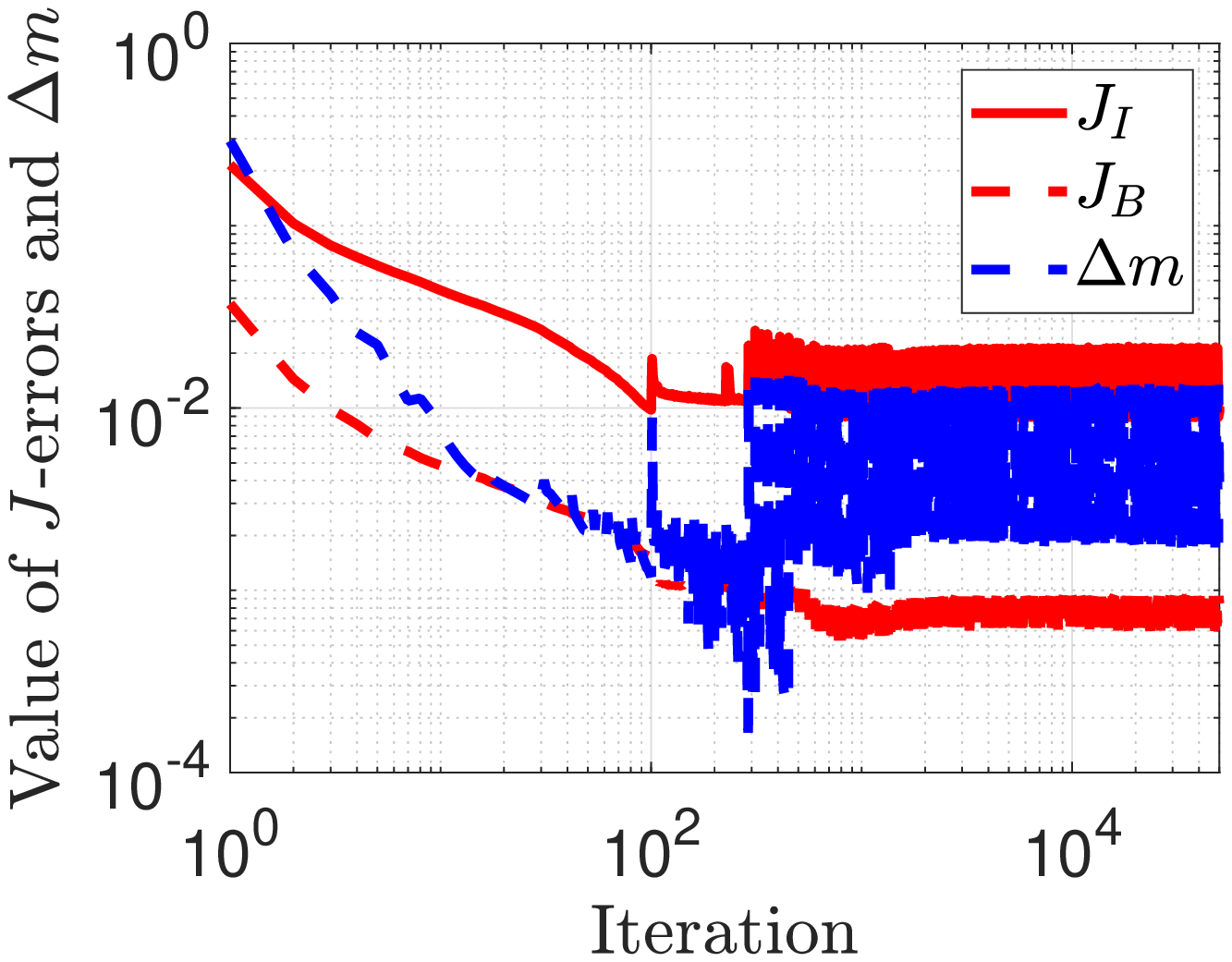"}
    \quad
    \includegraphics[width = 0.48\linewidth]{"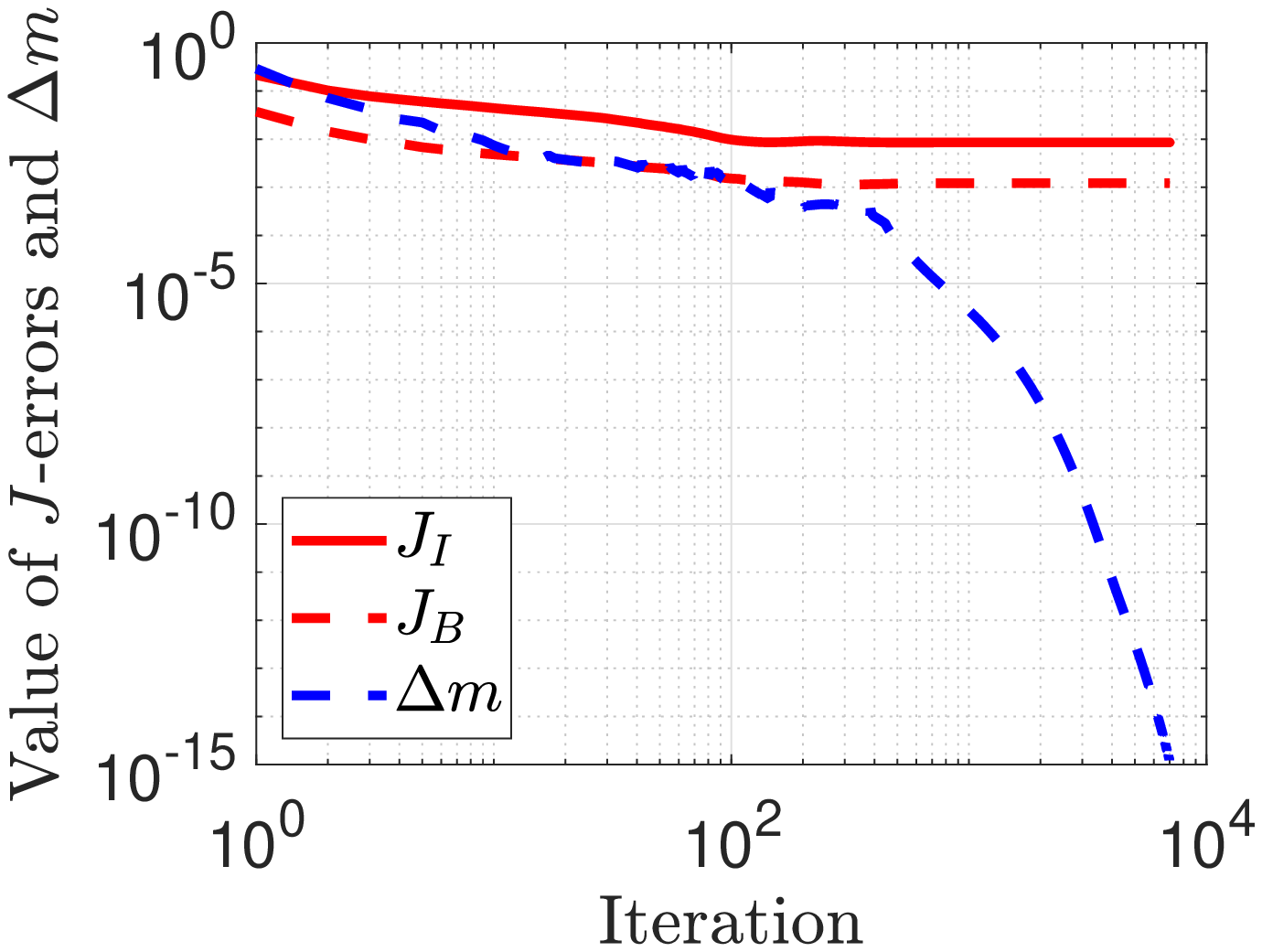"}
    \caption{ The errors $J_\textrm{I}$, $J_\textrm{B}$ and the update $\Delta m$ for PM with shock correction (left) and without (right). }
    \label{fig:JBerrors_case_3_PM}
\end{figure}

\begin{figure}[htbp]
	\centering
	\includegraphics[width = 0.345\linewidth]{"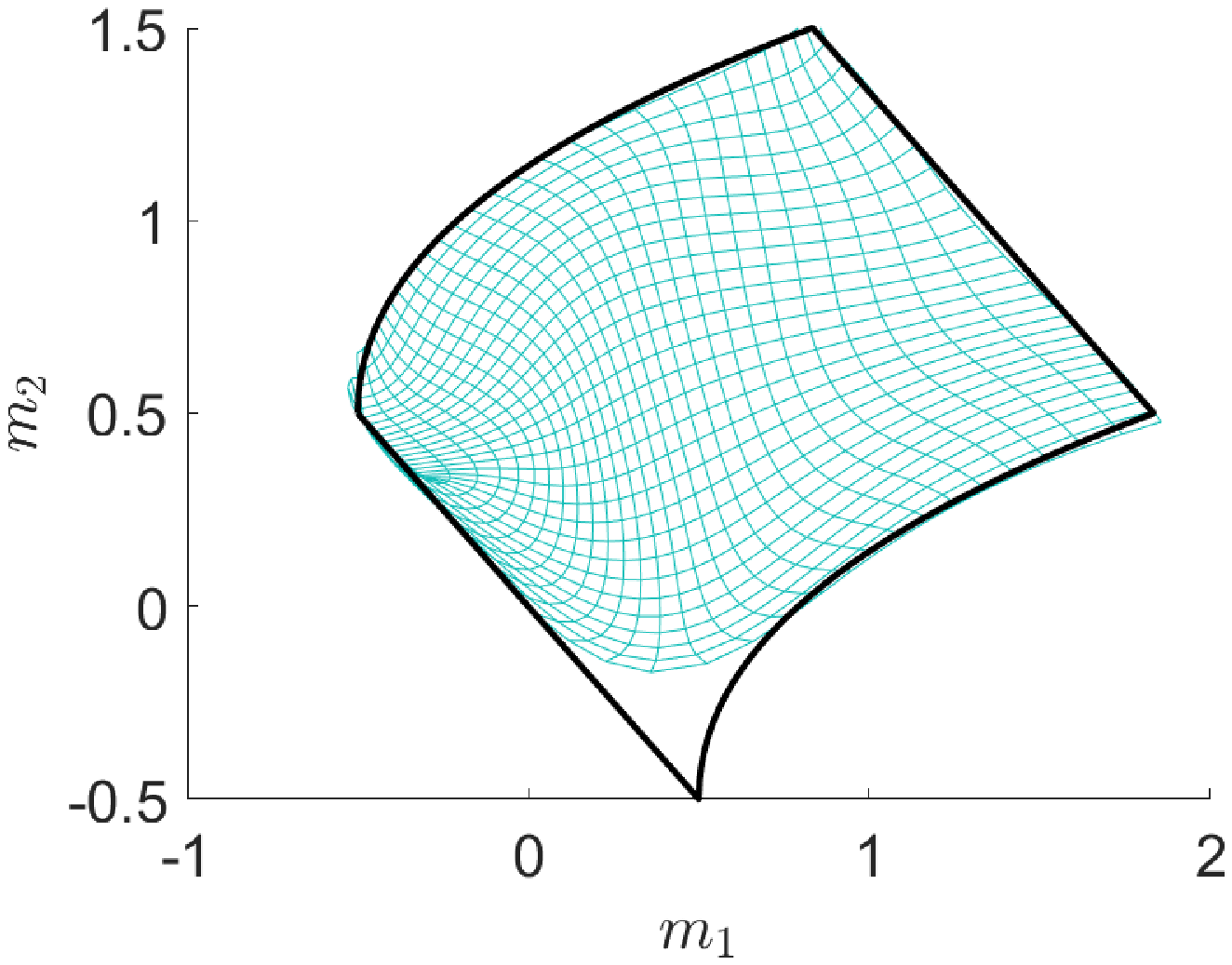"}
	\hspace{-15pt}
	\includegraphics[width = 0.345\linewidth]{"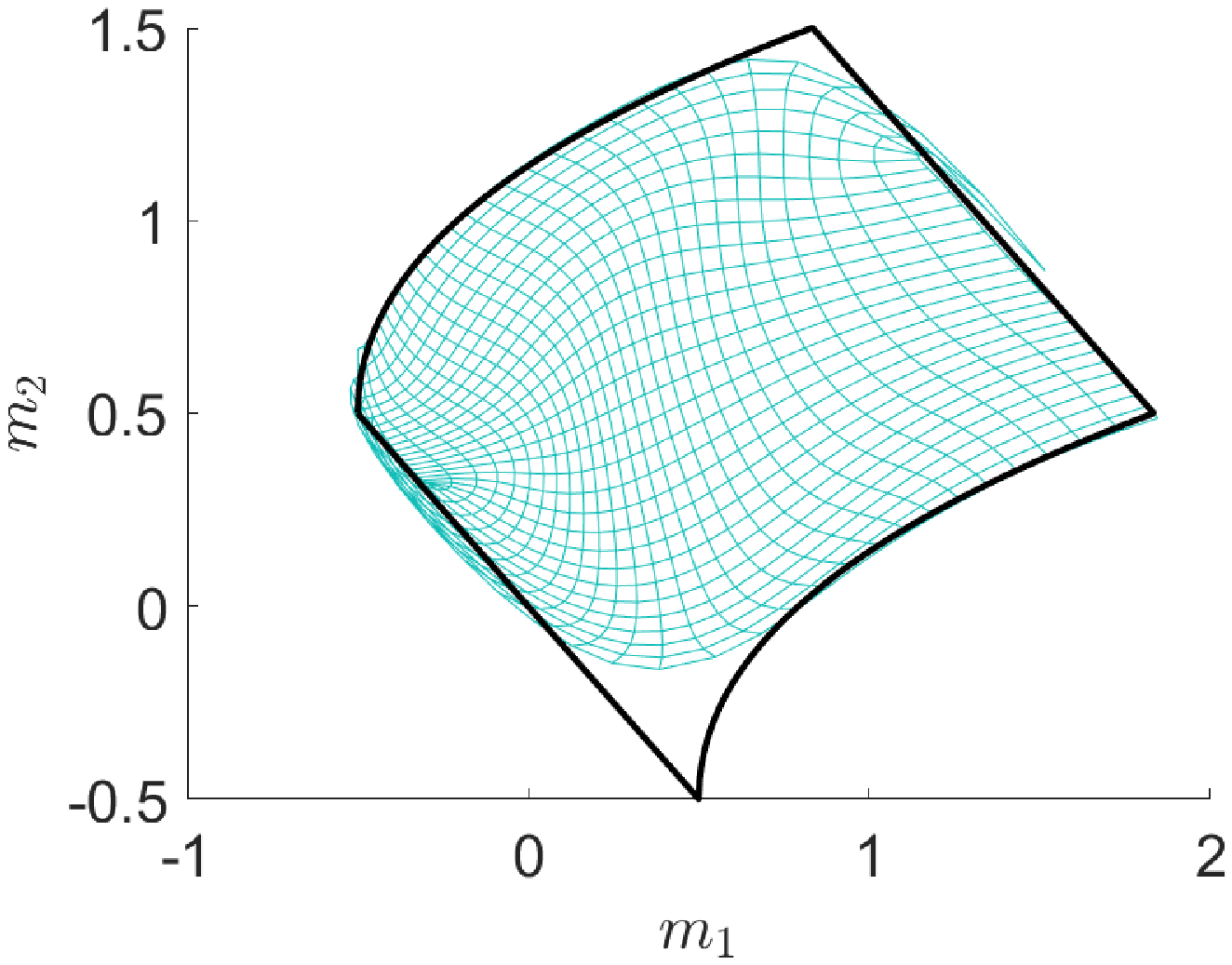"}
	\hspace{-15pt}
	\includegraphics[width = 0.345\linewidth]{"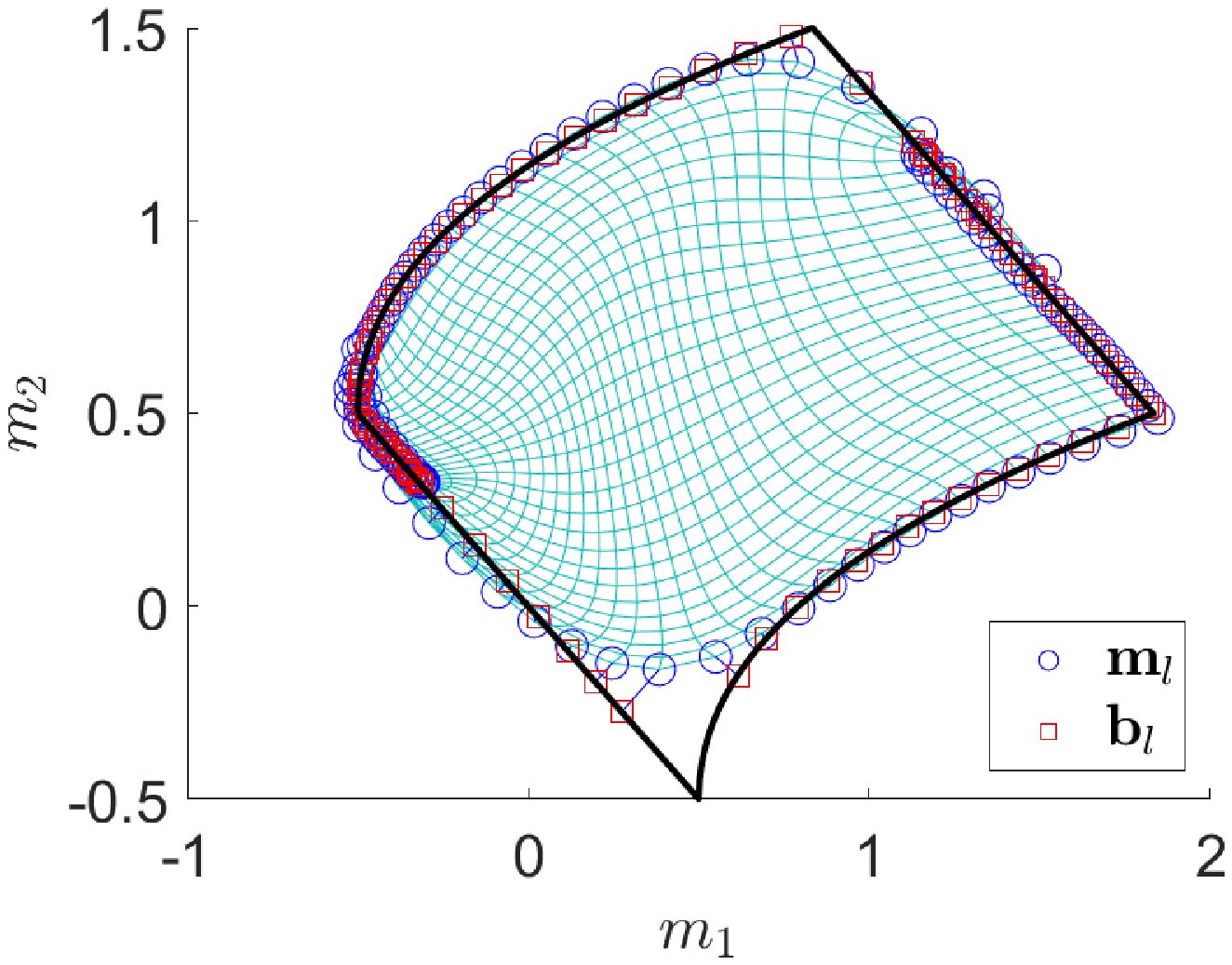"}
	\caption{ The mapping after 50,000 iterations for PM with shock correction (left) and without (middle) and the accompanying projection of $\mathbf{m}_\elll$ onto $\partial \mathcal{Y}$ for construction of $\mathbf{b}_\elll$ without shock correction on the right. }
	\label{fig:Mappings_case_3_PM}
\end{figure}

Lastly we discuss the results for PM. In Figure~\ref{fig:JBerrors_case_3_PM} the errors $J_\textrm{I}$ and $J_\textrm{B}$ are shown, both with grid shock correction (left) and without (right). Clearly the example with grid shock correction does not converge, the method actually oscillates between intermediate solutions. One may be tempted to think that without grid shock correction the method does work, as $\Delta m$ goes to machine precision, but this is not the case as shown in Figure~\ref{fig:Mappings_case_3_PM}. The two leftmost figures show the mapping for a $29\times29$ grid after 50,000 iterations. Clearly, neither of the methods work as intended as there are gaps between the mesh spanned by $\mathbf{m}_\elll$ and $\partial \mathcal{Y}$, i.e., the transport boundary condition has not been satisfied. The reason why the algorithm with PM does not converge is that the projection of $\mathbf{m}$ onto $\partial \mathcal{Y}$ does not distribute $\mathbf{b}$ well. In particular, no points $\mathbf{b}_\elll$ near $(0.5, -0.5) \in \partial \mathcal{Y}$ are obtained, as can be seen in Figure~\ref{fig:Mappings_case_3_PM} on the right, where the blue circles represent $\mathbf{m}_\elll$, the red squares $\mathbf{b}_\elll$ and the thin black lines connect $\mathbf{m}_\elll$ to $\mathbf{b}_\elll$ for $\elll = 1, \dots, N$.

\FloatBarrier
    
\subsection{Inward fold}
\begin{figure}[b!]
   	\includegraphics[width=0.45\linewidth]{"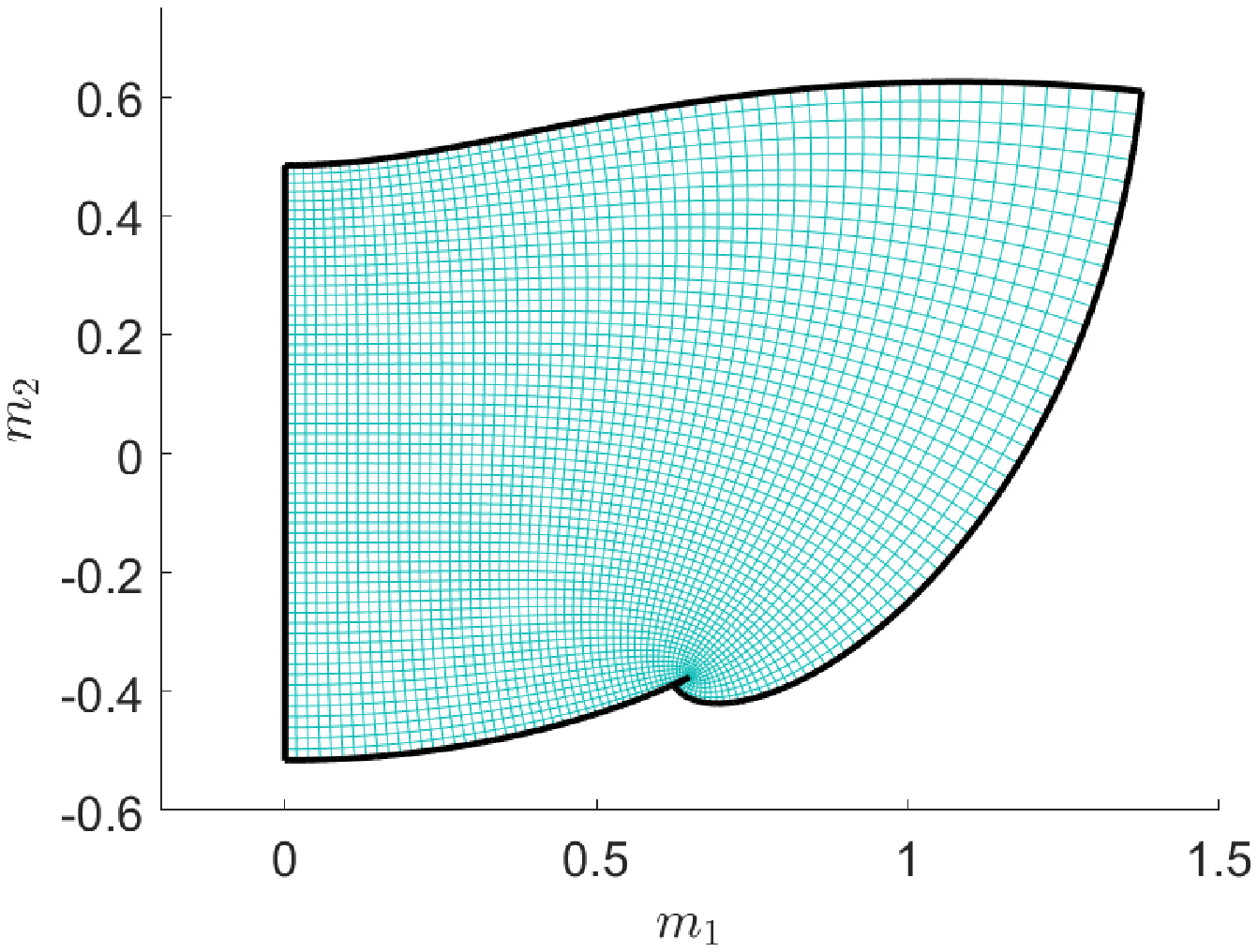"}
   	\quad
	\includegraphics[width=0.45\linewidth]{"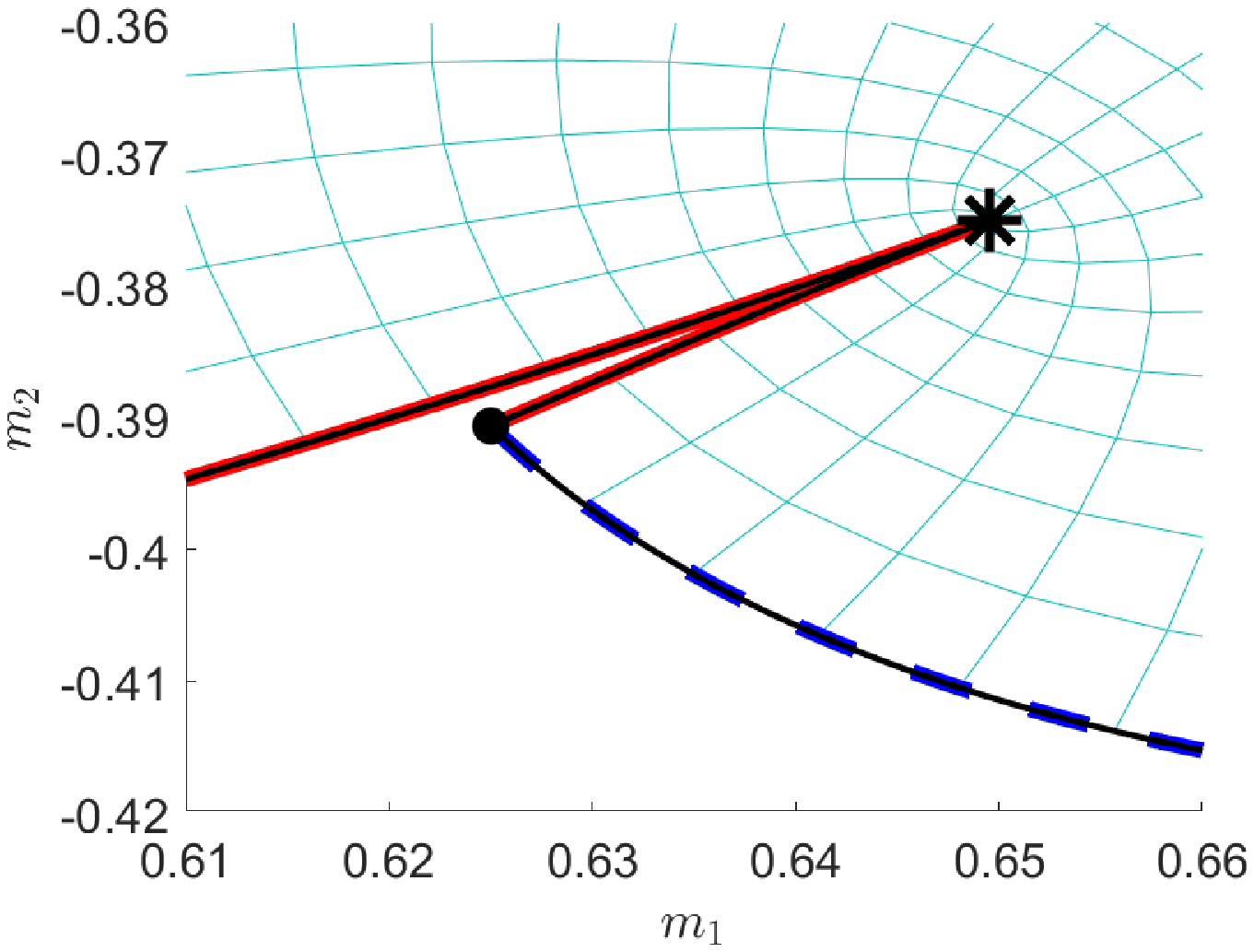"}
	\caption{The target domain and mapping, with a zoomed-in version of the Inward fold example.}
	\label{fig:case_11_mapping_and_zoomed}
\end{figure}

For this example we consider the target as illustrated in Figure~\ref{fig:case_11_mapping_and_zoomed}, for the exact solution on a $61\times61$ grid. On the right a zoomed-in version of the target is shown. In the figure we have marked two points, one by a solid circle, and one by an asterisk. The former is a point for which the boundary of the target is not differentiable, while for the latter it is. We will come back to this. 
The example shown in Figure~\ref{fig:case_11_mapping_and_zoomed} corresponds to $\mathcal{X} = [0, 1] \times[-1/2, 1/2]$ and $\partial \mathcal{Y} = \cup_{k=1}^4 \Gamma_k^\mathcal{Y}$ with
\begin{subequations}
\begin{align}
\Gamma_1^\mathcal{Y}(s) & = (0,-\tfrac{1}{4}s^4+\tfrac{1}{2}s^3-\tfrac{3 }{8}s^2-\tfrac{7}{8}s+\tfrac{31}{64}), \\
\Gamma_2^\mathcal{Y}(s) & = (\tfrac{9}{8}s-\tfrac{1}{2}s^3,-\tfrac{1}{4}s^4+\tfrac{3}{8}s^2-\tfrac{33}{64}) , \\
\Gamma_3^\mathcal{Y}(s) & = ( -s^3+\tfrac{3 }{2}s^2+\tfrac{1}{4}s+\tfrac{5}{8},-\tfrac{1}{4}s^4+\tfrac{1}{2}s^3+\tfrac{9 }{8}s^2-\tfrac{3 }{8}s-\tfrac{25}{64}) , \\
\Gamma_4^\mathcal{Y}(s) & = (-\tfrac{1}{2}s^3+\tfrac{3 }{2}s^2-\tfrac{19 }{8}s+\tfrac{11}{8},-\tfrac{1}{4}s^4+s^3-\tfrac{9 }{8}s^2+\tfrac{1}{4}s+\tfrac{39}{64}).
\end{align}
\end{subequations}
Furthermore we have
\begin{subequations}
    \begin{align}
    f^2(x_1,x_2) & = x_1^6 + 3 x_1^4 x_2^2 + 3 x_1^2 x_2 (x_2^3 - 2) + (1 + x_2^3)^2, \\
    u(x_1,x_2) & = \frac{x_1^2}{2} - \frac{x_1^4 x_2}{4} - \frac{x_2^2}{2} + \frac{x_1^2 x_2^3}{2} - \frac{x_2^5}{20}.
    \end{align}
\end{subequations}
Taking derivatives of $u$ yields the mapping, i.e.,
\begin{subequations}
\begin{align}
	m_1(x_1, x_2) & = u_{x_1}(x_1, x_2) = x_1 - x_1^3 x_2 + x_1 x_2^3, \\
	m_2(x_1, x_2) & = u_{x_2}(x_1, x_2) = -\frac{x_1^4}{4} - x_2 + \frac{3 x_1^2 x_2^2}{2} - \frac{x_2^4}{4}.
\end{align}
\end{subequations}
A straightforward calculation shows that $\mathbf{m}(1, 1/2) = (5/8, -25/64)$ which corresponds to the solid circle in Figure~\ref{fig:case_11_mapping_and_zoomed}. Henceforth, $\mathbf{m}$ is not differentiable in the point $(1, 1/2)$ as it is the image of a nondifferentiable (corner) point in $\mathcal{X}$ under a continuously differentiable map. The point depicted by the asterisk originates from the source boundary segment $\Gamma^\mathcal{X}_2 = [0,1] \times \{\tfrac{1}{2}\}$. Let $m_1$ and $m_2$ along the boundary be parametrized by $s$. Then in the point indicated by the asterisk, both $\tfrac{\diff m_1(s)}{\diff s}$ and $\tfrac{\diff m_2(s)}{\diff s}$ change sign. Henceforth, the location of the asterisk can be obtained by solving $\tfrac{\diff m_1(s)}{\diff s} = \tfrac{\diff m_2(s)}{\diff s} = 0$, which is equivalent to
\begin{align}
		\left.\pdiff{m_1}{x_1}\right\vert_{x_2 = 1/2} = 0, \qquad
		\left.\pdiff{m_2}{x_1}\right\vert_{x_2 = \tfrac{1}{2}} = 0, \qquad 
		0 \leq x_1 \leq 1.
\end{align}
Indeed, doing so one uniquely obtains $x_1 = \sqrt{3} / 2$ such that \linebreak$\mathbf{m}(\sqrt{3} / 2, 1/2) = (3 \sqrt{3} / 8, - 3/8)$, which corresponds to the point indicated by an asterisk in Figure~\ref{fig:case_11_mapping_and_zoomed}. Furthermore, smoothness of the boundary in said point is implied.

\begin{figure}[t!]
     \centering
     \includegraphics[width = 0.44\linewidth]{"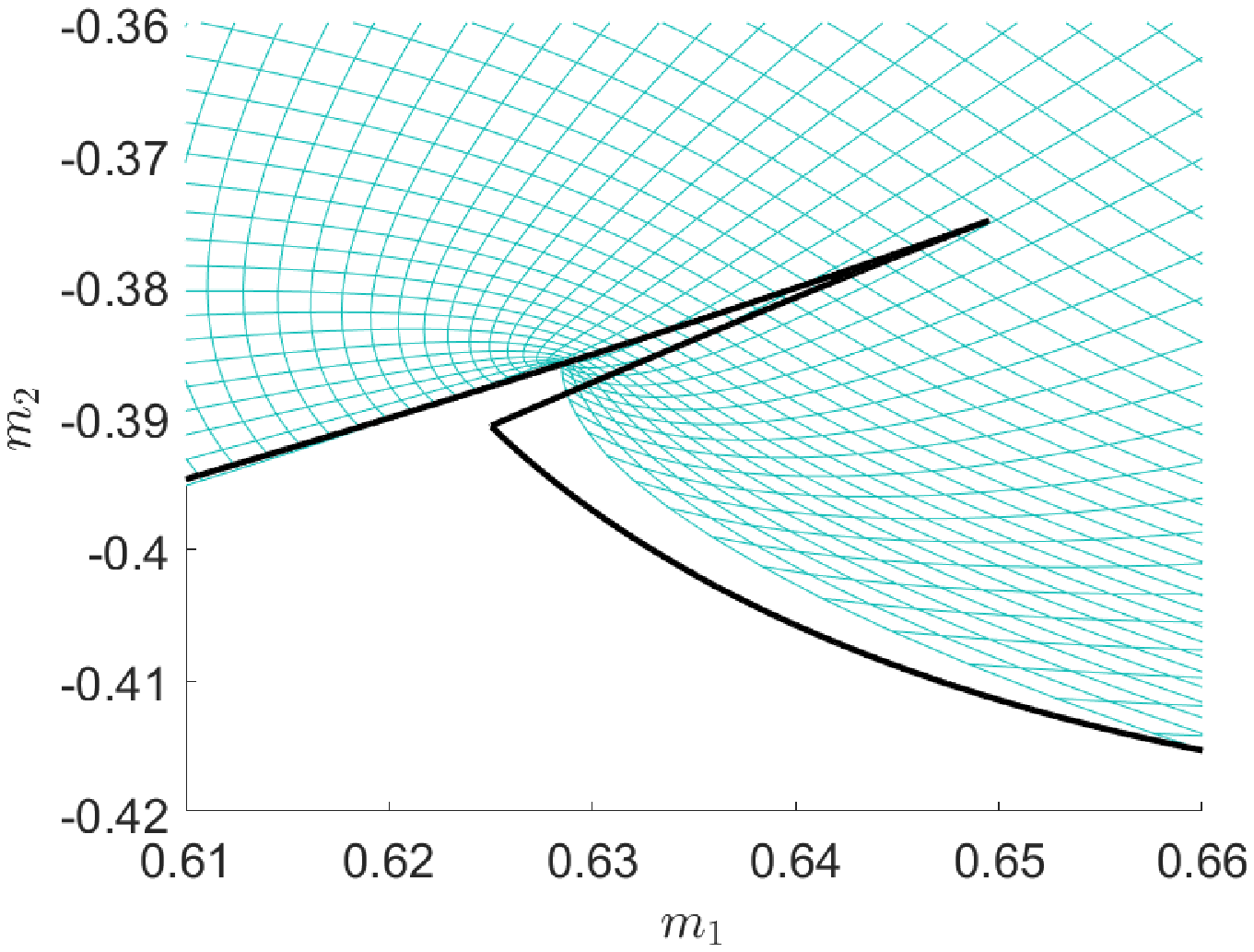"}
  \hspace{-5pt}
     \includegraphics[width = 0.44\linewidth]{"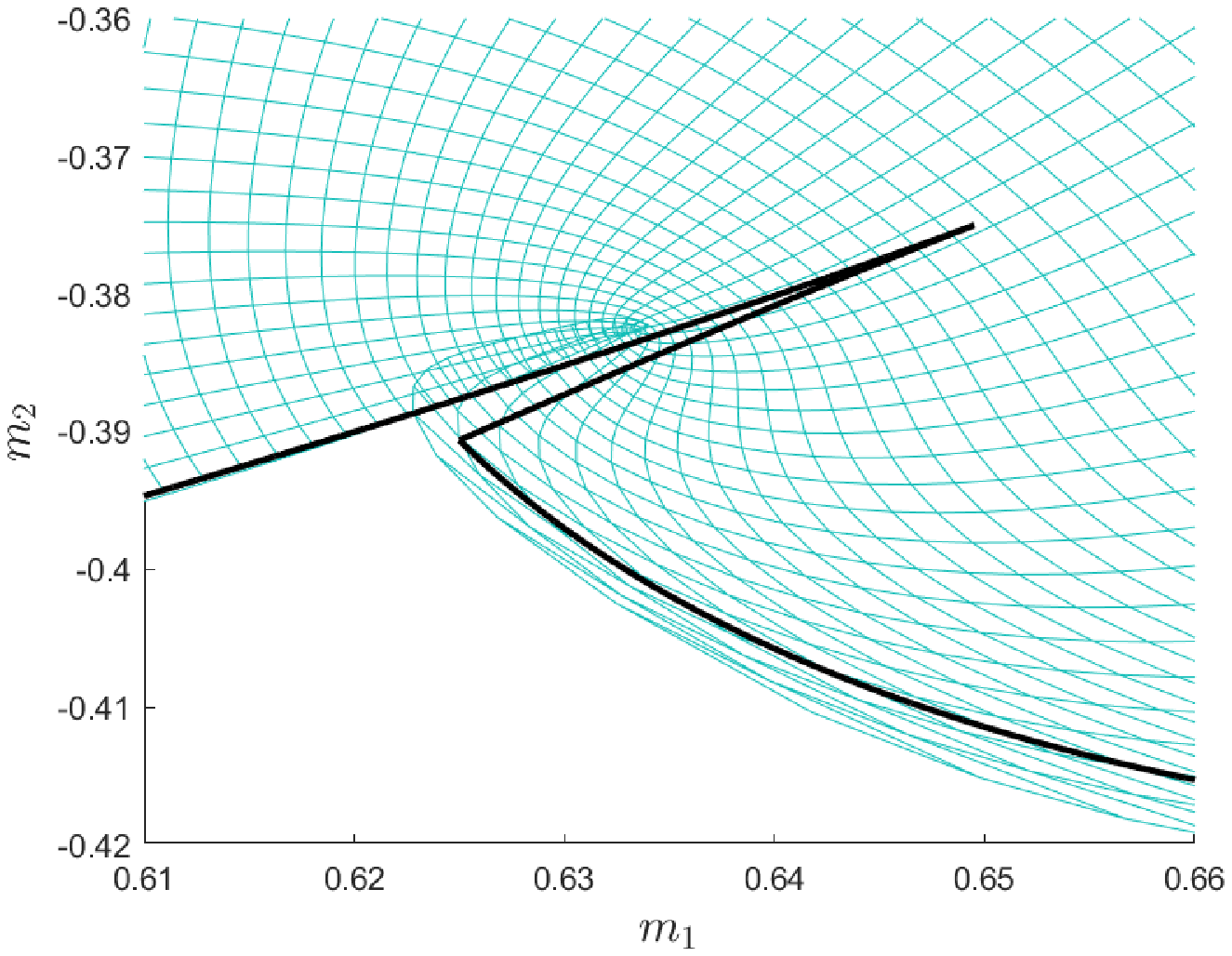"}
     \caption{Zoomed-in results for PM (left) and SPM  (right) on a $161 \times 161$ grid.}
     \label{fig:case_11_PM_SPM_mapping_zoomed}
\end{figure}

\begin{figure}[t!]
    \centering
    \includegraphics[width = 0.44\linewidth]{"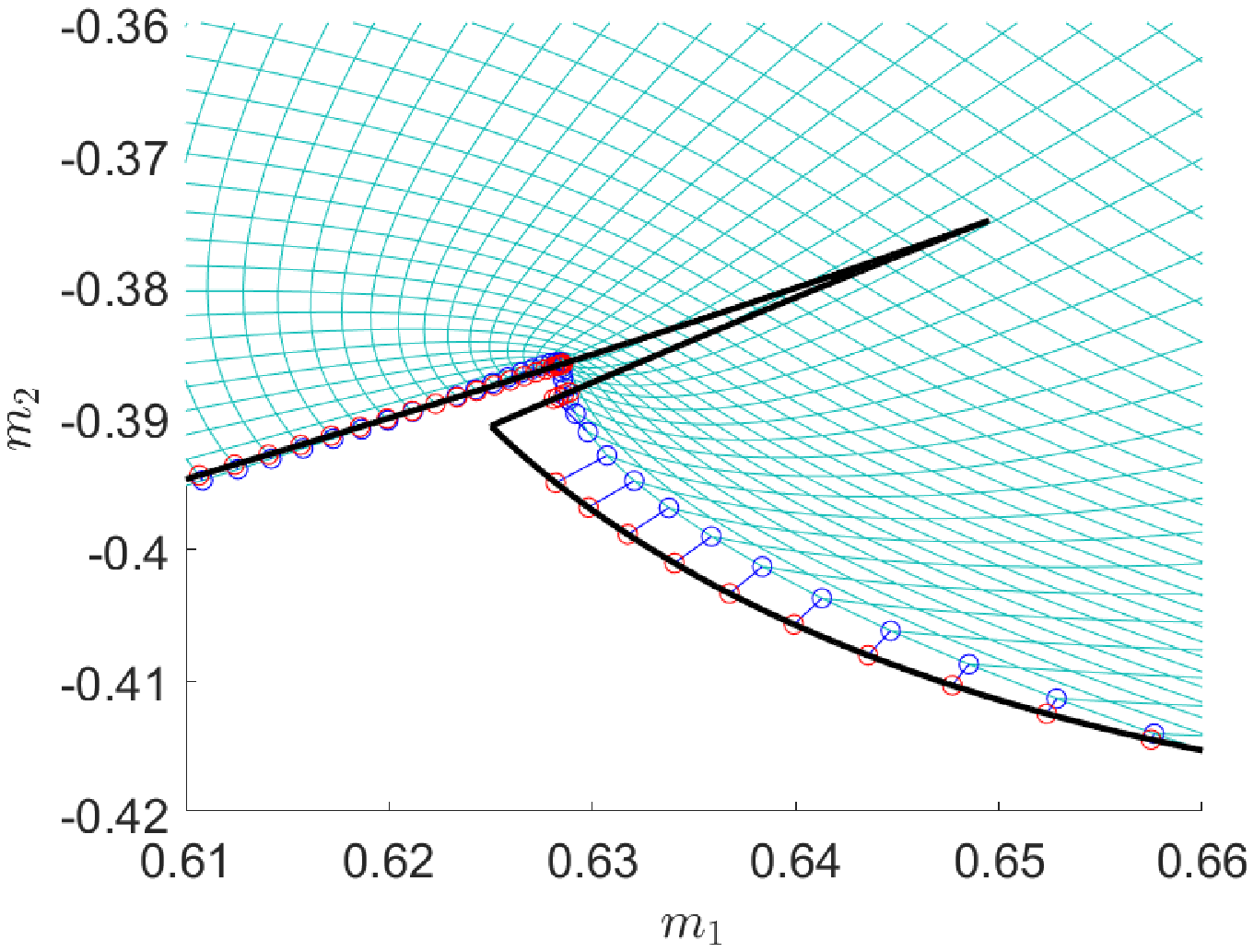"}
    \hspace{-5pt}
    \includegraphics[width = 0.44\linewidth]{"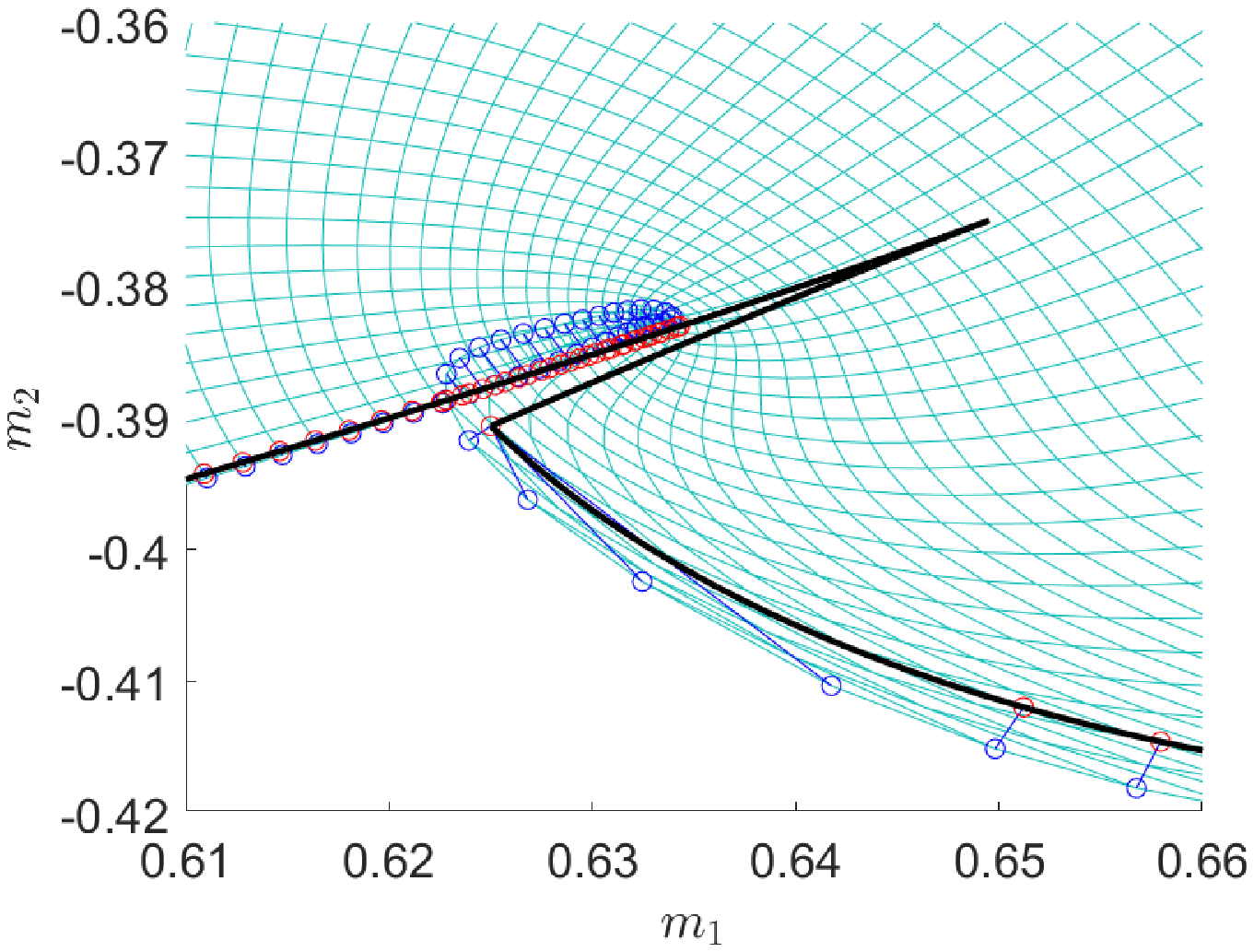"}
    \caption{Projection step after convergence for PM (left) and for SPM (right) on a $161 \times 161$ grid.}
    \label{fig:case_11_PM_SPM_mapping_zoomed_bProj}
\end{figure}

PM and SPM do not yield converging numerical approximations. Figure~\ref{fig:case_11_PM_SPM_mapping_zoomed} shows a zoomed-in version of two numerical solutions for $161 \times 161$ grids. The sharp inward fold seems to be the culprit for the boundary method, as is also seen in Figure~\ref{fig:case_11_PM_SPM_mapping_zoomed_bProj}, which shows the projection of $\mathbf{m}_\elll$ onto $\partial \mathcal{Y}$. The figure clearly shows that the method does not pick $\mathbf{b}_\elll$ deep within the fold and, consequently, the optimization for $\mathbf{m}$ does not produce a mapping with such a sharp fold.

\begin{figure}[!t]
    \centering
    \includegraphics[width = 0.345\linewidth]{"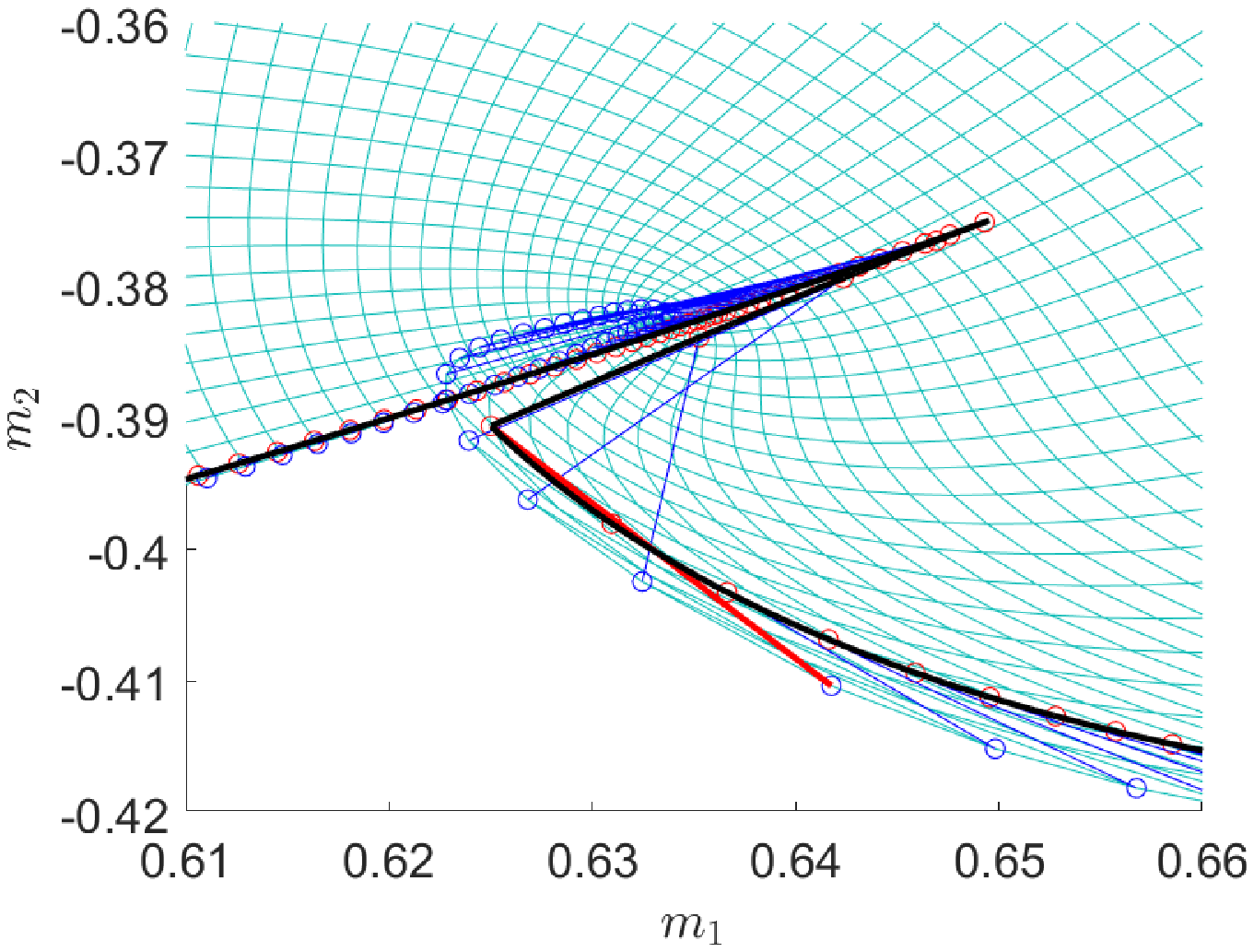"}
	\hspace{-15pt}
    \includegraphics[width = 0.345\linewidth]{"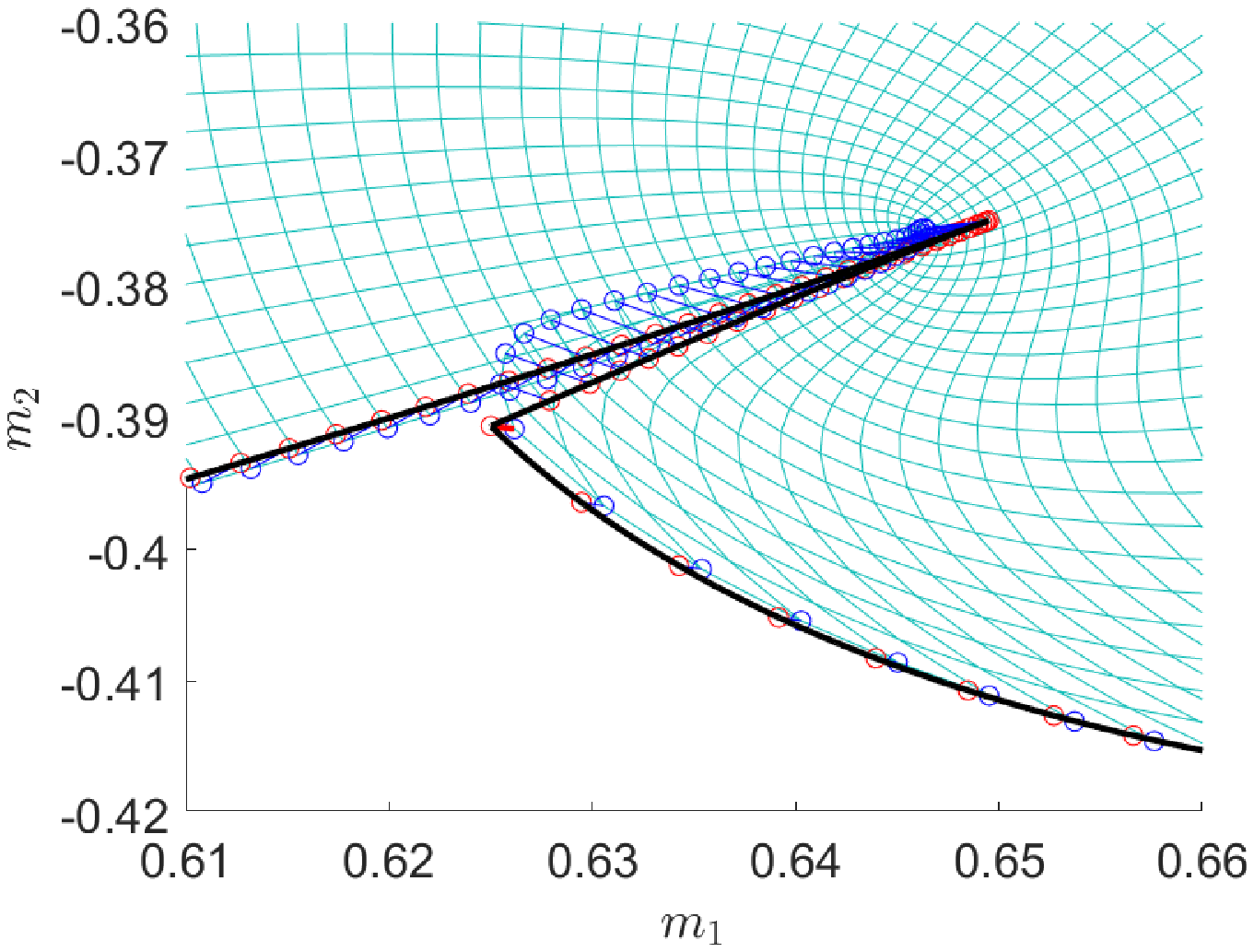"}
	\hspace{-15pt}
    \includegraphics[width = 0.345\linewidth]{"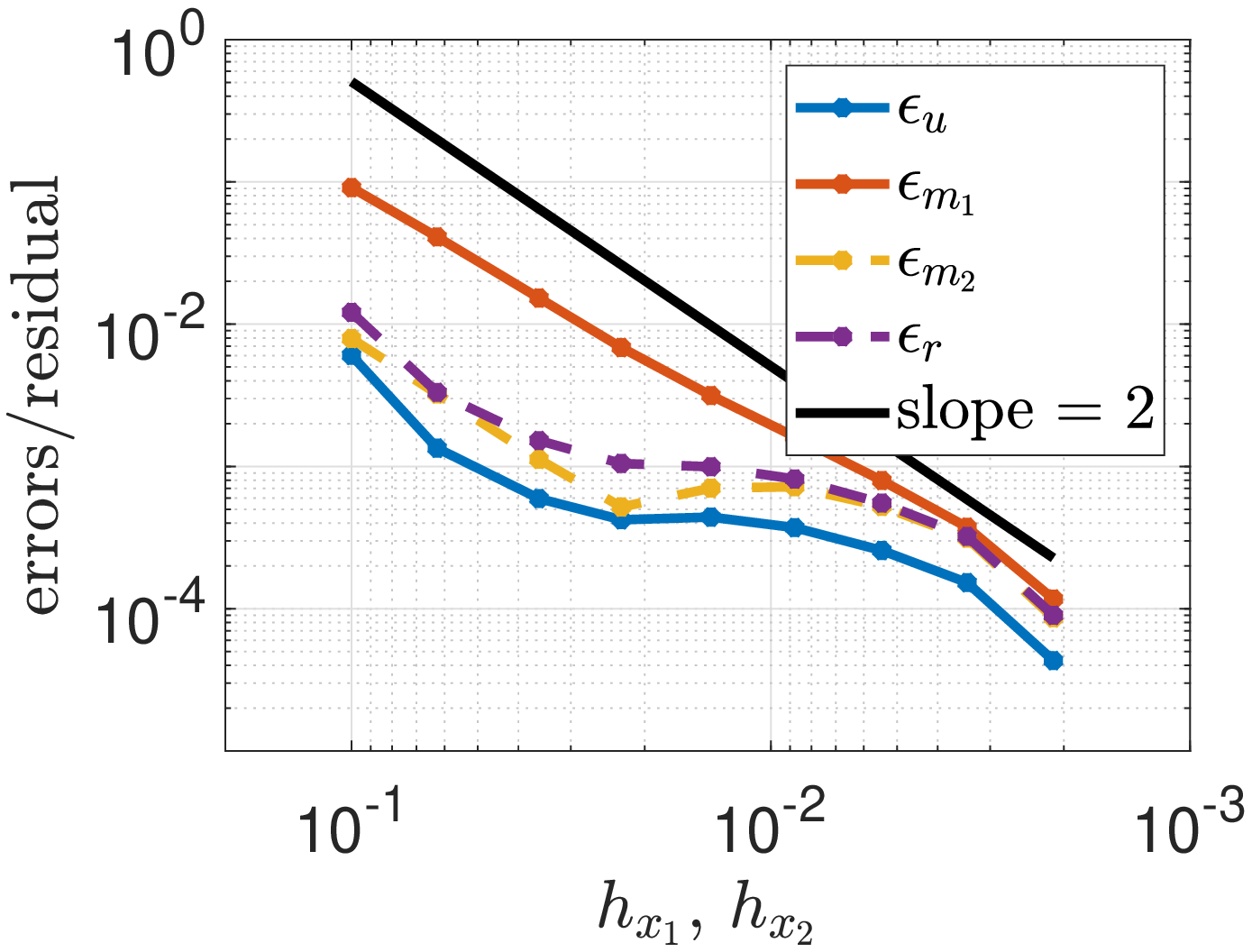"}
    \caption{The first (left) and 50$^{th}$ (middle) iteration of continuation by SALM after SPM has converged on a $161\times 161$ grid and the convergence of the errors and residual for SALM starting from a uniform initial guess (right).}
    \label{fig:case_11_SALM_mapping_zoomed_bProj}
\end{figure}

Because SALM does force points $\mathbf{b}_\elll$ to be located along the whole boundary, naturally points will end up in the fold. This can be seen in Figure~\ref{fig:case_11_SALM_mapping_zoomed_bProj}, where on the left the first iteration of applying SALM to the result of SPM is shown.
The $50^\text{th}$ iteration of SALM is shown in the middle, showing $\mathbf{m}(\partial \mathcal{X})$ being positioned in the fold.
Continuation using SALM yields similar results to using SALM starting from the default initial guess.
SALM shows approximately second-order convergence, as graphed on the right of Figure~\ref{fig:case_11_SALM_mapping_zoomed_bProj}, when using~\eqref{eqn:initial_guess_m} as initial guess. SALM does not show any visual distortions, in contrast to PM and SPM.

For the remaining results we will not discuss PM, as it performs, at best, as good as SPM while being more computationally expensive.

\FloatBarrier

\subsection{Annulus}
\label{sec:Annulus}
\begin{figure}[b!]
   	\centering
    \includegraphics[width = 0.45\linewidth]{"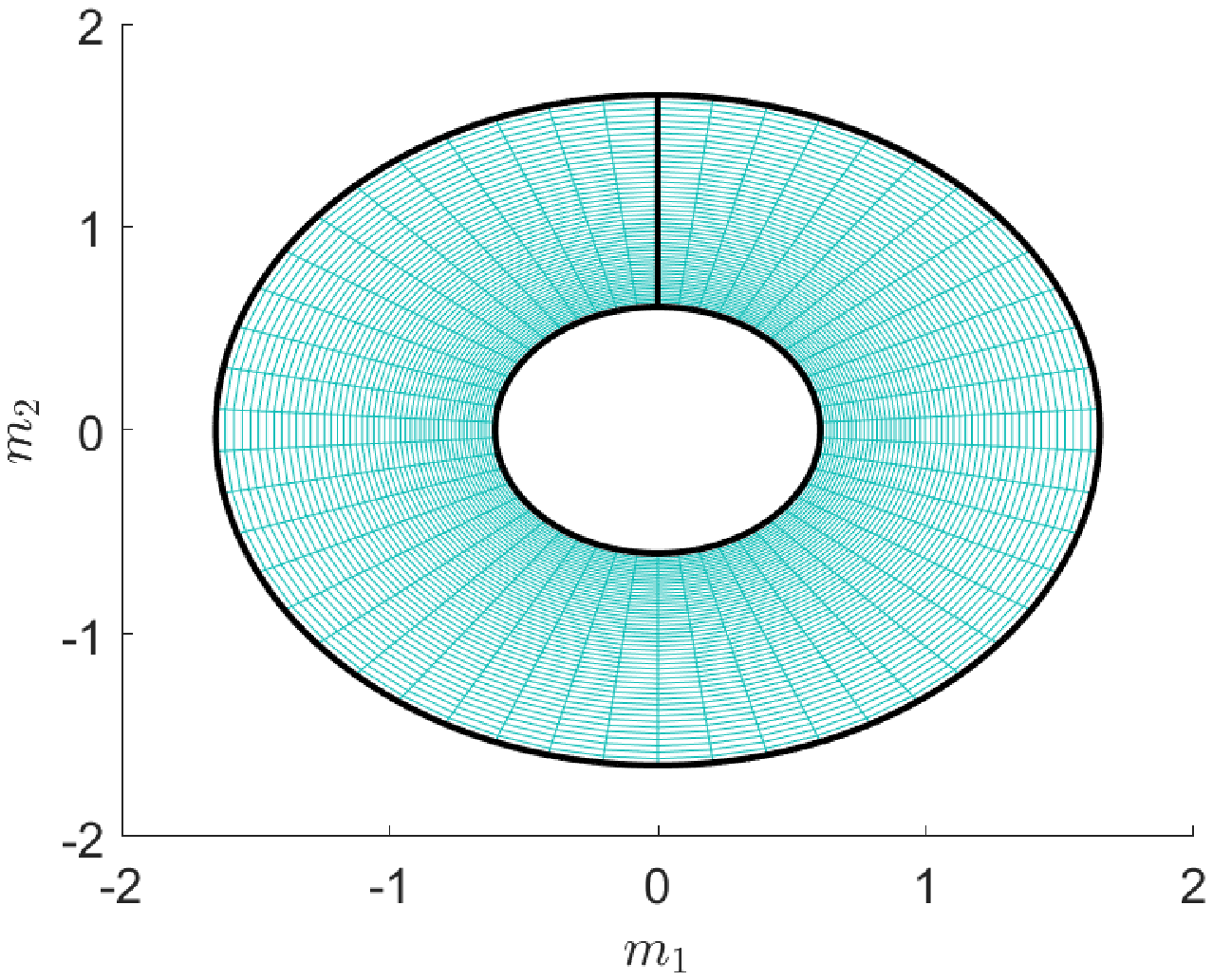"}
    \quad
    \includegraphics[width = 0.45\linewidth]{"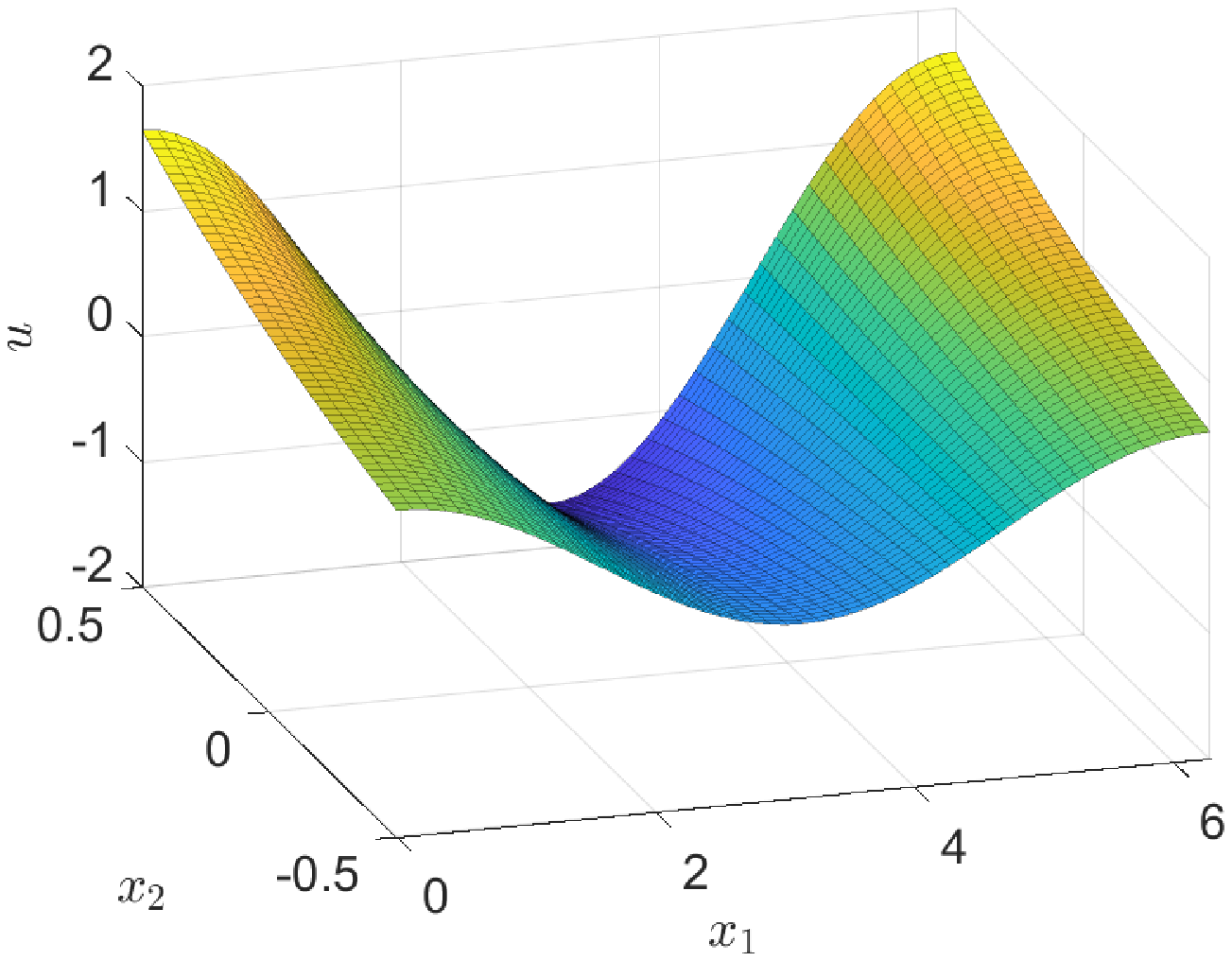"}        
    \caption{The target and the exact mapping on the left, and the solution surface on the right, both shown on a $51\times 51$ grid.}
    \label{fig:Case_9_exact}
\end{figure}
For this example we consider the target given in Figure~\ref{fig:Case_9_exact}, where the central part near $\mathbf{m} = (0,0)$ is not part of $\mathcal{Y}$.
Let $\mathcal{X} = [0, 2 \pi] \times[-1/2, 1/2]$, $\partial \mathcal{Y} = \cup_{k=1}^4 \Gamma_k^\mathcal{Y}$ with
\begin{subequations}
\begin{align}
\Gamma_1^\mathcal{Y}(s) & = (0, \, e^{s-\frac{1}{2}}), \\
\Gamma_2^\mathcal{Y}(s) & = \sqrt{e}(- \sin (2 \pi  s), \, \cos (2 \pi  s))  , \\
\Gamma_3^\mathcal{Y}(s) & = (0, \, e^{\frac{1}{2}-s}) , \\
\Gamma_4^\mathcal{Y}(s) & = \frac{1}{\sqrt{e}}(\sin (2 \pi  s), \, \cos (2 \pi  s)),
\end{align}
\end{subequations}
and $f^2(x_1,x_2) = e^{2x_2}$, such that the exact solution is given by
\begin{align}
    u(x_1,x_2) & = e^{x_2} \cos(x_1),
\end{align}
as shown on the right of Figure~\ref{fig:Case_9_exact}. Observe that $\mathbf{m}\vert_{\partial \mathcal{X}}$ is not bijective, as $\Gamma^\mathcal{Y}_1 = \Gamma^\mathcal{Y}_3$. Nevertheless, we introduce both $\Gamma^\mathcal{Y}_1$ and $\Gamma^\mathcal{Y}_3$ as the orientation, i.e., the parametrization of the segments, matters for SALM. 

\begin{figure}[!t]
    \centering
    \includegraphics[width = 0.45\linewidth]{"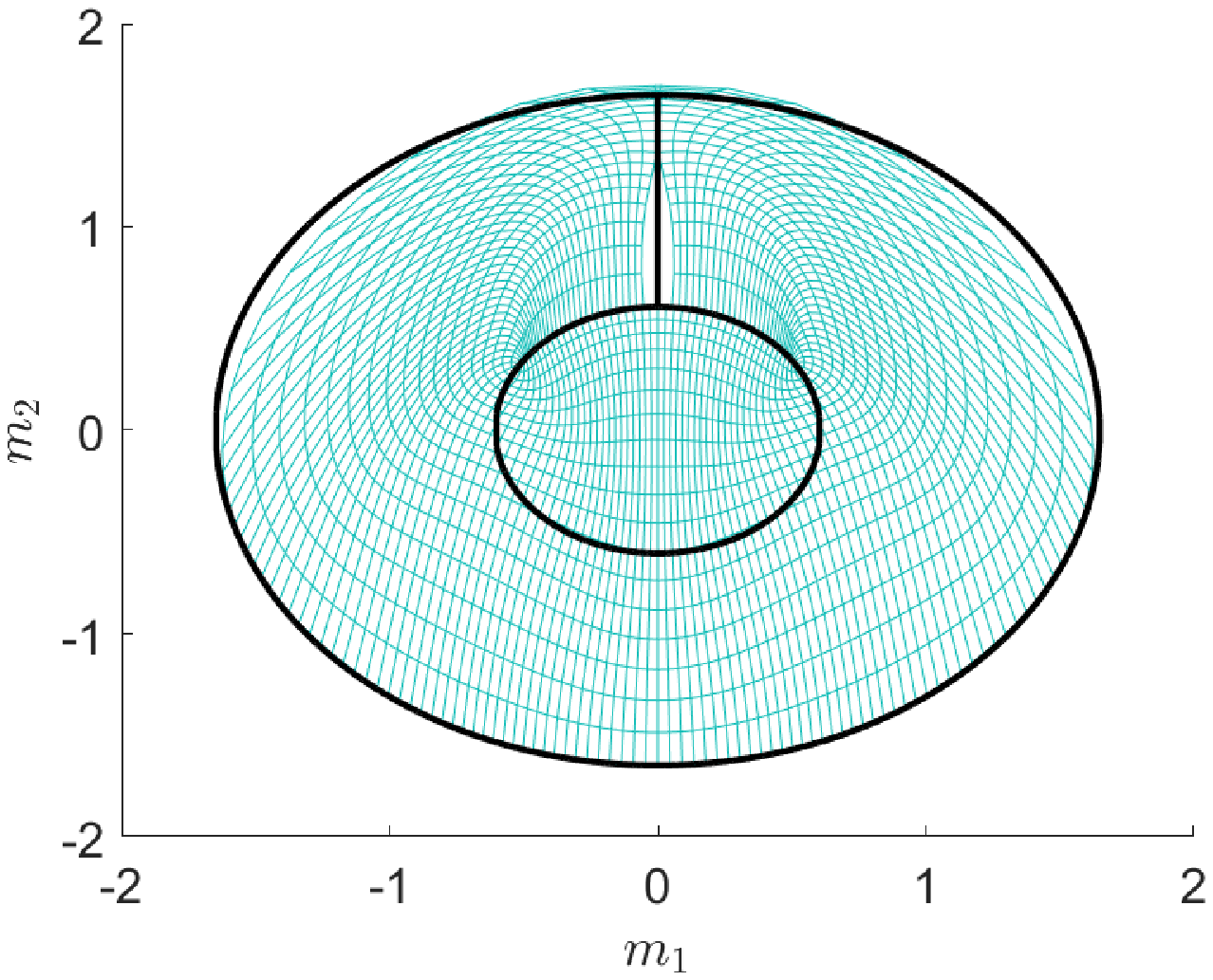"}
	\quad
    \includegraphics[width = 0.45\linewidth]{"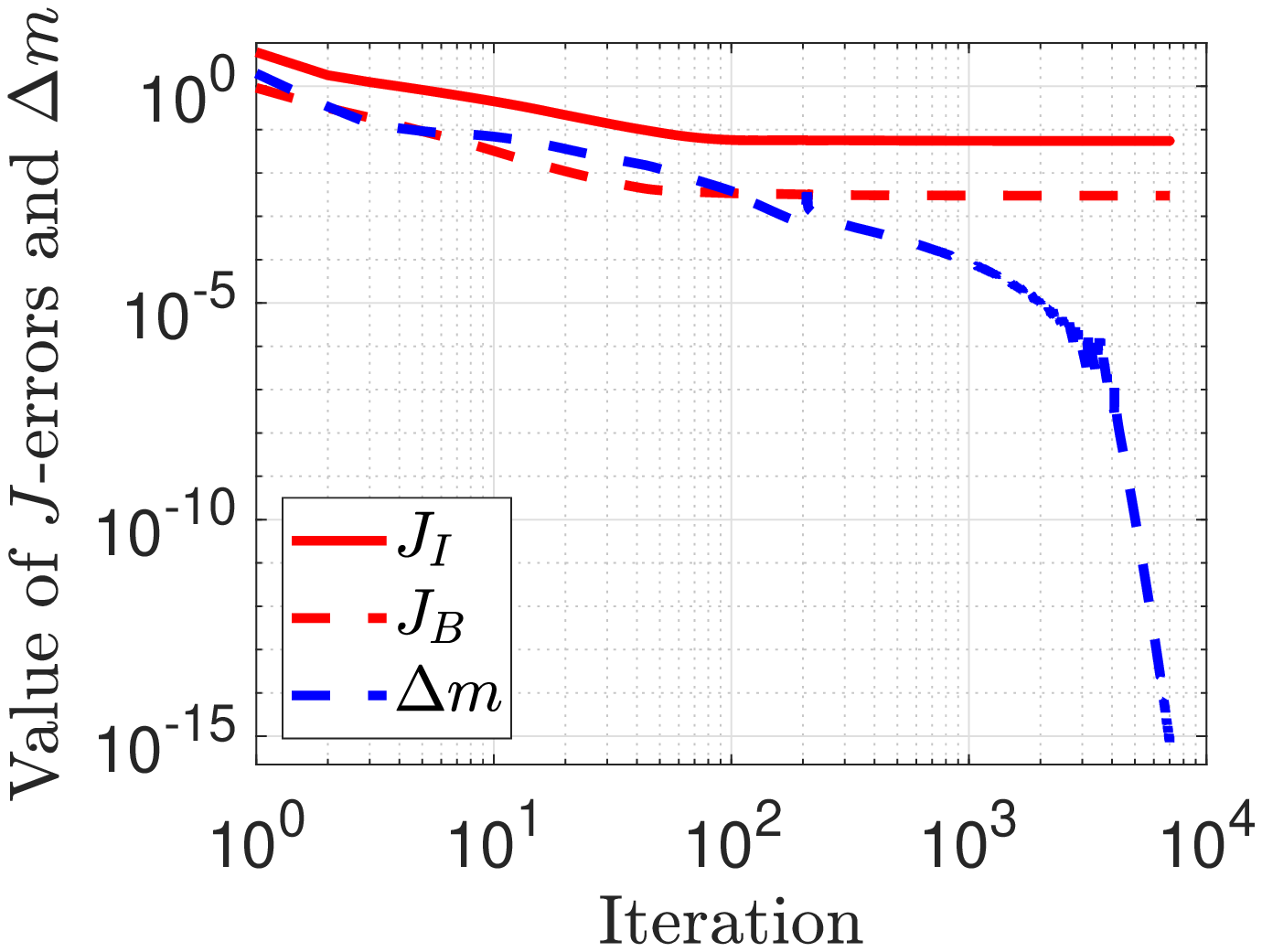"}
    \caption{ The numerical mapping $\mathbf{m}$ after convergence (left) and the history of $J_\textrm{I}$, $J_{\textrm{B}}$ and $\Delta m$ for SPM. }
    \label{fig:case_9_SPM}
\end{figure}

\begin{figure}[!t]
	\centering
	\includegraphics[width = 0.45\linewidth]{"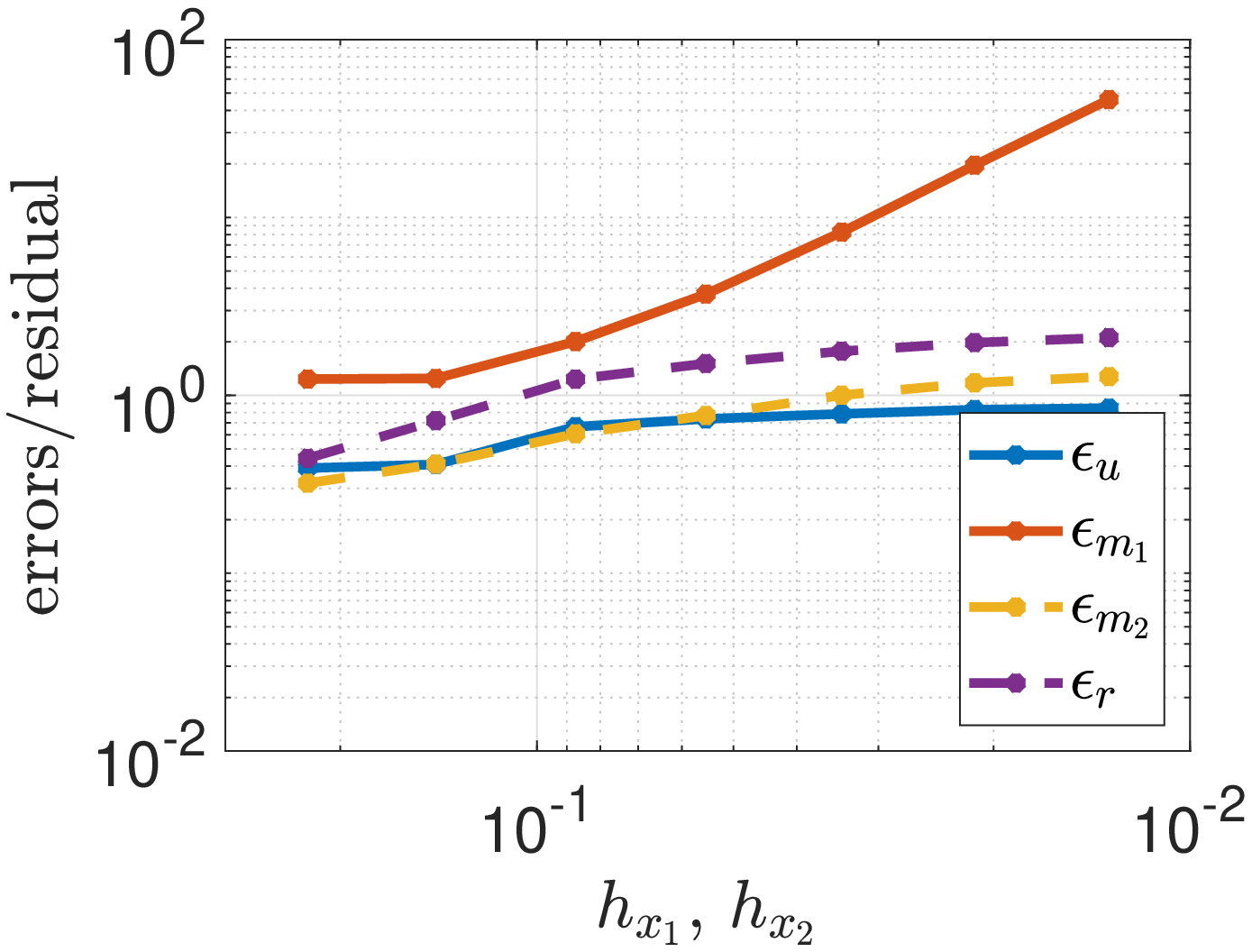"}
	\hspace{5pt}
	\includegraphics[width = 0.45\linewidth]{"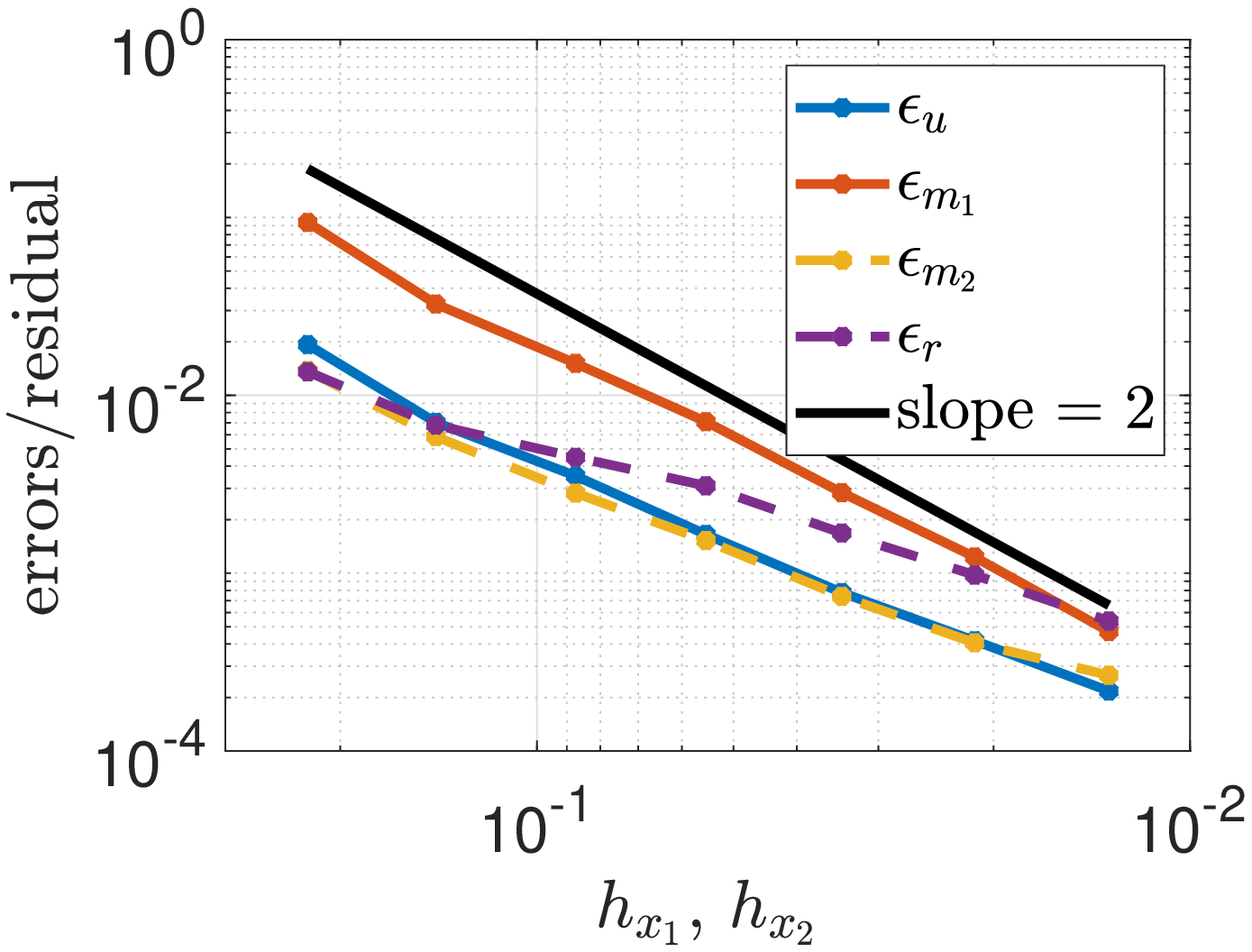"}
	\caption{ Residual for SPM (left) and SALM (right). No convergence for SPM, and second-order convergence for SALM is observed. }
	\label{fig:case_9_JIJB_SPM_SALM}
\end{figure}

Figure~\ref{fig:case_9_SPM} shows results for SPM. On the left the mapping after the algorithm has converged for a grid with $N_{x_1} = 115$ and $N_{x_2} = 19$. The grid parameters are chosen such that $h_{x_1} \approx h_{x_2}$ as $N_{x_1} / N_{x_2} \approx (x_1^M - x_1^m)/(x_2^M - x_2^m) = 2 \pi$. 
Although the figure on the right clearly shows $J_\textrm{I}$ and $J_{\textrm{B}}$ have converged, and that $\Delta m$ reached computer precision, the algorithm does not yield a correct solution, as it does not satisfy the transport boundary condition because there are points $\mathbf{m}_{ij}$ which lie outside $\mathcal{Y}$, nor does it solve the hyperbolic \MAe{} as is shown by the residual $\epsilon_r$ in Figure~\ref{fig:case_9_JIJB_SPM_SALM} on the left. 

\begin{figure}[!t]
	\centering
	\includegraphics[width = 0.33\linewidth]{"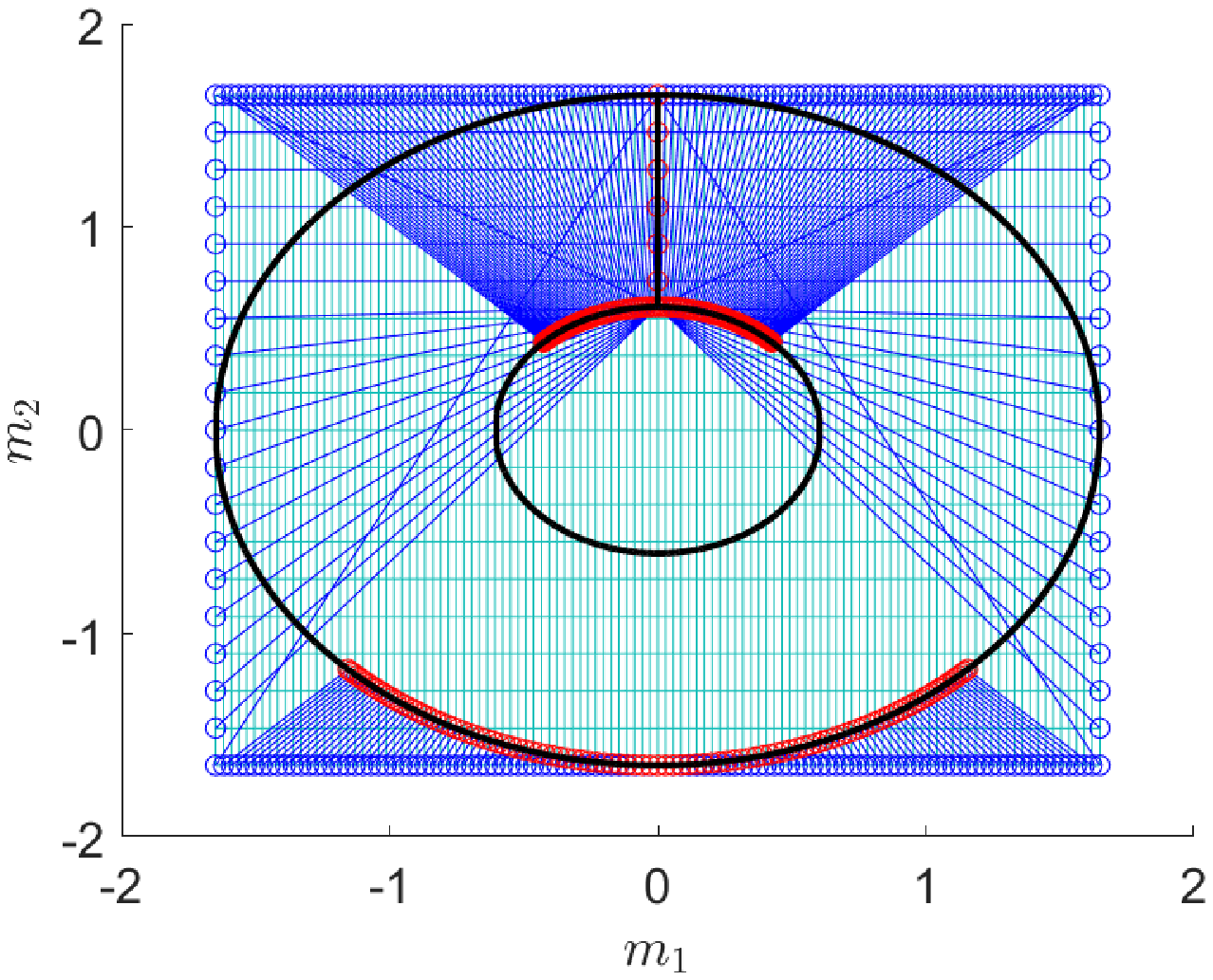"}
	\hspace{-5pt}
	\includegraphics[width = 0.33\linewidth]{"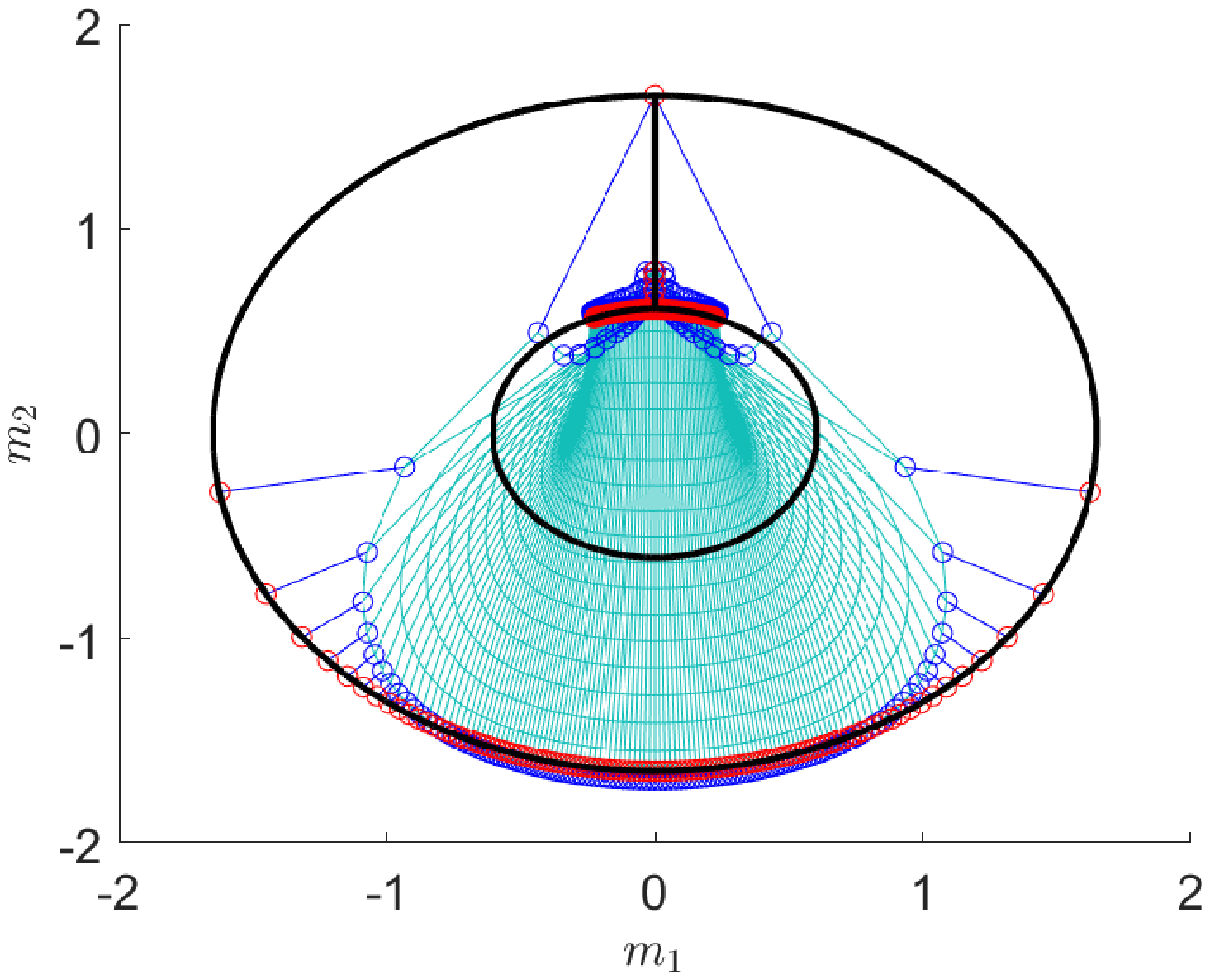"}
	\hspace{-5pt}
	\includegraphics[width = 0.33\linewidth]{"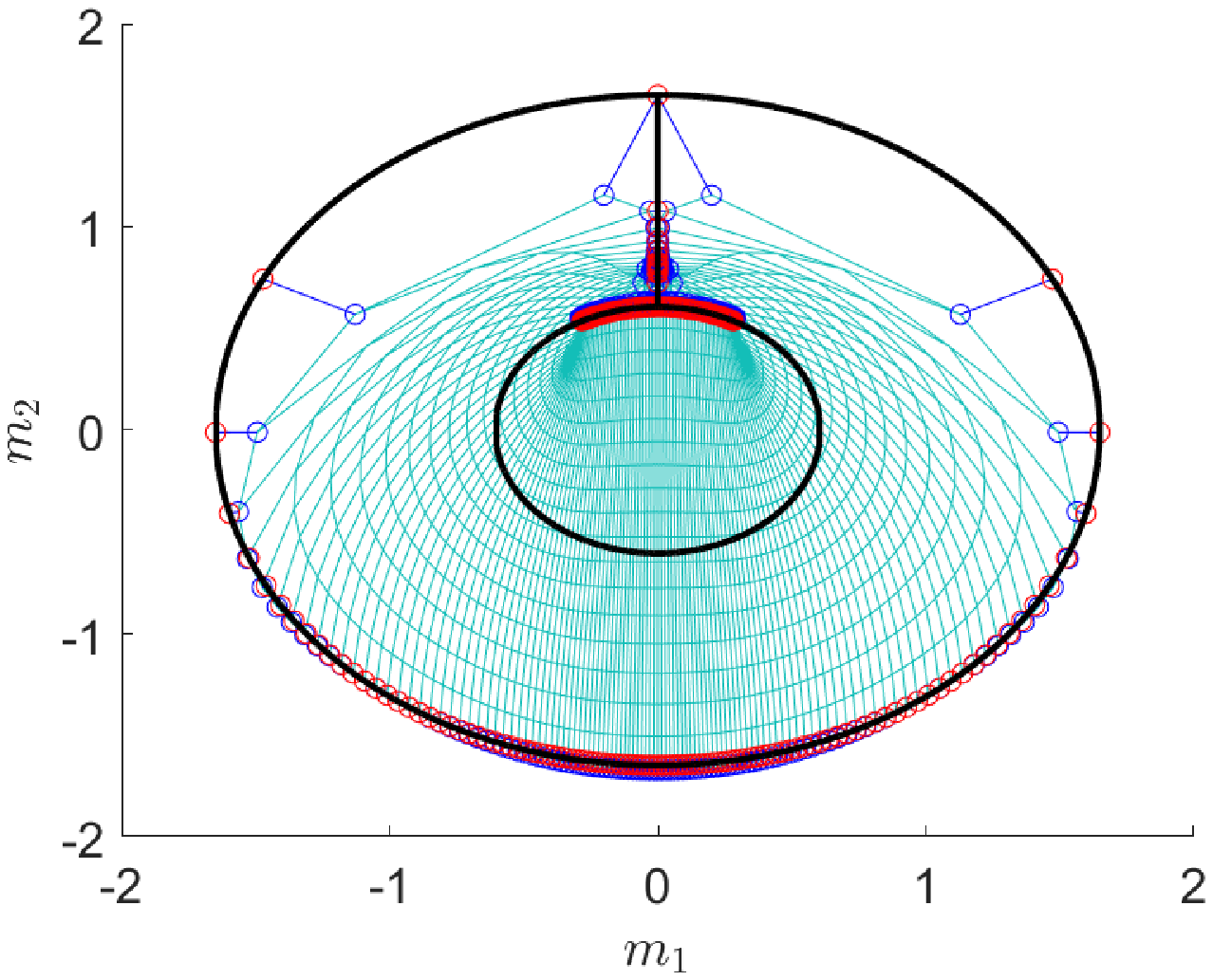"}
	\caption{From left to right, the first, third and tenth iteration for SPM. }
	\label{fig:Case_9_bndM32_firstIts}
\end{figure}

The reason why SPM does not produce accurate solutions is easiest demonstrated by visualizing a few iterations. To this end, consider the boundary routine for the first, third and tenth iteration as shown in Figure~\ref{fig:Case_9_bndM32_firstIts}. We focus on one segment of the mapping of the boundary, i.e, $\mathbf{m}_{ij}$ with $\mathbf{x}_{ij} \in \Gamma^\mathcal{X}_4 = [0, 2 \pi] \times \{-\tfrac{1}{2}\}$, which corresponds to $m_2(0, -\tfrac{1}{2}) = \exp(-1/2) \approx 0.6$ for the initial guess in Figure~\ref{fig:Case_9_bndM32_firstIts}. For the exact solution, $\Gamma^\mathcal{X}_4$ needs to be mapped to the entire inner circle of the target, i.e., $\Gamma^\mathcal{Y}_4$. As shown for the first iteration, $\Gamma^\mathcal{X}_4$ is mapped to only part of $\Gamma^\mathcal{Y}_4$, viz., the accompanying $\mathbf{b}_\elll$ lies on the northern part of $\Gamma^\mathcal{Y}_4$ (the inner circle). In subsequent iterations, shown in the middle and on the right in Figure~\ref{fig:Case_9_bndM32_firstIts}, $\Gamma^\mathcal{X}_4$ will again not be mapped to the whole of $\Gamma^\mathcal{Y}_4$, as the distance to the northern part of the inner circle remains minimal. This process continues indefinitely.

\begin{figure}[t!]
	\centering
	\includegraphics[width = 0.345\linewidth]{"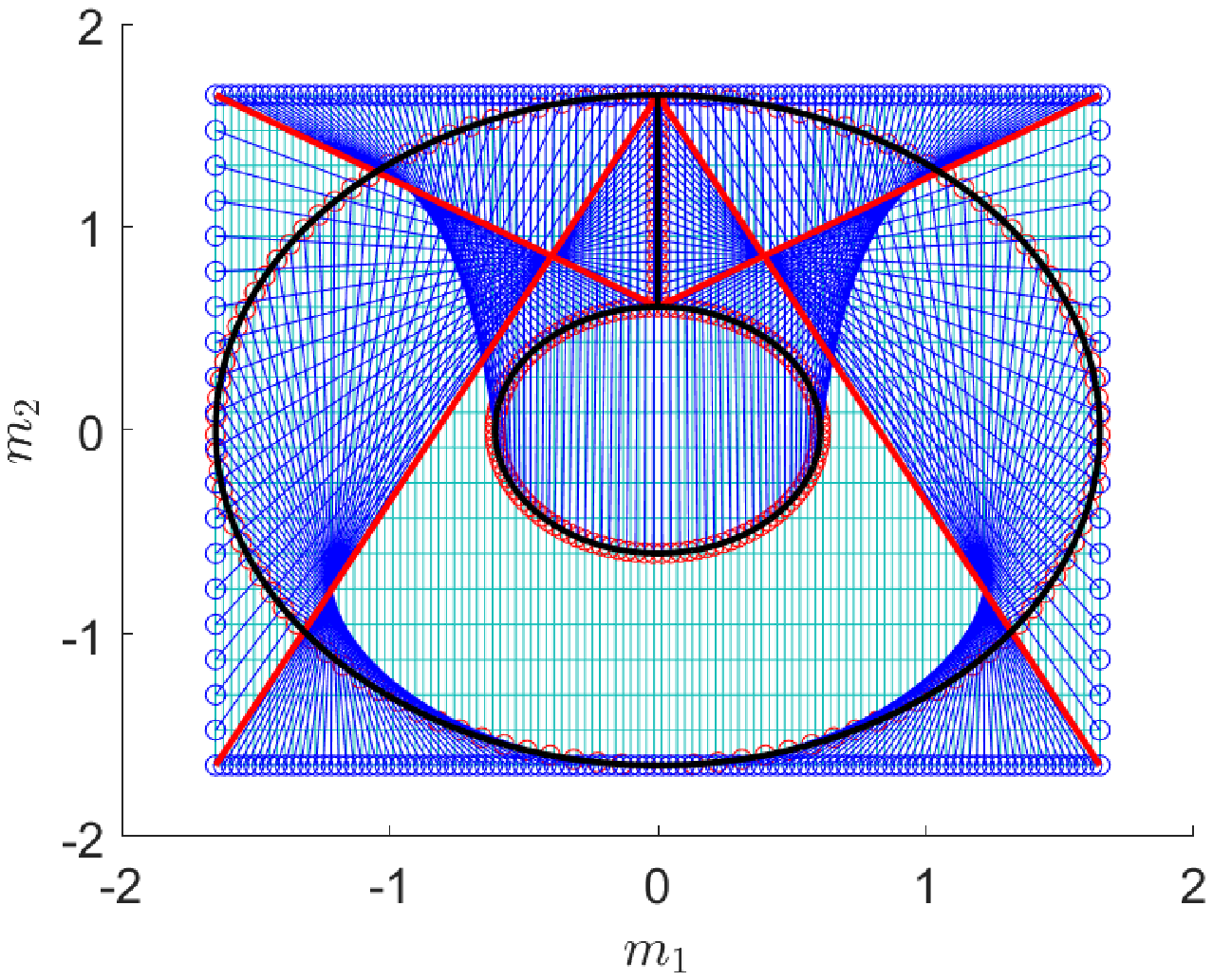"}
	\hspace{-15pt}
	\includegraphics[width = 0.345\linewidth]{"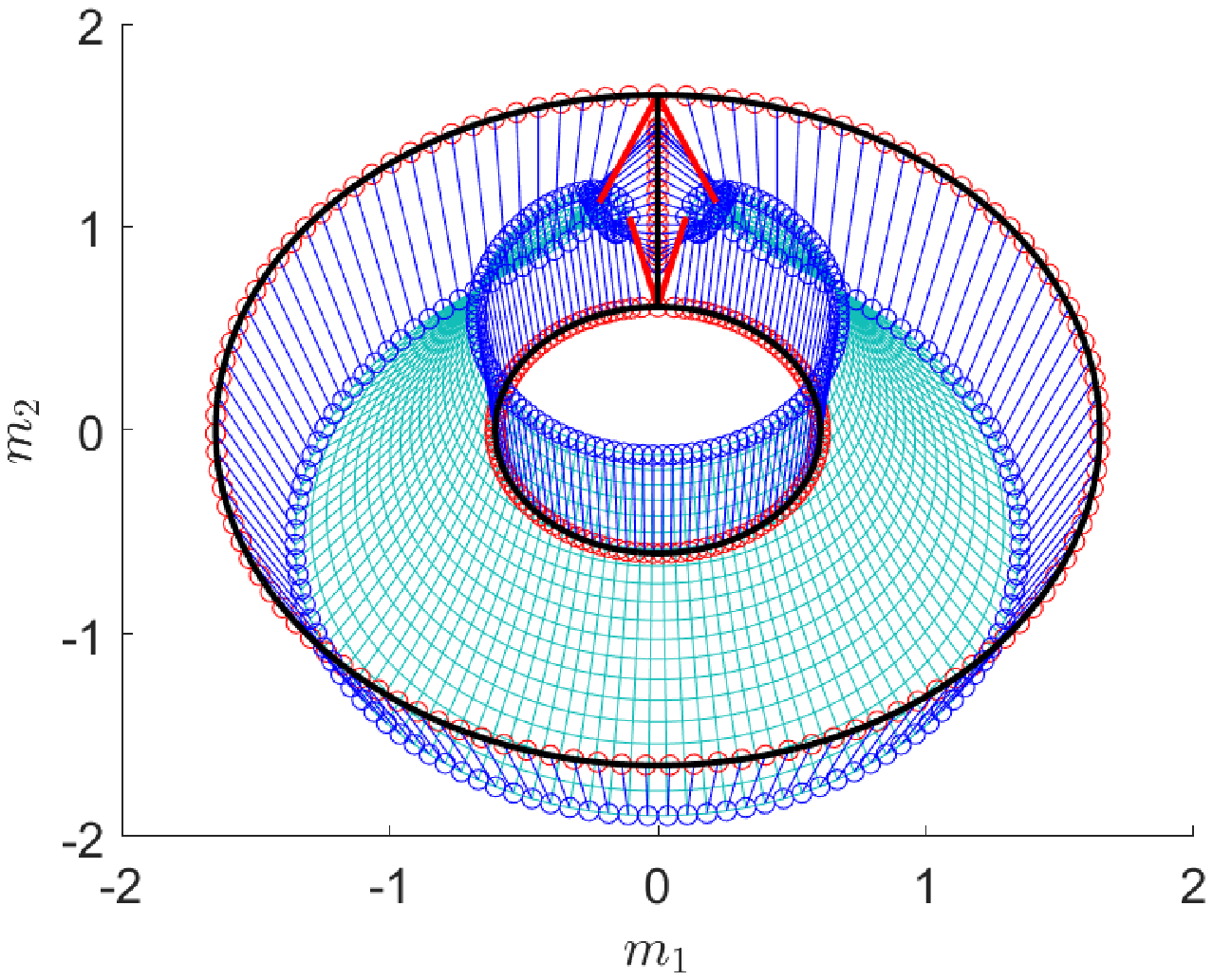"}
	\hspace{-15pt}
	\includegraphics[width = 0.345\linewidth]{"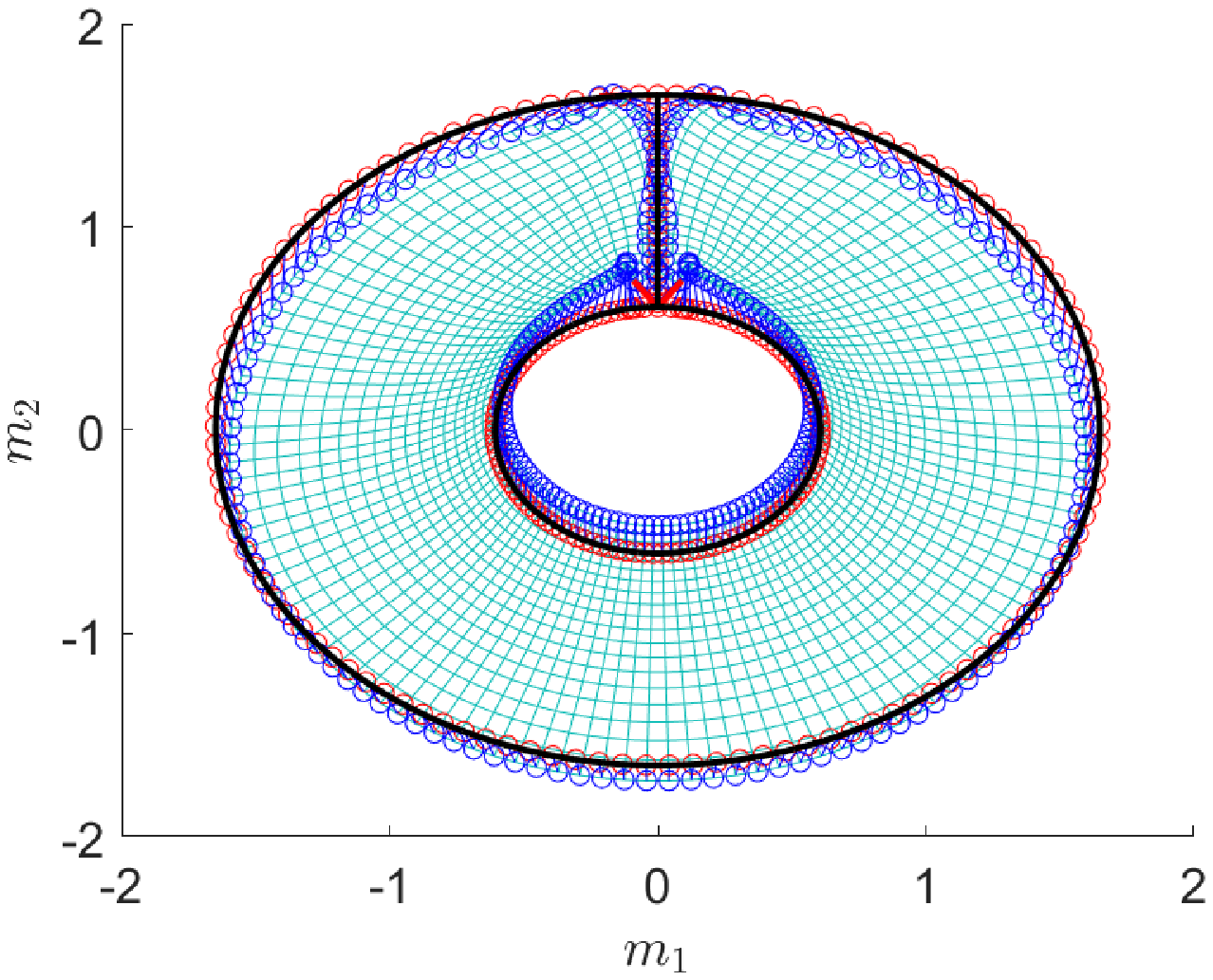"}
	\caption{From left to right, the first, second and third iteration for SALM. }
	\label{fig:Case_9_bndM10_firstIts}
\end{figure}

For SALM such accumulation of $\mathbf{b}_\elll$ does not occur, as by construction, $\mathbf{b}_\elll$ is distributed over the boundary segments. The results of the first, second and third iteration of the $\mathbf{b}$-minimization are shown in Figure~\ref{fig:Case_9_bndM10_firstIts}. Clearly, SALM does not suffer from the same flaws as SPM. As such, the convergence is expected to behave as for the other examples, which is confirmed by the results shown in Figure~\ref{fig:case_9_JIJB_SPM_SALM} on the right.


\subsection{Gradient dependent problem}
\begin{figure}[!t]
   	\centering
   	\includegraphics[width = 0.45\linewidth]{"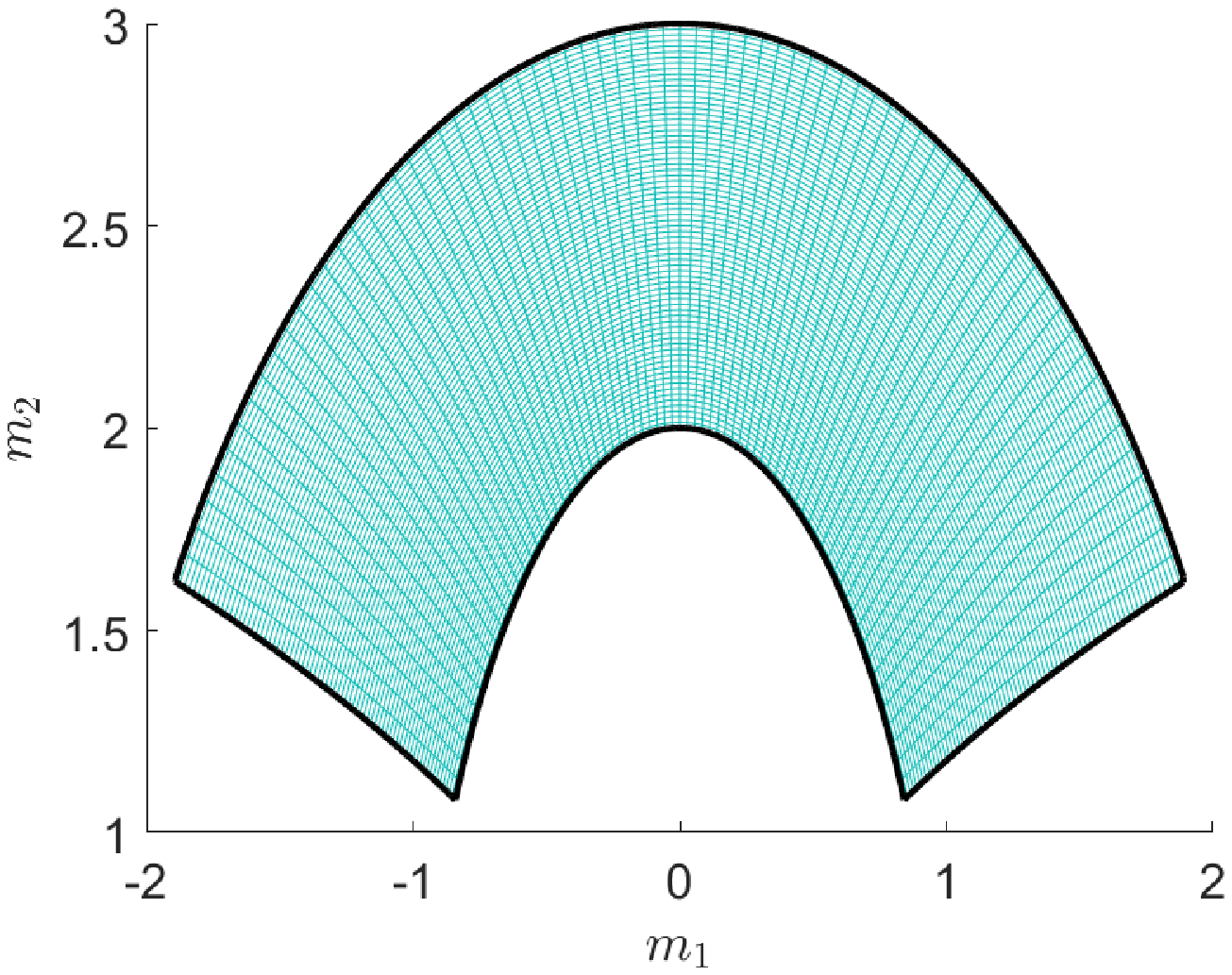"}
   	\quad
   	\includegraphics[width = 0.45\linewidth]{"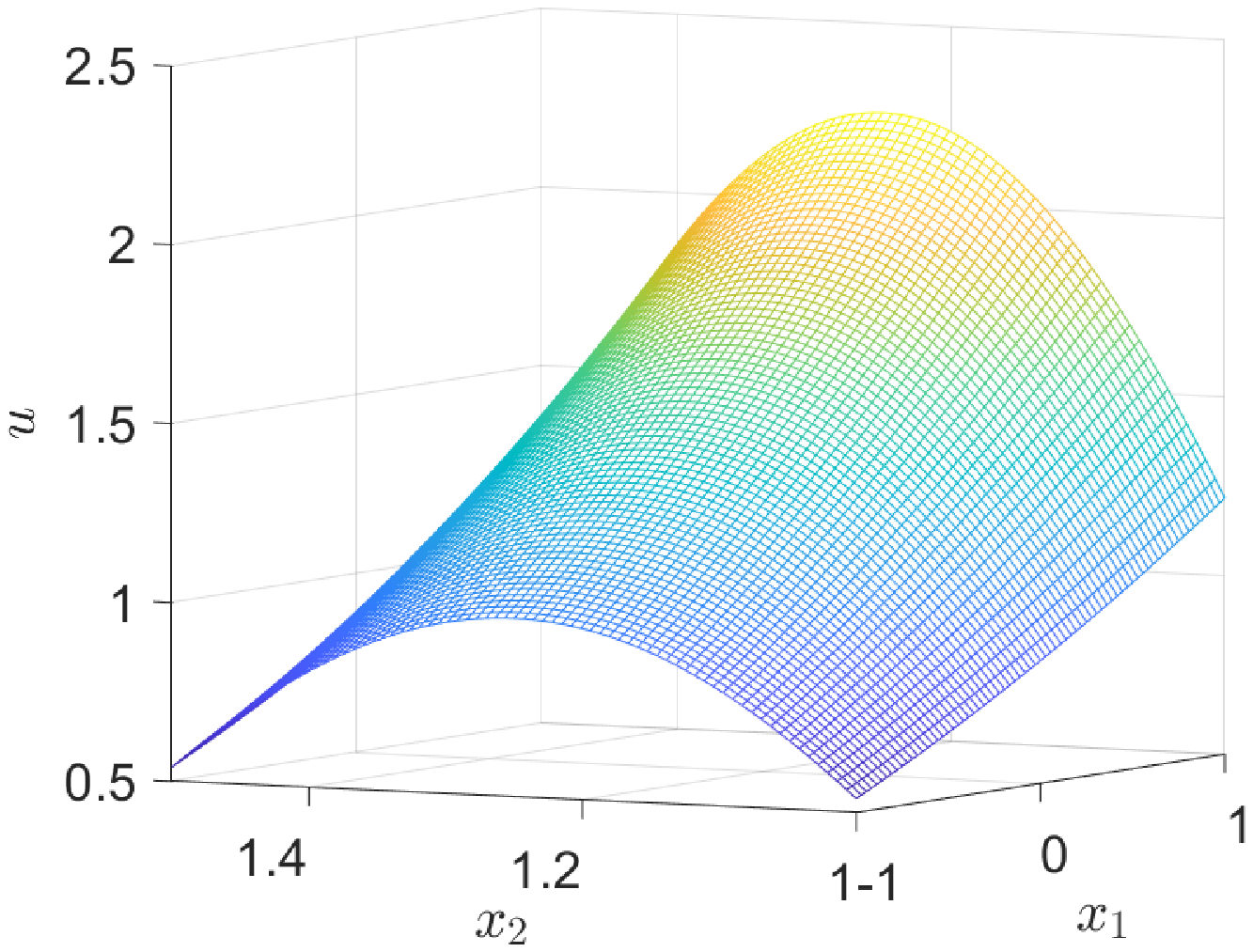"}
   	\caption{ The target and the exact mapping on the left, and the solution surface on the right, both shown on a $51\times 51$ grid. }
   	\label{fig:case_502_solutions}
\end{figure}

Lastly we consider an example with $f$ dependent on the gradient of the solution, i.e., $f = f(x_1, x_2, \nabla u)$, viz.
\begin{align}
   f^2(\mathbf{x}, \mathbf{m}) = 3 x_2^2 - m_1 \sin(x_1) - \frac{1}{4} m_2^2.
\end{align}
We consider the domain $\mathcal{X} = [-1, 1] \times[1, 3/2]$ and $\partial \mathcal{Y} = \cup_{k=1}^4 \Gamma_k^\mathcal{Y}$ with
\begin{subequations}
    \begin{align}
    \Gamma_1^\mathcal{Y}(s) & = (\tfrac{\sin(1)}{4} s^2 + \sin(1) s + \sin(1), \, \cos(1) s + 2 \cos(1))  , \\
    \Gamma_2^\mathcal{Y}(s) & = (-\tfrac{9}{4} \sin(2s-1), \, 3 \cos(2s-1)) , \\
    \Gamma_3^\mathcal{Y}(s) & = (-\tfrac{\sin(1)}{4} s^2 + \tfrac{3}{2} s - \tfrac{9 \sin(1)}{4}, \, -\cos(1) s + 3 \cos(1)) , \\
    \Gamma_4^\mathcal{Y}(s) & = (\sin(2s -1), \, 2 \cos(2s-1)).
    \end{align}
\end{subequations}    
The exact solution is given by
\begin{align}
u(x_1,x_2) & = x_2^2\cos(x_1), 
\end{align}
and is, together with the mapping and target domain, shown in Figure~\ref{fig:case_502_solutions}.

\begin{figure}[t!]
    \centering
    \includegraphics[width = 0.45\linewidth]{"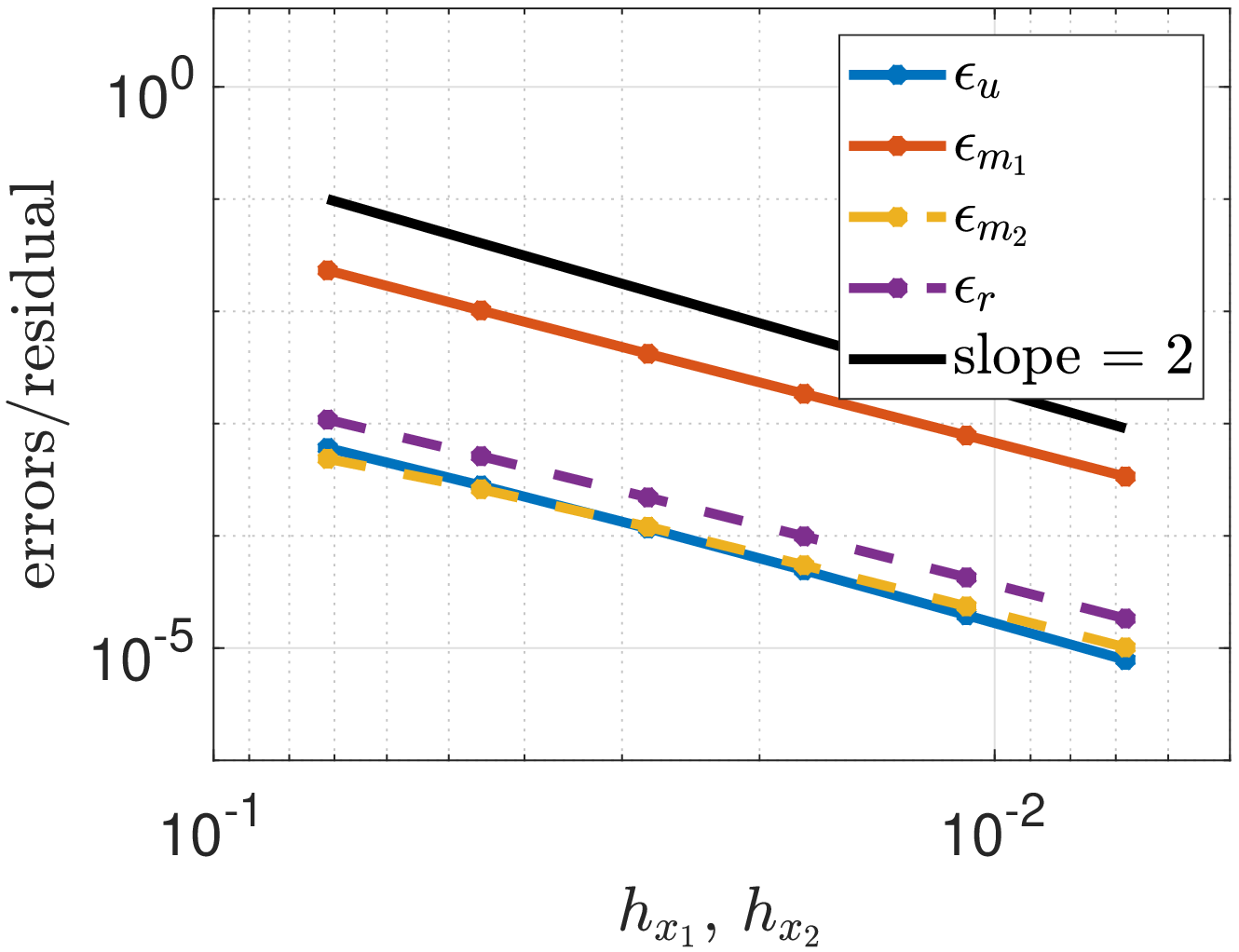"}
    \quad
    \includegraphics[width = 0.45\linewidth]{"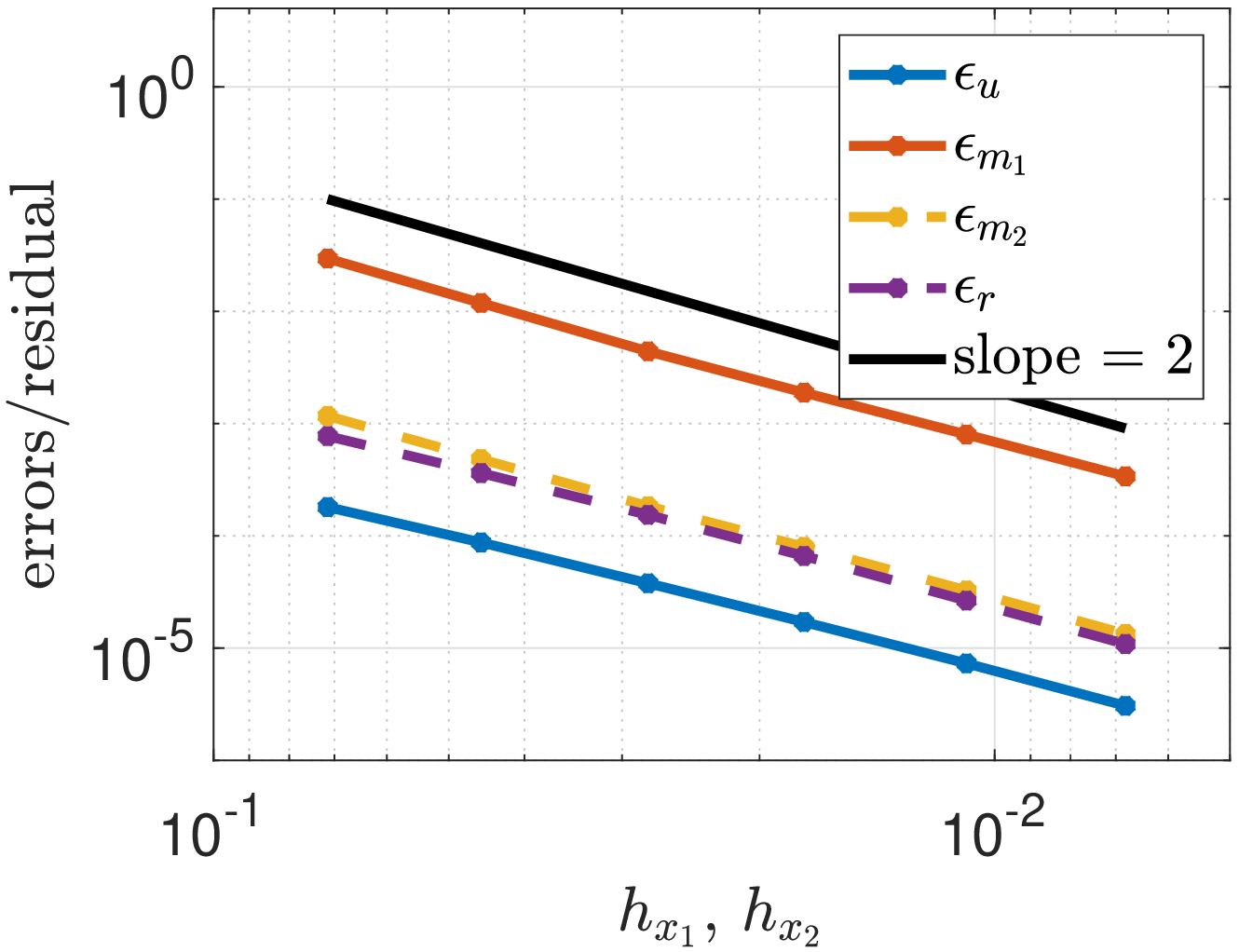"}
    \caption{ Convergence of SPM (left) and SALM (right). }
    \label{fig:case_502_SPM_SALM}
\end{figure}

By construction of the algorithm, little effort is required for $f$ to be dependent on the mapping $\mathbf{m}$. The difference being that during the $n^\text{th}$ iteration, $f^2(\mathbf{x}_{ij}, \mathbf{m}^n_{ij})$ has to be evaluated instead of $f^2(\mathbf{x}_{ij})$ in the optimization of $\mathbf{P}$.
The results for SPM and SALM with $N_{x_1} = N_{x_2}$ are given in Figure~\ref{fig:case_502_SPM_SALM}, showing second-order convergence for both methods. In this case grid shock correction is needed for SPM  to ensure proper convergence.

\FloatBarrier

\section{Conclusion}
\label{sec:conclusion}
We have introduced a least-squares solver for the hyperbolic \MAe{} with transport boundary condition. The algorithm, originally introduced by Prins et al.~\cite{Prins_2015} for the elliptic \MAe{}, has been improved to encompass a more complete description of the roots for the $P$-optimization. Furthermore, we introduced two new boundary methods. All three boundary methods, if convergent, show second-order convergence of the residual and the global discretization errors as function of the mesh size, and also second-order convergence as function of the number of boundary points. Of the three boundary methods, the \textit{segmented arc length method} is both the only method to converge for all examples and is most computationally efficient, both in terms of computation time per iteration, as in total number of iterations required.

As far as the authors are aware, the least-squares method paired with the \textit{segmented arc length method} for the boundary, is the first method to solve the hyperbolic \MAe{} with transport boundary conditions. 

\bibliographystyle{unsrt}
\bibliography{refs.bib}

\end{document}